\documentclass[11pt]{amsart}

\usepackage{amsmath,amsthm}
\usepackage{amssymb}
\usepackage{epsfig,tikz}
\usepackage{yhmath}
\usepackage{graphicx}
\usepackage{subcaption} 
\usepackage{geometry}
\usepackage{setspace}
\usepackage{enumerate}
\usepackage{tikz}
\usetikzlibrary{calc}
\usepackage{tikz-cd}
\usepackage{xcolor}
\usepackage{stmaryrd}
\usepackage{etoolbox} 
\usepackage[T1]{fontenc}
\usepackage{txfonts}

\usepackage[backref=page]{hyperref}
\hypersetup{
  colorlinks   = true,          
  urlcolor     = blue,          
  linkcolor    = purple,          
  citecolor   = violet             
}

\newtheorem*{maintheoremA*}{Theorem A}
\newtheorem*{maintheoremB*}{Theorem B}
\newtheorem*{mainpropositionA*}{Corollary A}
\newtheorem*{mainpropositionB*}{Corollary B}

\theoremstyle{definition}
\newtheorem{definition}{Definition}[section]
\newtheorem{lemma}[definition]{Lemma}
\newtheorem{proposition}[definition]{Proposition}
\newtheorem{theorem}[definition]{Theorem}

\newtheorem{corollary}[definition]{Corollary}
\newtheorem{example}[definition]{Example}

\newtheorem{notation}[definition]{Notation}

\newtheorem{remark}[definition]{Remark}

\let\tilde\widetilde

\newcommand{\nn}{\mathbb{N}}
\newcommand{\cc}{\mathbb{C}}
\newcommand{\qq}{\mathbb{Q}}

\newcommand{\zz}{\mathbb{Z}}
\newcommand{\Z}{\mathbb{Z}}
\newcommand{\pp}{\mathbb{P}}

\newcommand{\s}{\subset}

\newcommand{\gb}{\beta}

\newcommand{\cyc}[1]{\ensuremath{\langle #1\rangle}}

\usepackage{mathrsfs}

\def\Res{{\rm Res}}
\def\wt{{\rm wt}}

\def\ep{\epsilon}

\newcommand{\la}{\lambda}
\def\A{{\mathbb A}}
\def\C{{\mathbb C}}
\def\P{{\mathbb P}}

\def\V{{\mathbb V}}

\def\Z{{\mathbb Z}}

\def\Y{{\mathcal Y}}

\def\1{{\bf 1}}
\def\tr{{\rm tr}}
\def \End{{\rm End}}

\def \al{\alpha}
\def \b{\beta}

\def \h{\mathfrak{h}}

\def \ra {\rightarrow}
\def \g{\mathfrak{g}}

\def\qdim{{\rm qdim}}

\def \ssq{\subseteq}
\def \C {\mathbb{C}}
\def\ch{{\rm ch}}
\def\cw{\mathrm{cw}}
\def\Id{\mathrm{Id}}
\def\op{\oplus}

\def\vac{\mathbf{1}}
\def\spn{\mathrm{span}}
\def \vac{\mathbf{1}}
\def\rank{\mathrm{rank}}
\def\bs{\backslash}
\def\om{\omega}

\def\o{\otimes}

\def\<{\langle}
\def\>{\rangle}

\def\Om{\Omega}

\def\ds{\dots}

\def\L{\mathcal{L}}

\def\pp{\mathfrak{p}}
\def\qq{\mathfrak{q}}

\def\V{\mathbb{V}}

\def \b{\beta}

\def \h{\mathfrak{h}}

\def \ra {\rightarrow}
\def \g{\mathfrak{g}}

\def\qdim{{\rm qdim}}

\def\L{\mathcal{L}}
\def\pp{\mathfrak{p}}
\def\qq{\mathfrak{q}}

\def\L{\mathcal{L}}

\def\rN{\mathrm{N}}
\def\S{\mathcal{S}}
\def\L{\mathsf{L}}
\def\sR{\mathsf{R}}
\def\scrU{\mathscr{U}}
\def\si{\sigma}

\def\A{\mathsf{A}}
\def\O{\mathcal{O}}

\def\sl{\mathfrak{sl}}
\def\Ch{\mathrm{Ch}}

\def\bt{\boxtimes}
\newcommand{\fusion}[3]{{\binom{#3}{#1\;#2}}}

\def\<{\langle}
\def\>{\rangle}

\allowdisplaybreaks

\title{Global generation of coinvariant vector bundles  on $\overline{\mathrm{M}}_{0,n}$}
\author{\vspace{-1 em}Avik Chakravarty}
	\address{Avik Chakravarty, Department of Mathematics, University of Pennsylvania, Philadelphia, PA, 19104}
	\email{avikchak@sas.upenn.edu}

	\author{Jianqi Liu}
	\address{Jianqi Liu, Department of Mathematics, University of Pennsylvania, Philadelphia, PA, 19104}
	\email{jliu230@sas.upenn.edu}
\date{}
\begin{document}
\maketitle
\vspace{-3 em}
\begin{abstract}
Using fusion closed-ness and strong generation properties, we identify collections of modules over vertex operator algebras (VOAs) whose induced vector bundles of coinvariants on the moduli space $\overline{\mathrm M}_{0,n}$ of stable $n$-pointed curves of genus zero are globally generated.  These results specialize to \cite{Fak12} for affine VOAs at positive integer levels, and of \cite{DG23} for strongly rational VOAs that are generated in degree 1. To illustrate the results, and to show that fusion-closedness is essential, we give numerous examples, including Lattice VOAs, Simple Current Extensions of Affine VOAs, and various tensor products.  

\end{abstract}


\section{Introduction}

Vector bundles on an algebraic variety $X$ are fundamental objects in algebraic geometry. When such a bundle is \emph{globally generated}, its sections induce a morphism $X\to Y$ to another projective variety $Y$, and studying this morphism can reveal important geometric information about $X$.
There are several ways to construct globally generated vector bundles on the moduli space $\overline{\mathrm M}_{0,n}$ of stable genus zero curves with $n$ marked points; see \cite[Section 1.1]{Fed15} for an overview. A separate construction comes from the representation theory of certain $\nn$-graded infinite-dimensional algebras over $\cc$, called \emph{vertex operator algebras}. Under mild assumptions on a vertex operator algebra $V$ (a VOA), any choice of $n$ representations $M^1,\ldots,M^n$ determines a \emph{coinvariant vector bundle} $\mathbb{V}(V;M^1\otimes\cdots\otimes M^n)$ on $\overline{\mathrm{M}}_{g,n}$ for every genus $g$ 
\cite{TUY89,DGT24}. Fakhruddin \cite{Fak12} proved that the  bundle
\[\mathbb{V}:=\mathbb{V}(V;M^1\otimes\cdots\otimes M^n)\longrightarrow\overline{\mathrm M}_{0,n}\]
is globally generated for affine VOAs at positive integer levels, an important class of VOAs generated in degree $1$. This phenomenon was subsequently extended in \cite{DG23} to other VOAs generated in degree $1$, including Heisenberg VOAs. However, not all coinvariant bundles arising from VOAs generated in higher degrees are globally generated. See \cite{DG23} and \cite{ChakParafermion} for examples involving rank one lattice VOAs and parafermion VOAs, respectively. In this paper, we establish criteria for the global generation of coinvariant bundles $\mathbb{V}$ and their determinant line bundles $\det(\mathbb{V})$ associated with VOAs generated in arbitrary degrees. Our criteria impose suitable conditions on the choice of $V$-modules, and we recover the results of \cite{Fak12,DG23} as special cases.
\subsection{Main theorems}
\label{subsec: intro_main_theorems}
Theorem A builds on the notion of \emph{proportional pairings} between two different VOAs and their fusion subrings introduced in \cite{ChakParafermion} (see Definition~\ref{def: pairing}). The fusion ring $\mathscr {F}(V)$ of a VOA $V$ is a ring that encodes the tensor-product rules among its representations (see Definition~\ref{def:fusionsubring}).

\begin{maintheoremA*}
Let $V$ and $\tilde{V}$ be two strongly rational\footnote{See the definition of a strongly rational VOA in \cite[Section 2.2]{DG23}.} VOAs that form a proportional pairing given by a ring isomorphism $f:\cyc{W^1,\ldots,W^n}\overset{\sim}{\longrightarrow} \cyc{M^1,\ldots,M^n}$ with a proportionality constant $k\in\Z_{\geq1}$ between subrings $\cyc{M^1,\ldots,M^n}\subset \mathcal{F}(V)$ and $\cyc{W^1,\ldots,W^n}\subset \mathcal{F}(\tilde{V})$. Then, 
\[
\det \mathbb{V}:=\det\mathbb{V}(V;M^1\otimes\cdots\otimes M^n)\cong \det \mathbb{V}(V',W^1\otimes\cdots\otimes W^n)^{\otimes k}
\]
as line bundles on $\overline{\mathrm{M}}_{0,n}$. Moreover, if $\det\mathbb{V}(\tilde{V},W^1\otimes\cdots\otimes W^n)$ is globally generated, then so is $\det\mathbb{V}$. 
\end{maintheoremA*}

At first sight, the role of the fusion-subring structure in Theorem~A is somewhat mysterious, especially because the construction of the bundle $\mathbb{V}$ does not \emph{a priori} involve a choice of fusion subring. The criterion in Theorem B provides an intrinsic explanation for their appearance; see also Remark~\ref{rem: highlight}.
\begin{maintheoremB*}
Let $V$ be a strongly rational VOA, $\mathscr W$ the set of irreducible $V$-modules, and  $\S \ssq \mathscr W$ be a strongly geometrically reducible subset (see Definition~\ref{df:stronglyreducible}). For any $M^1,\ldots,M^n\in\mathcal S$, the natural morphism
$\Psi:\bigotimes_{i=1}^nM^i(0)\otimes\mathcal O_{\overline{\mathrm M}_{0,n}}
	\longrightarrow
	\mathbb V(V;M^1\otimes\cdots\otimes M^n)$
is surjective. In particular,
$\mathbb V(V;M^1\otimes\cdots\otimes M^n)$
is globally generated on $\overline{\mathrm M}_{0,n}$.
\end{maintheoremB*}

In the theorem above, $M(0)$ denotes the lowest-degree piece of a graded $V$-module $M=\bigoplus_{i\in\mathbb{N}}M(i)$. For positive integer level affine VOAs, a VOA generated in degree $1$, the restriction of the natural morphism $\Psi$ to the smooth locus $\mathrm{M}_{0,n}$ was introduced by Tsuchiya--Ueno--Yamada, who proved in \cite[Proposition 3.5.1]{TUY89} that it is surjective. Fakhruddin subsequently extended this surjectivity to the boundary of $\overline{\mathrm{M}}_{0,n}$ in \cite[Lemma 2.5]{Fak12}. More recently, Damiolini and Gibney proved that the morphism $\Psi$ is a surjection for any $n$-choices of $V$-modules of a vertex operation algebra $V$ strongly generated in degree $1$ \cite[Theorem 1]{DG23}. Although we use the same natural morphism $\Psi$, our proof of surjectivity is necessarily different from those in the existing literature, since we allow $V$ to be generated in higher degree and accordingly impose additional restrictions on the irreducible $V$-modules under consideration. A sketch of the argument is given in Section~\ref{subsec: sketch of Thm B proof}. Theorem B also recovers existing global-generation results of \cite{Fak12} and of \cite{DG23} for strongly rational VOAs, as shown in Example~\ref{exmp: GR and affine VOAs}.

\medskip

Neither Theorem A nor Theorem B requires the VOA to be generated in degree $1$. Instead, both results impose conditions on the $V$-modules under consideration. Theorem A uses the Chern character formula for coinvariant bundles from \cite{DGT3} to impose proportionality relations between the Chern character of $\mathbb{V}$ and that of $\mathbb{V}(V',W^1\otimes\cdots\otimes W^n)$. See Theorem~\ref{thm: K-class result} for the precise relation. Theorem A above is a corollary of Theorem~\ref{thm: K-class result}; see Corollary~\ref{cor: main corollary}. The argument relies crucially on \cite[Proposition 3.5]{ChakParafermion}, which shows that factorization for coinvariant bundles on $\overline{\mathrm M}_{0,n}$ preserves fusion-subring structures.

\medskip

Theorem~B establishes global generation of the coinvariant vector bundles $\V$ themselves, and hence of their determinant line bundles $\det(\V)$ using a new notion of geometrically reducible $V$-modules. A subset $\S$ of the set of irreducible $V$ -modules $ \mathscr W$ is called {\em geometrically reducible} if, for every $n\geq3$ and every choice of modules $M^1,\ldots,M^n\in\S$, the natural morphism $\Psi|_{\mathrm M_{0,n}}$ is surjective; see Definition~\ref{df:geored}. We call $\S$ {\em strongly geometrically reducible} if it is closed under the fusion tensor product, as in Definition~\ref{def:Sclosedunderfusion}, and satisfies the sufficient conditions for geometric reducibility specified in Definition~\ref{df:stronglyreducible}. These conditions reflect the two stages of the proof. On smooth curves, the conditions on the modules allow us to reduce coinvariants to the bottom degrees, using techniques from restricted conformal blocks \cite{Liu23,GLZ25,GLZ26}. On nodal curves, factorization \cite{DGT24}, expressed through the mode transition algebra $\Phi(\A)$, introduces additional module insertions at the preimages of each node. Fusion closedness ensures that only irreducible modules belonging to $\S$ contribute to $\Phi(\A)$  \cite{ChakParafermion}, while the stronger reduction conditions allow the resulting insertions to be reduced within the inductive argument. Together, these properties establish the surjectivity of $\Psi$ on every boundary fiber.

\subsection{Sketch of the proof of Theorem A}
A \emph{proportional pairing} $$f: \cyc{W^1,\ldots,W^n}\overset{\sim}{\longrightarrow} \cyc{M^1,\ldots,M^n},$$ given by $f(W^i)=M^i$, comes equipped with a \emph{proportionality constant} $k\in\mathbb{Z}_{\geq 1}$, encoded by the relation $\operatorname{cw}(f(W^i))=k\,\operatorname{cw}(W^i)$ between the conformal weights of the modules $W^i$ and $M^i$, which are both rational numbers. In \cite{ChakParafermion}, it was shown that $c_1(\mathbb{V})=k\,c_1\bigl(\mathbb{V}(\tilde{V},W^1\otimes\cdots\otimes W^n)\bigr)$ in $H^2(\overline{\mathrm M}_{0,n};\mathbb{Q})$. This follows from comparing the coefficients of the tautological classes in the first Chern class formula for coinvariant vector bundles given in \cite{DGT3} and using the fact that factorization theorem preserves the subring structure \cite{ChakParafermion}. We show that similar relations extend to the level of Chern characters. More precisely, Theorem~\ref{thm: K-class result} shows that the $K_0$-class of $\mathbb{V}$ in $K_0(\overline{\mathrm{M}}_{0,n})\otimes\mathbb{Q}$ is given by $[\mathbb{V}]=[\psi^k\,\mathbb{V}(\tilde{V},W^1\otimes\cdots\otimes W^n)]$, where $\psi^k$ denotes the $k$-th Adams endomorphism. As a consequence, comparing the first Chern classes of both bundles yields the isomorphism $\det \mathbb{V}\cong \det\bigl(\mathbb{V}(\tilde{V},W^1\otimes\cdots\otimes W^n)\bigr)^{\otimes k}$; see Corollary~\ref{cor: main corollary}.

\subsection{Sketch of the proof of Theorem B}
\label{subsec: sketch of Thm B proof}
 For any $n$ simple $V$-modules $W^1,\ldots,W^n\in\mathscr{W}$, we use $W^{\bullet}$ to denote the tensor product $W^1\otimes\cdots\otimes W^n$. We consider the sufficient conditions on the $V$-modules $W^1,\ldots,W^n$ so that the morphism 
$
\Psi:\bigotimes_{i=1}^n W^{\bullet}(0)\otimes \mathcal O_{ \overline{\mathrm M}_{0,n}} \to \mathbb{V}(V;W^{\bullet})
$
of vector bundles on $\overline{\mathrm M}_{0,n}$ is a surjection. 
To prove surjection of $\Psi$, it suffices to prove that the morphism 
\[
\widetriangle{\Psi}: \bigotimes_{i=1}^n W^{\bullet}(0)\otimes\widetriangle{\mathrm M}_{0,n}\longrightarrow \widetriangle{\mathbb{V}}(V;W^{\bullet})
\]
of sheaves on the moduli space $\widetriangle{\mathrm M}_{0,n}$ of coordinatized pointed curves is surjective. 
The moduli space $\widetriangle{\mathrm M}_{0,n}$ parametrizes triples $(C,P_{\bullet},t_{\bullet})$, where $(C,P_{\bullet})\in \overline{\mathrm M}_{0,n}$ and $t_{\bullet}$ is an $n$-tuple of formal coordinates $t_i$ at each of the marked points $P_i$. A detailed description of $\widetriangle{\mathrm M}_{0,n}$ can be found in \cite[Section 2.2.2]{DGT24}. The morphism $\widetriangle{\Psi}$ was used in \cite{DG23} to prove global generation of coinvariant bundles on $\overline{\mathrm M}_{0,n}$ arising from representations of VOAs generated in degree $1$. Even though we use the same morphism, our procedure for proving surjectivity of $\widetriangle{\Psi}$ is different from that used in \cite{DG23}.

\medskip

In Section~\ref{subsec: smooth case}, we prove surjectivity of the morphism $\widetriangle{\Psi}$ for the smooth coordinatized curves $(C,P_{\bullet}, t_{\bullet})$, that is, $C$ is a smooth rational curve. First,  Proposition~\ref{prop:smoothdecent} asserts that all but one of the modules appearing in the space of coinvariants $\widetriangle{\mathbb{V}}(V,W^{\bullet})_{(C,P_{\bullet},t_{\bullet})}$ may be replaced by their bottom-degree pieces.
\[
\widetriangle{\mathbb{V}}(V;W^{\bullet})_{(C,P_{\bullet},t_{\bullet})} = [W^1\otimes \bigotimes_{i=2}^n W^i(0)]. 
\]
See Notation~\ref{notn: smooth case} for details on the right hand side of the above identity. The result essentially follows from a Riemann-Roch calculation discussed in \cite[Section 1.10]{DGT24} together with the idea of restricted conformal blocks \cite[Theorem 5.19]{GLZ26}. This motivates the notion of a \emph{geometrically reducible} (cf. Definition~\ref{df:geored}) collection of modules: for such a collection, the remaining module $W^1$ can also be reduced to its bottom degree, that is, the natural map
\[
\widetriangle{\Psi}_{(C,P_{\bullet},t_{\bullet})}: W^{\bullet}(0)\longrightarrow \widetriangle{\mathbb{V}}(V;W^{\bullet})_{(C,P_{\bullet},t_{\bullet})}
\]
is surjective. Theorem~\ref{thm:criterion1} gives a representation-theoretic criterion for geometric reducibility in terms of strong generation from the bottom degree and the action of the nonnegative modes of a chosen set of strong generators. Specifically, the criterion requires that the $V$-module $W$ is strongly generated  by the chosen strong generators $U_+$ of the VOA $V$, and the nonnegative modes of these generators acting on $W(0)$ land inside the $V_1$-generated submodule $\mathscr{U}(\widehat{V}_1). W(0)$. 

\medskip



Section~\ref{subsec: nodal case} is at the heart of the proof of
Theorem~B. In this section, we prove that
$\widetriangle{\Psi}_{(C,P_{\bullet},t_{\bullet})}$ is surjective
for every nodal $n$-pointed rational curve $(C,P_{\bullet})$. Equivalently, denote the fiber $\widetriangle{\mathbb{V}}(V;W^{\bullet})_{(C,P_{\bullet},t_{\bullet})}$ by  $[W^\bullet]_{\mathcal L_{C\bs P_\bullet}(V)}=W^\bullet/\mathcal{L}_{C\bs P_\bullet}(V).W^\bullet$, we prove that it has a reduction 
\[
[W^\bullet]_{\mathcal L_{C\bs P_\bullet}(V)}=[W^\bullet(0)]=\spn\{[v^\bullet]: v^i\in W^i(0), 1\leq i\leq n\}.
\]
When $(C,P_\bullet)=(X,P_\bullet\sqcup Q_\bullet)$ is a rational curve with a single node $\qq$, the factorization
description inserts the induced modules
$\Phi^\L(\A)$ and $\Phi^\sR(\A)^{\theta},$
associated with the Zhu algebra $\A$ of $V$ \cite{Z96}, at the
preimages $\qq_+$ and $\qq_-$ of the node under the normalization map
\cite{DGT24,DGK25}. 
\medskip

If the marked modules $M^\bullet$ and $W^\bullet$ belong to a
fusion-closed subset $\mathcal S$, then modules outside $\mathcal S$
cannot contribute to the spaces of coinvariants on smooth curves
\cite[Proposition~3.5]{ChakParafermion}. Accordingly, we introduce
the ideal
\[
\A_{\mathcal S}
=
\prod_{W\in\mathcal S}\operatorname{End}(W(0))
\subset
\prod_{W\in\mathscr W}\operatorname{End}(W(0))
=
\A
\]
of the Zhu algebra. The induced modules $\Phi^\L(\A)$ and
$\Phi^\sR(\A)^\theta$ attached at $\qq_+$ and $\qq_-$ can then be
replaced by $\Phi^\L(\A_{\mathcal S})$ and
$\Phi^\sR(\A_{\mathcal S})^\theta$, respectively;
see Lemma~\ref{lm:PhiAandAS}.


\medskip 

Lemma~\ref{lm:redAS} is the key to our reduction process. Since we
have replaced the Zhu algebra $\A$ by
$\A_\S=\prod_{W\in \S} W(0)\o W(0)^\ast$, the generation conditions
$(\S2)$ and $(\S3)$ for irreducible modules $W\in \S$ allow us to
reduce the fiber
$[M^\bullet\o \Phi^\ast(\A_\S)]_{\mathcal{L}_{C\bs P_\bullet\sqcup Q}(V)}$,
where $\Phi^\ast(\A_{\mathcal S})$ denotes either
$\Phi^\L(\A_{\mathcal{S}})$ or $\Phi^\sR(\A_{\mathcal{S}})^{\theta}$
and is attached at a marked point obtained by normalizing a node,
to a much smaller subspace
\[
[M^\bullet(0)\o 1_{\A_\S}]
\]
in comparison with the bottom degree of
$M^\bullet\o \Phi^\ast(\A_\S)$, which is $M^\bullet(0)\o \A_\S$.
This explains why the definition of
\emph{strongly geometrically reducible} subsets in
Definition~\ref{df:stronglygen} has three conditions. Fusion
closedness of $\S$ helps us pick out a smaller set of $V$-modules
that could satisfy the good conditions $(\S2)$ and $(\S3)$, while
these conditions allow us to achieve a strong reduction of the
coinvariants to $1_{\A_\S}$.
Combining this reduction with Lemmas~\ref{lm:subring}
and~\ref{lm:PhiAandAS}, together with the fact that the isomorphism
between nodal coinvariants and the coinvariants on the normalized
curve is given by
$[u^\bullet\o v^\bullet]\mapsto
[u^\bullet\o 1_{\A_\S}]\o_\A[1_{\A_\S}\o v^\bullet]$,
Proposition~\ref{prop:onenode} proves the desired bottom-degree
surjection for a curve with one node.

\medskip

Building on this reduction to $[M^\bullet(0)\o 1_{\A_\S}]$,
Lemma~\ref{lm:onenodespecialred} and Proposition~\ref{prop:onenodeAS}
strengthen this one-node reduction by retaining an additional transition
module $\Phi^*(\A_{\mathcal S})$. This strengthened form serves as the
base case for the induction on the number of nodes of
$(C,P_\bullet)\in \partial \overline{\mathrm{M}}_{0,n}$.
Lemma~\ref{lm:notinSvanishing} extends the vanishing statement of
Lemma~\ref{lm:subring} to stable nodal rational curves: if all marked
modules except one belong to the fusion-closed subset $\S$, and the
remaining module lies outside $\S$, then the space of coinvariants
vanishes. Using this lemma, we prove in
Proposition~\ref{prop:arbitrarynodal} that, for
$(C,P_\bullet\sqcup \pp)\in \partial \overline{\mathrm{M}}_{0,n}$,
where $\pp$ is a smooth point, the space of coinvariants
$[M^\bullet\o \Phi^\ast(\A_\S)]_{\mathcal{L}_{C\bs P_\bullet\sqcup \pp}(V)}$
can be reduced to $[M^\bullet(0)\o 1_{\A_\S}]$ as well. This is the
final preparatory step for our induction argument in
Theorem~\ref{thm:main1}.

\medskip

Now, given any
$(C,P_\bullet)\in \partial \overline{\mathrm{M}}_{0,n}$, let $\qq$
be a node and $\tilde{C}=C_+\sqcup C_-$ be the normalization at $\qq$.
Assume that $\bullet=\star\sqcup \circ$, where $P_\star\subset C_+$
and $P_\circ\subset C_-$. Then, by the identification of nodal
coinvariants due to Damiolini--Gibney--Krashen \cite{DGK25}, the fiber
at $(C,P_\bullet)$ can be identified as
$$
[M^\star\o M^\circ]_{\mathcal{L}_{C\bs P_\bullet}(V)}
\cong
[M^\star \o \Phi^\L(\A_\S)]_{\mathcal{L}_{C_+\bs P_\star\sqcup \qq_+}(V)}
\o_\A
[\Phi^\sR(\A_\S)^\theta\o M^\circ]_{\mathcal{L}_{C_-\bs \qq_-\sqcup P_\circ}(V)},
$$
where the isomorphism is given by
$[u^\star\o u^\circ]\mapsto
[u^\star\o 1_{\A_\S}]\o_\A [1_{\A_\S}\o u^\circ]$.
Since $C_\pm$ are either smooth or have fewer nodes than $C$, the
induction hypothesis, together with
Proposition~\ref{prop:arbitrarynodal} and Lemma~\ref{lm:redAS},
allows us to reduce
$[M^\star \o \Phi^\L(\A_\S)]_{\mathcal{L}_{C_+\bs P_\star\sqcup \qq_+}(V)}$
and
$[\Phi^\sR(\A_\S)^\theta\o M^\circ]_{\mathcal{L}_{C_-\bs \qq_-\sqcup P_\circ}(V)}$
to $[M^\star(0)\o 1_{\A_\S}]$ and
$[1_{\A_\S}\o M^\circ(0)]$, respectively. The identification map
above then shows that
$[M^\star\o M^\circ]_{\mathcal{L}_{C\bs P_\bullet}(V)}$
can be reduced to its bottom degree. This completes the proof that
$\widetriangle{\Psi}_{(C,P_{\bullet},t_{\bullet})}$ is surjective;
see Theorem~\ref{thm:main1}. Finally, if $V=L_k(\g)$ is the affine VOA at a positive integral level $k$, then we can take $\S=\mathscr W$. In other words, the full set
of irreducible $V$-modules is strongly geometrically reducible.
Hence, Theorem~\ref{thm:main1} recovers the main global generation
results of \cite{Fak12,DG23}; see Corollary~\ref{cor:degreeone}.

\subsection{Examples}
We illustrate Theorems~A and~B with examples arising from lattice VOAs and simple-current extensions of affine VOAs. Together, these theorems yield globally generated coinvariant bundles of both rank one and higher rank for VOAs that are not generated in degree one. We also give a counterexample showing that the fusion-closure condition on $\S$ in Theorem~B cannot be omitted.

\subsubsection{Globally generated coinvariants bundles from strongly geometrically reducible sets}\label{exmp: intro GRS}
Consider the simple-current extension $V_e=L_8(\mathfrak{sl}_2)\oplus L_8(8)$ \cite{Li97}, which is generated in degrees $1$ and $2$. Its irreducible modules are $\mathscr{W}=\{V_e,X,W^+,W^-\}$, where both $W^\pm$ restrict to the affine module $L_8(4)$ and satisfy conditions $(\S2)$ and $(\S3)$. Since $W^\pm\boxtimes_{V_e}W^\pm\cong V_e\oplus W^\pm$, each subset $\mathcal{S}^\pm=\{V_e,W^\pm\}$ is strongly geometrically reducible; see Example~\ref{ex:VeGG}. Theorem~B therefore shows that $\mathbb{V}(V_e;(W^\pm)^{\otimes n})$ is globally generated over $\overline{\mathrm{M}}_{0,n}$ for every $n\geq3$. Moreover, these bundles have rank greater than one for $n\geq4$. 

The same VOA provides a counterexample showing the importance of fusion closure. Although $\mathcal{S}^+$ and $\mathcal{S}^-$ are strongly geometrically reducible, their union $\mathcal{S}=\{V_e,W^+,W^-\}$ is geometrically reducible but not fusion-closed, since $W^+\boxtimes_{V_e}W^-\cong X\notin\mathcal{S}$. The natural morphism $\Psi$ is therefore surjective on the smooth locus for marked modules in $\mathcal{S}$, but the bundle defined by mixed modules $\mathbb{V}:=\mathbb{V}_{0,8}(V_e;(W^+)^{\otimes4}\otimes(W^-)^{\otimes4})$ is not globally generated. Indeed, for the $F$-curve whose four tails each carry one $W^+$ and one $W^-$ insertion, factorization gives
$
c_1(\mathbb{V})\cdot F
=\deg\mathbb{V}_{0,4}\bigl(V_e;X^{\otimes4}\bigr)=-1.
$
Thus surjectivity on the smooth locus need not extend across the boundary of the moduli space. Geometric reducibility nevertheless implies effectiveness of the first Chern class \cite{DbPrivate}, a perspective that the authors will explore further in forthcoming joint work with Daebeom Choi.

Let $L=\mathbb{Z}\alpha$ with $(\alpha\mid\alpha)=2k$, where $k\geq1$, and let $V_L$ be the associated lattice VOA \cite{FLM88,dong93}. By Proposition~\ref{prop:rankonelattice}, the subset $\mathcal{S}=\{V_L,V_{L+\frac{\alpha}{2}}\}$ is strongly geometrically reducible. Hence Theorem~B gives global generation of $\mathbb{V}(V_L;M^1\otimes\cdots\otimes M^n)$ for every choice of $M^1,\ldots,M^n\in\mathcal{S}$. For $k\geq2$, the VOA $V_L$ is not generated in degree one. See Examples~\ref{ex:latticesqrtKA2} and~\ref{ex:latticeA1+A1} for additional examples given by rank-two lattices. 


We also obtain an infinite family from the diagonal simple-current extensions $V^{(m)}=(V\otimes V)\oplus(J_m\otimes J_m)$, where $V=L_{2m}(\mathfrak{sl}_2)$, $J_m=L_{2m}(2m)$, and $m\geq2$. The induced module $X^{(m)}=\operatorname{Ind}_{V\otimes V}^{V^{(m)}}(V\otimes J_m)$ satisfies $X^{(m)}\boxtimes_{V^{(m)}}X^{(m)}\cong V^{(m)}$, and $\{V^{(m)},X^{(m)}\}$ is strongly geometrically reducible by Proposition~\ref{prop:digonalsimplecurrect}. Theorem~B therefore applies to all choices of marked modules in this subset. Strong geometric reducibility is also preserved under VOA-tensor products \cite{FHL93} by Proposition~\ref{prop:tensorproduct}, giving further examples. These applications are collected in Section~\ref{exmp: list of examples}.

\subsubsection{Globally generated determinant bundles from proportional pairing with affine VOAs}
We study subrings $\mathcal{R}$ of the fusion ring $\mathcal{F}(V)$ that are proportionally paired with subrings of the fusion rings of affine VOAs, with proportionality constant $k\in\zz_{\geq1}$. By Theorem~A, for any simple $V$-modules $M^1,\ldots,M^n$ whose classes lie in $\mathcal{R}$, the determinant bundle $\det\mathbb{V}(V;\bigotimes_i M^i)$ is isomorphic to $\det(\mathbb{V}')^{\otimes k}$ for a corresponding affine coinvariant bundle $\mathbb{V}'\to\overline{\mathrm M}_{0,n}$. Since $\det(\mathbb{V}')$ is globally generated by \cite{Fak12,DG23}, so is $\det\mathbb{V}(V;\bigotimes_i M^i)$.

We apply this result to an infinite family of lattices $\{L(k,r)\mid k,r\in\zz_{\geq1}\}$, where $k$ determines the bilinear pairing and $r$ denotes the rank. For each pair $(k,r)$, the corresponding lattice VOA $V_{L(k,r)}$ is generated in degrees $1$ and $k$. 

For $L=L(k,r)$, Proposition~\ref{prop: subring str} identifies a subring $\mathcal{S}_r$ of $\mathcal{F}(V_L)$ which is proportionally paired with the fusion ring of $L_1(\mathfrak{sl}_{r+1})$ with proportionality constant $k$. This is the full fusion ring when $k=1$ and a proper subring when $k>1$. Since every nonzero coinvariant bundle arising from a lattice VOA has rank one, Theorem~A gives global generation of $\mathbb{V}(V_L;\bigotimes_{i=1}^n W^i)$ for any choice of simple modules $W^i$ whose classes lie in $\mathcal{S}_r$.

\begin{remark}\label{rem: highlight}
	Setting $r=1$ above, the set $\mathcal{S}_1$ coincides with the set $\mathcal{S}=\{V_L,V_{L+\alpha/2}\}$ of strongly geometrically reducible $V_L$-modules from Section~\ref{exmp: intro GRS}. This illustrates the main role of Theorem B: it provides an intrinsic vertex operator algebraic explanation for the cohomological phenomenon detected by Theorem A. On the other hand, the principal advantage of Theorem A is that its criterion is combinatorial and can often be verified through relatively elementary algebraic computations, which allows us to relate global generation of (often in practice) an infinite collection of determinant bundles coming from VOAs. 
\end{remark}

\subsubsection{Positivity for rank $1$ lattice VOAs}
If a vector bundle $\mathbb{V}$ is globally generated, then its first Chern class $c_1(\mathbb{V})$ is represented by a semi-ample divisor. It is therefore natural to ask whether Theorems~A and~B recover all semi-ample divisors arising from the VOAs considered above. To investigate this question, we study positive divisors associated with rank $1$ lattice VOAs. Rank $1$ lattice VOAs form a discrete family indexed by positive integers $k\in\mathbb{Z}_{\geq 1}$. The corresponding lattice is determined by the bilinear pairing $(\alpha\mid\alpha)=2k$. In Section~\ref{subsec: positivity rank 1 lattice}, we exhibit a collection of nontrivial semi-ample bundles, for all integers $k\geq1$, that cannot be obtained from Theorems~A and~B. 

The collection does not give all nontrivial coinvariant bundles arising from rank $1$ lattice VOAs. In Section~\ref{subsec: classification rank 1 lattice k is 2}, we give a complete classification (Proposition~\ref{prop: classification in k2 rank1}) of the semi-ample divisors on $\overline{\mathrm M}_{0,n}$ arising from the representation theory of the rank $1$ lattice VOA with $k=2$. An analogous classification for larger values of $k$ appears to be inaccessible because of the rapidly increasing combinatorial complexity. 


\subsection{A reader's guide}
The paper is organized as follows. In Section~\ref{sec: prelim}, we recall the necessary background on coinvariant bundles, factorization, fusion rings, proportional pairings, and chiral Lie algebras. In Section~\ref{sec: VOA approach}, we develop the VOA theoretic approach to global generation and prove Theorem~B, including the smooth and nodal reduction arguments. Section~\ref{sec: examples} studies geometric reducibility for lattice VOAs, simple current extensions, and tensor product VOAs, and provides both examples and a non-example illustrating the role of fusion closedness. In Section~\ref{sec: cohom approach}, we prove the $K$-theoretic comparison underlying Theorem~A and study its applications to type $A$ lattice VOAs and to the positivity of rank $1$ lattice VOA line bundles.

\subsection*{Acknowledgments}
We are grateful to Angela Gibney for valuable discussions and helpful feedback on an earlier draft of this paper. We used torsorlab \cite{torsorlab} with Claude to compare the paper with the literature and review the draft. The authors checked all suggestions and are responsible for the content. We used ChatGPT (GPT-6) and Claude to explore examples arising from our main Theorem B and to assist in identifying and correcting minor errors. We emphasize that the proofs of the main theorems, the overall structure of the paper, and the vast majority of the examples are the result of the authors' independent work and collaboration over the past ten months. The authors take full responsibility for the correctness of all theorems, proofs, and examples presented in this paper.


\section{Preliminaries}
\label{sec: prelim}
\subsection{Vertex operator algebras}
We refer to \cite{FLM88,FBZ04,FHL93,Z96,DLM98,DGT1,DGT24} for the definition and properties of vertex operator algebras (VOAs) and the sheaves of VOA-coinvariants on the moduli space of stable curves $\overline{\mathrm M}_{g,n}$.

Throughout this paper, we assume that $V$ is a strongly rational VOA, and we write
$$\mathscr{W}=\{V=W^0,W^1,\ldots,W^N\}$$
for the set of all irreducible $V$-modules up to isomorphism. For $V$-modules $M^1,\dots,M^n$ we write $\mathbb{V}(V;M^1\otimes\cdots\otimes M^n)$, or $\mathbb{V}_{g,n}(V,M^1\otimes\cdots\otimes M^n)$ when we wish to record the genus $g$ and the number $n$ of marked points, for the sheaf of coinvariants on $\overline{\mathrm{M}}_{g,n}$. We abbreviate $M^\bullet:=M^1\otimes\cdots\otimes M^n$ and $M^\bullet(0):=M^1(0)\otimes\cdots\otimes M^n(0)$, and for a subset $\mathcal{S}\subseteq\mathscr{W}$ we write $M^\bullet\subset \mathcal{S}$ to denote a tuple $(M^1,\ldots,M^n)\in \mathcal{S}^{\oplus n}$. For a vector bundle $\mathbb{V}$ we write $\mathbb{V}^{\vee}$ for its dual; the dual of a bundle of coinvariants is the bundle of conformal blocks. We write $L_k(\g)=L_{\hat{\g}}(k,0)$ for the level $k$ simple affine VOA associated with a simple Lie algebra $\g$, and $\cw(M)$, or $h_M$, for the conformal weight of a simple module $M$.

\subsection{Coinvariant bundles}
The Chern character of the bundles of coinvariants is computed in \cite{DGT3}; we recall the formula in Section~\ref{subsec: proof of Theorem A}, where it is used.
\begin{theorem}[\cite{DGT3}]\label{thm: strongly rational}
For any $n$-tuple of irreducible modules $M^1,\dots, M^n\in \mathscr{W}$ (possibly with repetitions), the sheaf of coinvariants $\mathbb{V}_{g,n}(V,M^1\otimes\cdots\otimes M^n)$ is a vector bundle on $\overline{\mathrm M}_{g,n}$.
\end{theorem}

\begin{theorem}[Propagation of Vacua \cite{Cod19,DGT1}]
\label{thm: POV}
Let $M^1,\dots, M^n\in \mathscr{W}$. Then
\[
\mathbb{V}_{g,n+1}(V,M^1\otimes \dots \otimes M^n\otimes V) \cong \pi^*\mathbb{V}_{g,n}(V, M^1\otimes \dots \otimes M^n),
\]
where $\pi:\overline{\mathrm{M}}_{g,n+1}\to\overline{\mathrm{M}}_{g,n}$ is the forgetful map that forgets the $(n+1)$-th marked point.
\end{theorem}
\begin{theorem}[Factorization theorem \cite{TUY89, BFM, DGT1,DGT24,DGT3}]
\label{thm: factorization}
For any $n$ representations $M^1,\ldots,M^n$ of a VOA $V$, any integer $0\leq g_1\leq g$, and any partition $I_1\sqcup I_2=\{1,\ldots,n\}$ with $\min(|I_1|,|I_2|)\geq2$ if $g_1=0$ or $g_1=g$, let $\xi:\overline{\mathrm{M}}_{g_1,|I_1|+1}\times\overline{\mathrm{M}}_{g-g_1,|I_2|+1}\to\overline{\mathrm{M}}_{g,n}$ be the gluing map. Then
\[
\xi^*\mathbb{V}_{g,n}(V,M^1\otimes\cdots\otimes M^n)\cong\bigoplus_{W} \mathbb{V}_{g_1,|I_1|+1}(V,M^{I_1}\otimes W)\boxtimes \mathbb{V}_{g-g_1,|I_2|+1}(V,M^{I_2}\otimes W'),
\]
where $M^{I_j}:=\bigotimes_{i\in I_j}M^i$, the sum runs over all simple $V$-modules $W$, and $W'$ denotes the contragredient dual of $W$.
\end{theorem}

\begin{definition}\label{def:fusionsubring}
   Let $\mathscr{F}(V)=\bigoplus_{i=0}^N \mathbb Z [W^i]$, where $[W^i]$ is a formal symbol associated to the irreducible $V$-module $W^i$.  Define a multiplication on $\mathscr{F}(V)$ on the generators by
	\begin{equation}\label{eq:fusionmult}
		[W^i]\cdot [W^j]:=\sum_{k=0}^N N_{W^iW^j}^{W^k}\, [W^k],\quad 0\le i,j\le N,
	\end{equation}
and extend it bilinearly. It follows from \cite{FHL93,H05,H08,GL25,DGT3} that $\mathscr F(V)$ is a commutative unital ring, which is referred to as the \textbf{fusion ring} of $V$.

Throughout this paper, when we say $\mathcal S$ is a {\bf subring of the fusion ring $\mathscr F (V)$}, it means
\[
\mathcal S=\<V,W^{i_1},\dots, W^{i_m}\>=\mathbb{Z}[V]\op \mathbb Z [W^{i_1}]\op \ds \op \mathbb Z [W^{i_m}],
\]
for some subset $\mathcal S^{\mathrm{sim}}=\{V=W^0,W^{i_1},\dots, W^{i_m}\}\subseteq \mathscr W$, and $\mathcal S$ is closed under the fusion tensor product \eqref{eq:fusionmult}. In particular, $\mathcal S$ is not an arbitrary subring of $\mathscr{F}(V)$.
\end{definition}

\begin{proposition}[{\cite[Proposition 3.5]{ChakParafermion}}]
\label{prop: factorization respects subrings}
If the $V$-modules $M^1,\ldots, M^n$ lie in a subring $\mathcal{S}$ of the fusion ring $\mathscr{F}(V)$, then it suffices to consider the sum over simple modules $W$ in $\mathcal{S}$ for the factorization theorem.
\end{proposition}

\begin{definition}\label{def: pairing}
Let \(V^1\) and \(V^2\) be two VOAs, and let \(\mathcal{S}^1\) (respectively, \(\mathcal{S}^2\)) denote a subring of the fusion ring \(\mathscr F (V^1)\) (respectively, \(\mathscr F(V^2)\)). Fix a rational number $\eta\in\mathbb{Q}_{>0}$. We say that the pair \((\mathcal{S}^1, \mathcal{S}^2)\) forms a \emph{proportional pairing} with \emph{proportionality} $\eta$  if there exists a ring homomorphism \(f : \mathcal{S}^1 \to \mathcal{S}^2\) satisfying the following properties
\begin{enumerate}[i.]
\item $f$ is injective and $f(M)\in (\mathcal{S}^2)^{\mathrm{sim}}$ for all $M\in(\mathcal{S}^1)^{\mathrm{sim}}$, and
\item $\cw(f(M))=\eta\cdot \cw(M)$ for all $M\in (\mathcal{S}^1)^{\mathrm{sim}}$.
\end{enumerate}
\end{definition}

Recall that a divisor class $D$ on $\overline{\mathrm{M}}_{0,n}$ is called \emph{$F$-positive} (or \emph{$F$-nef}) if $D\cdot F\geq 0$ for every $F$-curve $F\subset\overline{\mathrm{M}}_{0,n}$, i.e., for every one-dimensional boundary stratum; see \cite{KM13,Fed15}. Nef divisors are $F$-positive, and semi-ample divisors are nef. The divisor $c_1\mathbb{V}$ of a globally generated vector bundle $\mathbb{V}$ is semi-ample. The following result is important for Section~\ref{sec: cohom approach}.

\begin{theorem}[\cite{Fak12}]\label{thm: Fakh}
Any coinvariant bundle on $\overline{\mathrm{M}}_{0,n}$ corresponding to $n$ choices of admissible representations of the affine VOA $L_k(\g)$ is globally generated for all $k\in\mathbb{Z}_{\geq1}$.
\end{theorem}

\subsection{Chiral Lie algebra and a consequence of Riemann--Roch}
We follow the notation in \cite{DGT24,DG23} for chiral Lie algebras.
Let $(C,P_\bullet,t_\bullet)$ be a nodal curve, with marked points $P_\bullet$ and formal coordinates $t_\bullet$ about $P_\bullet$. The ancillary Lie algebra of $V$ at $P_i$ is
\[
L_{P_i}(V)=\frac{V\o \C(\!(t_i)\!)}{\mathrm{Im} \nabla},\quad \nabla=L(-1)\o 1+1\o\frac{d}{dt_i}.
\]
The Lie bracket of $L_{P_i}(V)$ is the same as that of the Borcherds Lie algebra $L(V)$:
\[
[a_{[m]},b_{[n]}]=\sum_{j\geq 0}\binom{m}{j} (a_{(j)}b)_{[m+n-j]},\quad a,b\in V,\ m,n\in \Z,
\]
where $a_{[m]}=a\o t_i^{m}+\mathrm{Im}\nabla\in L_{P_i}(V)$. In particular, the $V$-module $M^i$ attached at $P_i$ is a module over the Lie algebra $L_{P_i}(V)$.
Let $\mathscr{V}_C$ be the vertex algebra  bundle  on the curve $C$ \cite{FBZ04}, let $\om_C$ be the dualizing sheaf on $C$, and let $\nabla$ be a logarithmic connection on $ \mathscr{V}_C$, which is a $\C$-linear map
$\nabla: \mathscr{V}_C\ra \mathscr{V}_C\o \om_C$ satisfying
$\nabla(fs)=s\o df+f\nabla(s)$ for any local sections $f$ of $\mathscr{O}_C$ and $s$ of $\mathscr{V}_C$. Then
$H^0(D^\times_{P_i},\mathscr{V}_C\o \om_C/\mathrm{Im} \nabla)\cong L_{P_i}(V),$ where $D_{P_i}^\times$ is a formal disk at $P_i$.

\begin{definition}[\cite{DGT24}]\label{def:chiralLie}
	Let $(C,P_\bullet)$ be a nodal curve, with marked points $P_\bullet=(P_1,\ds, P_n)$, and let $V$ be a VOA. The chiral Lie algebra associated to $(C,P_\bullet)$ is
$	\mathcal{L}_{C\bs P_\bullet}(V):=H^0(C\bs P_\bullet, \mathscr{V}_C\o \om_C/\mathrm{Im} \nabla).$

There exists a homomorphism of Lie algebras
\begin{equation}\label{chiralLie}
\begin{aligned}
\mathcal{L}_{C\bs P_\bullet}(V)&\ra \bigoplus_{i=1}^n H^0(D^\times_{P_i},\mathscr{V}_C\o \om_C/\mathrm{Im} \nabla)\cong \bigoplus_{i=1}^n L_{P_i}(V),
\end{aligned}
\end{equation}
given by $\si\mapsto (\si_{P_1},\ds, \si_{P_n}),$
where $\si_{P_i}\in L_{P_i}(V)=V\o \C(\!(t_i)\!)/\mathrm{Im}\nabla$ for all $i$. In particular, let $M^\bullet=(M^1,\ds, M^n)$ be a tuple of $V$-modules attached to $P_\bullet$. Then the chiral Lie algebra
	$\mathcal{L}_{C\bs P_\bullet}(V)$ acts diagonally on the tensor product space $M^1\o \ds \o M^n$ via the homomorphism \eqref{chiralLie}:
	\begin{equation}\label{chiralLieact}
		\si.(u^1\o\ds \o u^n)=\sum_{i=1}^n u^1\o\ds \o \si_{P_i}.u^i\o\ds \o u^n,
	\end{equation}
	where $\si\in \mathcal{L}_{C\bs P_\bullet}(V)$ and $u^i\in M^i$ for all $i$. The space of coinvariants with respect to the datum $(C,P_\bullet, t_\bullet, M^\bullet)$ is defined by
	\[
	\V(V;M^\bullet)_{(C,P_\bullet,t_\bullet)}=\frac{M^1\o \ds \o M^n}{\mathcal{L}_{C\bs P_\bullet}(V).(M^1\o\ds\o M^n)},
	\]
	where the chiral Lie algebra action is defined by \eqref{chiralLie} and \eqref{chiralLieact}.
\end{definition}

The action of the chiral Lie algebra can be proved to be coordinate-independent \cite{FBZ04,DGT24}, which leads to the notion of the space of coinvariants $\V(V;M^\bullet)_{(C,P_\bullet)}$.
	This space is the fiber of the sheaf of coinvariants $\V(V;M^\bullet)$ on $\overline{\mathrm{M}}_{g,n}$ at $(C,P_\bullet)$, which is proved to be a vector bundle when $V$ is strongly rational \cite{DGT24,DGK25}.

	\subsubsection{Sections of the chiral Lie algebra}
	We recall a consequence of the Riemann--Roch theorem; see \cite[Section 1.10]{DGT24}.
	Let $C$ be a smooth curve, possibly disconnected, with two sets of marked points $P_\bullet$ and $Q_\bullet$. Assume that each connected component contains at least one of the marked points $P_\bullet$, so that $C\bs P_\bullet$ is affine. Let $s_i$ be the formal coordinate at $Q_i$ for all $i$.

	Fix $r\in \Z$. Then, for any integers $d$ and $N\gg 0$, there exists a section $\mu\in H^0(C\bs P_\bullet \sqcup Q_\bullet,\om^{\o r}_C)$ such that, for a fixed $i$,
\begin{equation*}
		\mu_{Q_i}=\left(s_i^d + s_i^N\sum_{n\geq 0} \la_n s_i^n \right)(ds_i)^r,\qquad \text{and}\qquad 
		\mu_{Q_j}=s_j^N\left(\sum_{n\geq 0} \la'_n s_j^n \right)(ds_j)^r,\quad \text{for all }j\neq i.
\end{equation*}
	With the setup as above, since $C\bs P_\bullet$ is affine, we have the following description of the chiral Lie algebra; see \cite[Section 5.1]{DG23}:
	\begin{equation}\label{eq:descriptionofchiralLie}
		\mathcal{L}_{C\bs P_\bullet\sqcup Q_\bullet}(V)\cong \frac{\bigoplus_{k\geq 0} V_k\o H^0(C\bs P_\bullet\sqcup Q_\bullet, \om_C^{1-k})}{\nabla(H^0(C\bs P_\bullet\sqcup Q_\bullet, \mathscr{V}_C))}.
	\end{equation}
	As a consequence of the existence of the section $\mu$, we have the following lemma:

	\begin{lemma}\label{lmRR}
		Let $V$ be a VOA and let $(C,P_\bullet\sqcup Q_\bullet)$ be a smooth curve. Let $a\in V$ be a homogeneous element of weight $k\geq 0$, and let $d,N\in \Z$. Then there exists a section $\si=a\o \mu\in \bigoplus_{l\geq 0} V_l\o H^0(C\bs P_\bullet\sqcup Q_\bullet, \om_C^{1-l})$ such that its equivalence class in the chiral Lie algebra $\mathcal{L}_{C\bs P_\bullet\sqcup Q_\bullet}(V)$ satisfies
		\begin{equation}\label{eq:RRlocalsigma}
			\begin{aligned}
				\si_{Q_i}&=a_{[d]}+\sum_{n\geq 0} \la_n a_{[n+N]}\in L_{Q_i}(V)=V\o \C(\!(s_i)\!)/\mathrm{Im}\nabla,\\
				\si_{Q_j}&=\sum_{n\geq 0} \la'_n a_{[n+N]}\in L_{Q_j}(V)=V\o \C(\!(s_j)\!)/\mathrm{Im}\nabla,\quad j\neq i.
			\end{aligned}
		\end{equation}
	\end{lemma}

\section{VOA approach to global generation}
\label{sec: VOA approach}
In this section, we prove Theorem B; see Section~\ref{subsec: sketch of Thm B proof} for a summary of the key ideas used in the proof. The proof is detailed in Sections~\ref{subsec: smooth case} and~\ref{subsec: nodal case}. The main result is stated in Theorem~\ref{thm:main1}.

\subsection{Smooth case}
\label{subsec: smooth case}
\subsubsection{Definition of geometric reducibility}
We introduce some notation first.

\begin{notation}\label{notn: smooth case}
	Let $W^1,\ds, W^n$ be $V$-modules, and let $N^i\ssq W^i$ be a subspace for all $i$. Let $(C,P_\bullet)\in \overline{\mathrm{M}}_{0,n}$ be a stable curve with smooth marked points $P_\bullet=(P_1,\ds, P_n)$.

	In the space of coinvariants $\V(V;W^\bullet)_{(C,P_\bullet,t_\bullet)}=W^\bullet/\mathcal{L}_{C\bs P_\bullet}(V). W^\bullet$, we denote the subspace spanned by the equivalence classes of elements in $N^1\o\ds \o N^n\ssq W^1\o\ds\o W^n$ by
	\begin{equation}\label{eq:classincoinv}
		[N^1\o \ds \o N^n]:=\frac{N^1\o \ds \o N^n+\mathcal{L}_{C\bs P_\bullet}(V). W^\bullet}{\mathcal{L}_{C\bs P_\bullet}(V). W^\bullet}\ssq \V(V;W^\bullet)_{(C,P_\bullet, t_\bullet)}.
	\end{equation}

\end{notation}

Based on the idea of restricted VOA-conformal blocks \cite{Liu23,GLZ25,GLZ26}, we have the following simplification of coinvariants on smooth curves.
\begin{proposition}\label{prop:smoothdecent}
	Let $(C,P\sqcup Q_\bullet)$ be a smooth genus-zero curve with one marked point labeled by $P$ and a set of marked points labeled by $Q_\bullet=(Q_1,\ds, Q_m)$. Let $(W,M^1,\ds, M^m)$ be irreducible $V$-modules attached to $(P,Q_1,\ds, Q_m)$, respectively. Then
	\begin{equation}\label{coinvbottom}
		\V(V;W\o M^\bullet)_{(C,P\sqcup Q_\bullet, t\sqcup s_\bullet)}= [W\o M^1(0)\o\ds\o M^m(0)].
	\end{equation}
\end{proposition}
\begin{proof}
	We first show that
	$\V(V;W\o M^\bullet)_{(C,P\sqcup Q_\bullet,t\sqcup s_\bullet)}= [W\o M^1\o\ds\o M^{m-1}\o M^m(0)].$
	Since $M^m$ is irreducible, it is spanned by elements of the form $a_{(d)}v^m$, where $a\in V,d\in \Z$, and $v^m\in M^m(0)$ \cite{LL04}. Consider a general spanning element
	\[
	w\o u^1\o\ds\o u^{m-1}\o a_{(d)}v^m\in W\o M^\bullet,\quad w\in W,\ u^j\in M^j,\ v^m\in M^m(0).
	\]
	We may choose $N>0$ large enough so that, when $r\geq N$, we have $a_{(r)}u^j=0$ for all $1\leq j\leq m-1$ and $a_{(r)}v^m=0$.
	By Lemma~\ref{lmRR}, there exists $\si\in \mathcal{L}_{C\bs P\sqcup Q_\bullet}(V)$ such that
\[
\si_{Q_m}=a_{[d]}+\sum_{n\geq 0} \la_n a_{[n+N]}\in L_{Q_m}(V),\qquad \si_{Q_j}=\sum_{n\geq 0} \la'_n a_{[n+N]}\in L_{Q_j}(V),\quad \text{for all }1\leq j\leq m-1.
\]
	In particular, $\si_{Q_j}.u^j=0$ for all $1\leq j\leq m-1$ and $\si_{Q_m}.v^m=a_{(d)}v^m$.
	Then, by the definition of the chiral Lie algebra action \eqref{chiralLieact}, we have
	\[
	[w\o u^1\o\ds\o u^{m-1}\o a_{(d)}v^m]=-[(\si_P.w)\o u^1\o\ds \o u^{m-1}\o v^m].
	\]
	This shows that $\V(V;W\o M^\bullet)_{(C,P\sqcup Q_\bullet,t\sqcup s_\bullet)}= [W\o M^1\o\ds\o M^{m-1}\o M^m(0)].$ Next, we consider a general spanning element
	\[
	w\o u^1\o\ds \o u^{m-2}\o b_{(k)}v^{m-1}\o v^m\in W\o M^1\o\ds \o M^{m-1}\o M^m(0),
	\]
	where $v^{m-1}\in M^{m-1}(0)$. Using Lemma~\ref{lmRR} again, there exists a section $\tau\in  \mathcal{L}_{C\bs P\sqcup Q_\bullet}(V)$ such that
	$$\tau_{Q_j}.u^j=0\ \  (1\leq j\leq m-2);\quad \tau_{Q_{m}}.v^m=0;\quad \tau_{Q_{m-1}}.v^{m-1}=b_{(k)}v^{m-1},$$
	where we take $\tau_{Q_m}\in L_{Q_m}(V)_{<0}$ and $\tau_{Q_{m-1}}\equiv b_{[k]}\pmod{ L_{Q_{m-1}}(V)_{<0}}$ in \eqref{eq:RRlocalsigma}. Here, $L_{Q}(V)_{<0}$ is the Lie subalgebra of $L_Q(V)$ spanned by elements $a_{[n]}$ such that $\deg(a_{[n]})<0$.
	Then we have
	\[
	[w\o u^1\o\ds \o u^{m-2}\o b_{(k)}v^{m-1}\o v^m ]=-[(\tau_P.w)\o u^1\o\ds \o u^{m-2}\o v^{m-1}\o v^m].
	\]
	Hence, $[W\o M^1\o\ds\o M^{m-1}\o M^m(0)]=[W\o M^1\o\ds\o M^{m-2}\o M^{m-1}(0)\o M^m(0)]$. Proceeding inductively, we can show \eqref{coinvbottom}.
\end{proof}

We wish to further reduce \(W\) in \eqref{coinvbottom} to its bottom degree \(W(0)\) inside the space of coinvariants. However, this is a subtle problem in general, which motivates the following key notion for modules over an arbitrary VOA:
\begin{definition}\label{df:geored}
Let \(V\) be a VOA, and let $\mathscr W$ be the set of isomorphism classes of irreducible admissible \(V\)-modules. A subset of irreducible modules $\S\subseteq \mathscr{W}$ is called {\bf geometrically reducible} if, for every {\em smooth curve} \(C\simeq \P^1\) with $n\geq 3$ marked points $P_\bullet=(P_1,\dots,P_n)$, and for every choice of modules $W^1,\dots,W^n\in \S$ (possibly with repetitions), the canonical map to the quotient
	\begin{equation}\label{eq:surjfrombottom}
		W^\bullet(0)
		\to
		\mathbb V(V;W^\bullet)_{(C,P_\bullet,t_\bullet)},
		\qquad
		u^\bullet\mapsto [u^\bullet],
	\end{equation}
	is surjective. In other words, the space of coinvariants has the following reduction: $$\V(V;W^\bullet)_{(C,P_\bullet,t_\bullet)}=[W^\bullet(0)]=[W^1(0)\o\ds \o W^n(0)].$$
\end{definition}

\begin{example}
\label{exmp: GR and affine VOAs}
If $V=L_k(\g)$ is an affine VOA of positive integral level $k\in \Z_{>0}$, then the irreducible $V$-modules are the same as the integrable highest-weight $\hat{\g}$-modules $\mathcal{H}_{\la}$ of level $k$ \cite{Kac90,FZ92}.
It was proved in \cite{TUY89,Uen08} that the following canonical map is surjective:
$$L(\la_1)\o\ds\o L(\la_n)\twoheadrightarrow \mathbb V(V;\mathcal{H}_{\la_1}\o\ds \o \mathcal{H}_{\la_n})_{(C,P_\bullet,t_\bullet)},\quad v^\bullet\mapsto [v^\bullet],$$
where $\la_1,\ds, \la_n\in P^k_+$, and $\mathcal{H}_{\la}\cong L_{k}(\g,\la)$ as an $L_k(\g)$-module. The bottom degree of $L_{k}(\g,\la)$ is the $\g$-module $L(\la)$.
Therefore, for the affine VOA $L_k(\g)$, the full set of irreducible modules $\mathscr{W}$ is geometrically reducible. More generally, if $V$ is strongly generated by $V_1$, then it was proved in \cite[Section 5.1]{DG23} that
all irreducible admissible $V$-modules are geometrically reducible. See Example~\ref{exmp: GR examples}(1) for further discussion.
\end{example}

\begin{remark}
\label{rem: commutative zhu algebra}
We should point out that the quotient map \eqref{eq:surjfrombottom} is not necessarily surjective for an arbitrary choice of irreducible modules over a general VOA. For instance, let $V$ be a discrete series Virasoro VOA $L(c_{p,q},0)$ \cite{FZ92,W93,DMZ94} or a principal nilpotent $\mathcal{W}$-algebra $W_k(\g)$ at a non-degenerate admissible level $k$ \cite{A15}. Then $V$ is strongly rational, and every irreducible $V$-module $W$ has a one-dimensional bottom degree $W(0)$, since $A(V)$ is commutative. Hence, $W^\bullet(0)$ is one-dimensional as well. In particular, if the rank of the bundle $\V(V;W^\bullet)$ is greater than one for some choice of irreducible modules $W^\bullet$, then \eqref{eq:surjfrombottom} cannot be surjective; see also \cite[Section 8]{DG23}. 
\end{remark}

It is difficult to identify a geometrically reducible subset \(\mathcal{S}\) for a general VOA. We provide a criterion for identifying such a subset by combining the strong generation property of the VOA with the degrees of line bundles on the smooth curve \(C\).

\subsubsection{A criterion for geometric reducibility}
The criterion is given in Theorem~\ref{thm:criterion1}. Before stating the result, we need to discuss some notation and the following definition.
\begin{notation}
\begin{enumerate}
\item
Recall that $C_1(V)=\spn\{a_{(-1)}b: a,b\in V_+\}+L(-1)V_+$, where $V_+=\bigoplus_{n\geq 1} V_n$. Let $U\subset V$ be a subspace such that $V=U+ C_1(V)$. It was proved in \cite{KL99} that $V$ is {\em strongly generated} by $U$. In other words, a $C_1$-complementary subspace $U$ in $V$ is a space of {\em strong generators} of $V$. There is a similar relation between a $C_1$-complementary subspace and strong generation for $V$-modules $W$ \cite{H05,Liu25}.
\item
Recall that the degree-one subspace $V_1$ of a VOA $V$ is a Lie algebra with respect to
\[
[a,b]=a_{(0)}b,\quad a,b\in V_1.
\]
Let $\widehat{V_1}$ be the subspace of the universal enveloping algebra $\mathscr{U}(V)$ \cite{FZ92,FBZ04} spanned by elements $a_{(n)}$, where $a\in V_1$ and $n\in \Z$, and let $\mathscr{U}(\widehat{V_1})$ be the subalgebra of $\mathscr{U}(V)$ generated by $\widehat{V_1}$. Note that any admissible $V$-module $W$ is a module over the associative algebra $\mathscr{U}(\widehat{V_1})$, and  $$\mathscr{U}(\widehat{V_1})=\spn\{b^1_{(m_1)}\ds b^k_{(m_k)}: b^i\in V_1, m_i\in \Z\}.$$
\end{enumerate}
\end{notation}

\begin{definition}\label{df:stronglygen}
	Let $V=U+ C_1(V)$ for some homogeneous subspace $U$, and let $U_+:=U\cap V_+$. An irreducible $V$-module $M\in \mathscr{W}$ is said to be {\em strongly generated by $U_+$ from its bottom degree $M(0)$} if
	\begin{equation}\label{eq:Wspanning}
		M=\spn\{b^1_{(-n_1)}\ds b^r_{(-n_r)} u: b^i\in U_+,\ n_i\geq 0,\ \forall i,\  u\in M(0)\}.
	\end{equation}
\end{definition}

\begin{example}
\label{exmp: GR examples}
\begin{enumerate}
\item If $V$ is generated by $V_1$, then we can take $U=\C\vac+ V_1$, which satisfies $V=U+C_1(V)$. In this case, all irreducible $V$-modules $M\in\mathscr{W}$ are strongly generated by $U_+=V_1$ from their bottom degrees $M(0)$ \cite{FZ92,DG23}.
\item If $V=L(c_{p,q},0)$ is the discrete series Virasoro VOA \cite{FZ92,W93}, then we can take $U=\C\vac+ \C\om$, with $U_+=\C\om$. All irreducible $V$-modules $M\in \mathscr{W}$ are strongly generated by $\om$ from their bottom degrees $M(0)$.
\item However, if $V=V_L$ is the lattice VOA \cite{FLM88,DL93}, then, in general, not all irreducible $V$-modules are strongly generated from their bottom degrees. We will discuss this in detail in the next section.
\end{enumerate}
\end{example}

\begin{theorem}\label{thm:criterion1}
	Let $V=U+C_1(V)$ for some homogeneous subspace $U$ with $L(1)U=0$, and let $U_+:=U\cap V_+$.
	Let $\mathcal{S}\subset \mathscr{W}$ be a subset of irreducible modules such that each $M\in \mathcal{S}$ is strongly generated by $U_+$ from its bottom degree $M(0)$ and satisfies
	\begin{equation}\label{eq:assumptioncriterion1}
		b_{(n)} M(0)\ssq \mathscr{U}(\widehat{V_1}).M(0),\quad \forall b\in U\cap V_k,\ k\geq 1,\ n\geq 0.
	\end{equation}
	Then $\mathcal{S}$ is a geometrically reducible subset.
\end{theorem}
\begin{remark}
The proof given below depends on a series of lemmas. For ease of presentation, we have collected them in Section~\ref{subsec: lemmas}.
\end{remark}
\begin{proof}
	For any $m\geq 2$, let $W,M^1,\ds, M^m\in \mathcal{S}$, possibly with repetitions. By Proposition~\ref{prop:smoothdecent}, for any smooth curve $(C,P\sqcup Q_\bullet)$, we have $$\V(V;W\o M^\bullet)_{(C,P\sqcup Q_\bullet, t\sqcup s_\bullet)}= [W\o M^1(0)\o\ds\o M^m(0)].$$
	Let $W=\bigcup_{p=0}^\infty F_p(W)$, where $F_p(W)=W(0)+\ds +W(p)$. By Definition~\ref{df:geored}, it suffices to show that
	\begin{equation}\label{eq:ntsCriterion1}
		[F_p(W)\o M^1(0)\o\ds\o M^m(0)]=[W(0)\o M^1(0)\o\ds \o M^m(0)].
	\end{equation}
	We use induction on $p$ to show \eqref{eq:ntsCriterion1}. If $p=0$, this is clear. Let $p\geq 1$.
	By assumption, $W=\spn\{b^1_{(-n_1)}\ds b^r_{(-n_r)} u: b^i\in U_+,\ n_i\geq 0, u\in W(0) \}$ \eqref{eq:Wspanning}. Consider a general spanning element in $[F_p(W)\o M^1(0)\o\ds\o M^m(0)]$:
	\[
	\b=[b^1_{(-n_1)}\ds b^r_{(-n_r)} u\o u^1\o\ds \o u^m],\quad u\in W(0),\ u^i\in M^i(0),\ 1\leq i\leq m,\ b^i\in U_+,
	\]
	where $\deg(b^1_{(-n_1)}\ds b^r_{(-n_r)} u)=\deg(b^1_{(-n_1)})+\ds +\deg(b^r_{(-n_r)})\leq p$.

	Assume that $b^1\in U\cap V_k$. If $\deg(b^1_{(-n_1)})=k+n_1-1=0$, i.e., if $k=1$ and $n_1=0$, then, by adopting an argument similar to that in Proposition~\ref{prop:reductiononeptleaveout}, we can show that
	\[
	\b=\sum [b^2_{(-n_2)}\ds b^r_{(-n_r)}u\o v^1\o\ds\o v^m ],\quad u\in W(0),\ v^j\in M^j(0), \ 1\leq j\leq m.
	\]
	Moreover, in this case, $\deg (b^2_{(-n_2)}\ds b^r_{(-n_r)}u)=\deg(b^1_{(0)}b^2_{(-n_2)}\ds b^r_{(-n_r)}u)$. We may replace $\b$ by a summand $[b^2_{(-n_2)}\ds b^r_{(-n_r)}u\o v^1\o\ds\o v^m ]$.

	Therefore, we may assume that $\deg(b^1_{(-n_1)})>0$.
	Since $k\geq 1$ and $n_1\geq 0$, there exists a nonzero section $\mu\in H^0(C,\omega_C^{1-k}(n_1P))\bs H^0(C,\omega_C^{1-k}((n_1-1)P))$,
	which gives rise to
	$$b^1\o \mu \in V_k\o H^0(C, \om_C^{1-k}(n_1 P))\subset \bigoplus_{l\geq 0} V_l\o H^0(C\bs P\sqcup Q_\bullet, \om_C^{1-l}).$$
	Let $\si=b^1\o \mu +\mathrm{Im}\nabla\in  \mathcal{L}_{C\bs P\sqcup Q_\bullet}(V)$ \eqref{eq:descriptionofchiralLie}. Since $\mu$ has a pole of order exactly $n_1$ at $P$, we can write $\si_{P}=b^1_{(-n_1)}+\sum_{j\geq 1} \la_j b^1_{(-n_1+j)}$. Then

	\begin{equation}\label{eq:beta}
			\b=-\sum_{j\geq 1} \la_j[(b^1_{(-n_1+j)}\ds b^r_{(-n_r)} u)\o u^1\o \ds \o u^m]
			 \ - \sum_{i=1}^m [(b^2_{(-n_2)}\ds b^r_{(-n_r)} u)\o u^1\o\ds \o\si_{Q_i}u^i  \o \ds \o u^m ].
	\end{equation}
	Since $\deg (b^1_{(-n_1+j)}\ds b^r_{(-n_r)} u)=\deg(b^1_{(-n_1)}\ds b^r_{(-n_r)} u)-j<p$ for all $j\geq 1$, by the induction hypothesis, $[(b^1_{(-n_1+j)}\ds b^r_{(-n_r)} u)\o u^1\o \ds \o u^m]\in [W(0)\o M^1(0)\o\ds \o M^m(0)].$ On the other hand, note that $\deg (b^2_{(-n_2)}\ds b^r_{(-n_r)} u)<p$, as $\deg(b^1_{(-n_1)})>0$.
	Since $\mu$ is regular at each $Q_i$, we have $\si_{Q_i}=\sum_{n\geq 0} \gamma_n b^1_{(n)}$. By assumption \eqref{eq:assumptioncriterion1}, $\si_{Q_i}.u^i \in 	\mathscr{U}(\widehat{V_1}).M^i(0)$. Then, by Lemma~\ref{lm:UV1spanning}, we may write $\si_{Q_i}.u^i =\sum a^1_{(-m_1)}\ds a^s_{(-m_s)} v^i$, where $a^1,\ds, a^s\in V_1$, $m_1,\ds, m_s\geq 1$, and $v^i\in M^i(0)$.
	Hence,
	\begin{align*}
		&[(b^2_{(-n_2)}\ds b^r_{(-n_r)} u)\o u^1\o\ds \o\si_{Q_i}u^i  \o \ds \o u^m ]\\
		&=\sum [(b^2_{(-n_2)}\ds b^r_{(-n_r)} u)\o u^1\o\ds \o (a^1_{(-m_1)}\ds a^s_{(-m_s)} v^i)  \o \ds \o u^m ].
	\end{align*}

	Now proceed as in the proof of Lemma~\ref{lm:smoothredDG}. Take a nonzero section $$\eta=a^1\o \mu'\in V_1\o H^0(C,\O_C(m_1Q_i))\bs H^0(C,\O_C((m_1-1)Q_i)).$$ We have $\eta_{Q_j}.u^j\in M^j(0)$, since $a^1\in V_1$ and $\eta$ is regular at $Q_j$, for all $1\leq j\leq m$ with $j\neq i$. Moreover, $\eta$ is regular at $P$ as well, and
	\[
	\deg(a^1_{(n)}b^2_{(-n_2)}\ds b^r_{(-n_r)} u )\leq \deg(b^2_{(-n_2)}\ds b^r_{(-n_r)} u)<p,\quad \forall n\geq 0.
	\]
	Proceeding repeatedly, we have
	\begin{align*}
		& [(b^2_{(-n_2)}\ds b^r_{(-n_r)} u)\o u^1\o\ds \o (a^1_{(-m_1)}\ds a^s_{(-m_s)} v^i)  \o \ds \o u^m ]\\
		&=\sum [w^t\o v^1\o\ds \o v^i  \o \ds \o v^m ],
	\end{align*}
	where $w^t\in F_{p-1}(W)$ and $v^j\in M^j(0)$ for all $1\leq j\leq m$. Then, by the induction hypothesis, all the terms on the right-hand side of $\b$ in \eqref{eq:beta} are contained in $[W(0)\o M^1(0)\o\ds \o M^m(0)]$. This proves \eqref{eq:ntsCriterion1} for $p$ and completes the induction step.
\end{proof}

\begin{remark}
	Although condition \eqref{eq:assumptioncriterion1} may appear difficult to verify, it is in fact satisfied by a substantial class of VOAs. Indeed, if $V$ is generated by $V_1$ as in \cite{DG23}, then we can take $U_+=V_1$. It is clear that $b_{(n)}M(0)\ssq M(0)\ssq \mathscr{U}(\widehat{V_1}).M(0)$ for all $b\in V_1$ and $n\geq 0$.

	On the other hand, the module $M=V$ obviously satisfies condition \eqref{eq:assumptioncriterion1}, since $M(0)=\C\vac$ and $b_{(n)}\vac=0$ for all $n\geq 0$.
	In Section~\ref{sec:3.3}, we will show that certain modules over lattice VOAs also satisfy this condition.
\end{remark}

\begin{example}
	We first give an example of a VOA whose nontrivial modules do not satisfy condition \eqref{eq:assumptioncriterion1} for the standard choice of generators $U$. Then, in Sections~\ref{sec:3.3}--\ref{sec:3.5}, we give several examples of VOAs with a subset of irreducible modules satisfying condition \eqref{eq:assumptioncriterion1}.

	Let $V=L(\frac{1}{2},0)$ be the Ising model Virasoro VOA \cite{FZ92,W93}. Note that $V_1=0$ and that a canonical choice of generators $U$ is $U=\C\vac+ \C \om$, where $\om \in V_2$ is the Virasoro element. Moreover, $\mathscr{W}=\{V=L(\frac{1}{2},0), L(\frac{1}{2}, \frac{1}{2}), L(\frac{1}{2},\frac{1}{16})\}$.
	For $M=L(\frac{1}{2}, \frac{1}{2})$, we have $M(0)=\C v_{\frac{1}{2},\frac{1}{2}}$ and $L(-1)v_{\frac{1}{2},\frac{1}{2}}\neq 0$. Hence,
	\[
	\om_{(0)} M(0) =\C L(-1)v_{\frac{1}{2},\frac{1}{2}}\not\subseteq M(0),\quad \mathrm{but}\quad \mathscr{U}(\widehat{V_1}).M(0)=M(0).
	\]
	In other words, the irreducible modules over $L(\frac{1}{2},0)$ satisfy the strong generation property from their bottom degrees \eqref{eq:Wspanning}, but not property \eqref{eq:assumptioncriterion1}.
\end{example}

\subsubsection{Lemmas necessary to prove Theorem~\ref{thm:criterion1}}
\label{subsec: lemmas}

\begin{lemma}\label{lm:UV1spanning}
	Let $W$ be an admissible $V$-module, and let $N\subset W$ be a subspace such that $a_{(n)}N\ssq N$ for all $a\in V_1$ and $n\geq 0$. Let $	\mathscr{U}(\widehat{V_1}).N$ be the $\mathscr{U}(\widehat{V_1})$-submodule generated by $N$. Then we have
	\begin{equation}\label{eq:V1generationspann}
		\mathscr{U}(\widehat{V_1}).N=\spn\{a^1_{(-n_1)}\ds a^r_{(-n_r)} u: a^i\in V_1, n_i\geq 1,u\in N \}.
	\end{equation}
\end{lemma}
\begin{proof}
	By the definition of $\mathscr{U}(\widehat{V_1})$, we have
	\[
	\mathscr{U}(\widehat{V_1}).N=\spn\{b^1_{(m_1)}\ds b^k_{(m_k)} u: b^i\in V_1, m_i\in \Z, u\in N\}.
	\]
	Denote the space on the right-hand side of \eqref{eq:V1generationspann} by $S$. Then $S\ssq 	\mathscr{U}(\widehat{V_1}).N$. Conversely, we use induction on the length $k$ of $w=b^1_{(m_1)}\ds b^k_{(m_k)} u\in 	\mathscr{U}(\widehat{V_1}).N$ to show that $w\in S$.

	If $k=1$, then $b^1_{(m_1)}N\ssq N$ for $m_1\geq 0$ by assumption. Hence, the nonzero length-one elements in $\mathscr{U}(\widehat{V_1}).N$ are of the form $b^1_{(-m_1)}u$, where $b^1\in V_1$, $u\in N$, and $m_1\geq 1$. This proves the base case.
	Now let $k\geq 2$. If $m_1<0$, then, by the induction hypothesis, we may assume
	\[
	b^2_{(m_2)}\ds b^k_{(m_k)} u=\sum c^1_{(-l_1)}\ds c^s_{(-l_s)}v^t\in S,\quad l_i\geq 1,\ v^t\in N.
	\]
	Hence, $w=\sum  b^1_{(m_1)}c^1_{(-l_1)}\ds c^s_{(-l_s)}v^t\in S$. Now let $m_1\geq 0$. Then, by assumption, $b^1_{(m_1)}u\in N$. Furthermore, it follows from the induction hypothesis that
	\[
	w=b^2_{(m_2)}\ds b^k_{(m_k)} (b^1_{(m_1)}u)+\sum_{j=2}^k b^2_{(m_2)}\ds \left(([b^1,b^j])_{(m_1+m_j)}+m_1(b^1_{(1)}b^j)_{(m_1+m_j-1)}\right) \ds b^k_{(m_k)}u\in S,
	\]
	since $[b^1,b^j]\in V_1$ and $b^1_{(1)}b^j\in V_0=\C\vac$. Therefore, $	\mathscr{U}(\widehat{V_1}).N=S$.
\end{proof}

The following lemma is essentially the main result of \cite[Section 5.1]{DG23}; we skip the details.

\begin{lemma}\label{lm:smoothredDG}
	Let $W^1,\ds, W^n\in \mathscr{W}$. For each $1\leq i\leq n$, let $N^i\subset W^i$ be a subspace such that $a_{(m)}N^i\ssq N^i$ for all $a\in V_1$ and $m\geq 0$. With the notation \eqref{eq:classincoinv}, we have the following reduction property in the space of coinvariants $\V(V;W^\bullet)_{(C,P_\bullet,t_\bullet)}$:
	\begin{equation}
		[	\mathscr{U}(\widehat{V_1}).N^1\otimes \ds \o 	\mathscr{U}(\widehat{V_1}).N^n]=[N^1\o\ds \o N^n].
	\end{equation}
	In particular, $N^i=W^i(0)$ satisfies the property that $a_{(m)}W^i(0)\ssq W^i(0)$ for all $a\in V_1$ and $m\geq 0$. Then $[	\mathscr{U}(\widehat{V_1}).W^1(0)\otimes \ds \o 	\mathscr{U}(\widehat{V_1}).W^n(0)]=[W^\bullet(0)]$.
\end{lemma}

The following proposition generalizes Lemma~\ref{lm:smoothredDG} and is essential in the proof of Theorem~\ref{thm:criterion1}.
\begin{proposition}\label{prop:reductiononeptleaveout}
	Let $(C,P_\bullet\sqcup Q)$ be a smooth curve. Let $W^1,\ds, W^n\in \mathscr{W}$ and let $M$ be an admissible $V$-module, attached at $P_\bullet$ and $Q$, respectively. Assume that $n\geq 2$.

	Fix an element $\b\in \Om(M)$, i.e., $a_{(m)}\b=0$ for all $a\in V$ and $m\in \Z$ such that $\deg (a_{(m)})=\wt a-m-1<0$. Then, in the space of coinvariants $\V(V;W^\bullet\o M)_{(C,P_\bullet\sqcup Q,t_\bullet\sqcup s)}$, we have
	\begin{equation}\label{eq:reductiononeaddtionalpoint}
		[\scrU(\widehat{V_1}).W^1(0)\o\ds \o \scrU(\widehat{V_1}).W^n(0)\o \b]=[W^1(0)\o\ds\o W^n(0)\o \b],
	\end{equation}
	where $[N^\bullet\o \b]$ is defined by $\spn\{[v^\bullet\o \b]: v^\bullet\in N^\bullet\}$.
\end{proposition}
\begin{proof}
	As in Lemma~\ref{lm:smoothredDG}, we consider a spanning element on the left-hand side of \eqref{eq:reductiononeaddtionalpoint}:
	\[
	B=[(a^1_{(-n_1)}\ds a^r_{(-n_r)} u)\o w_2\o \ds \o w_n\o \b ],\quad u\in W^1(0),\ w_i\in \mathscr{U}(\widehat{V_1}).W^i(0),\ 2\leq i\leq n.
	\]
	Since $h^0(C,\om_C^{1-1}(n_1P_1-Q))=n_1$ and $h^0(C,\om_C^{1-1}((n_1-1)P_1-Q))=n_1-1$, there exists a nonzero section $\mu'\in H^0(C,\om_C^{1-1}(n_1P_1-Q))\bs H^0(C,\om_C^{1-1}((n_1-1)P_1-Q))$, which has a pole of order exactly $n_1$ at $P_1$. Now $\mu'$ gives rise to an element
	$$\tau=a^1\o\mu'+\mathrm{Im}\nabla\in \left(V_1\o H^0(C,\om^{1-1}_C(n_1P_1-Q))+\mathrm{Im}\nabla\right)/\mathrm{Im}\nabla\ssq \mathcal{L}_{C\bs P_\bullet\sqcup Q}(V),$$
	such that $\tau_Q=\sum_{m\geq 1} \la_m\cdot  a^1_{(m)}$, which annihilates $\b\in \Om(M)$ since $\deg(a^1_{(m)})=1-m-1<0$ for all $m\geq 1$; $\tau_{P_i}w_i$ is a sum of elements in $W^i$ of degree less than or equal to $\deg w_i$ for all $2\leq i\leq n$; and $\tau_{P_1}=a^1_{(-n_1)}+\sum_{k\geq 1} \la_k a^1_{(-n_1+k)}$. Thus,
	\begin{align*}
		B&=-\sum_{k\geq 1} \la_k[(a^1_{(-n_1+k)}\ds a^r_{(-n_r)} u)\o w_2\o \ds \o w_n\o \b]\\
		&\ \ - \sum_{i=2}^n [(a^2_{(-n_2)}\ds a^r_{(-n_r)} u)\o w_2\o\ds \o\tau_{P_i}w_i  \o \ds \o w_n \o \b],
	\end{align*}
	where the right-hand side is a sum of elements of degree strictly less than $\deg B$. Now \eqref{eq:reductiononeaddtionalpoint} follows by induction on the degree of $B$.
\end{proof}

\subsection{Nodal case}
\label{subsec: nodal case}
In this section, we denote the space of coinvariants $\V(V;M^\bullet)_{(C,P_\bullet,t_\bullet)}$ by $[M^\bullet]_{\mathcal{L}_{C\bs P_\bullet}(V)}$ for simplicity of notation, where $C$ is a smooth or singular curve of genus zero. Throughout this section, we let $V=U+C_1(V)$ for some homogeneous subspace $U$ with $L(1)U=0$, and let $U_+=U\cap V_+$.
The following lemma is essentially proved in \cite{Liu25}:
\begin{lemma}
	The Zhu algebra $\A=A(V)=V/O(V)$ \cite{Z96} is generated by the classes $[a]=a+O(V)$ with $a\in U$. In other words,
	\begin{equation}\label{eq:Aspanning}
		\A=\spn\{[a^1]\cdot[a^2] \ds \cdot [a^r]: r\geq 0,\ a^i\in U_+,\ \forall i \}.
	\end{equation}
\end{lemma}
\begin{proof}
	Denote the image of $U\subset V$ in $\A=V/O(V)$ by $[U]$.
	Recall that $\A$ is a filtered algebra $\A=\bigcup_{p=0}^\infty F_p\A$, where $F_p\A$ is the image of $\bigoplus_{n=0}^p V_n\subset V$ in $\A=V/O(V)$. Since $V=U+C_1(V)$, by the proof of \cite[Theorem 3.1]{Liu25}, the associated graded algebra $\mathrm{gr}(\A)$ is commutative and is generated by $\{[a]+F_{p-1}\A\in \mathrm{gr}(\A): [a]\in F_p\A\cap [U]\} $. Therefore, $\A$ itself is generated by $[U]$, with identity $1_\A=[\vac]$.
\end{proof}

\subsubsection{Coinvariants of modules that are closed under fusion}\label{sec:4.1}
\begin{definition}\label{def:Sclosedunderfusion}
Let $\S\ssq \mathscr{W}$ be a subset of irreducible $V$-modules. We say that $\S$ is {\em closed under the fusion tensor product} if 
\begin{equation}
	M^i\bt_V M^j=\bigoplus_{ W\in \S} N_{M^iM^j}^{W} W,\quad \forall M^i,M^j\in \S.
\end{equation}
Equivalently, for any $M^i,M^j\in \S$ and $W\notin \S$, the fusion rule $N_{M^iM^j}^{W}=0$; in other words, the free abelian group generated by $\S$ is a subring of the fusion ring $\mathscr{F}(V)=\bigoplus_{W\in \mathscr{W}}\Z [W]$ in the sense of Definition~\ref{def:fusionsubring}; note that $V\in\S$ if $\S\neq\emptyset$, by Lemma~\ref{lm:rigidity}.
\end{definition}

\begin{lemma}\label{lm:subring}
	Let $(C,P\sqcup Q_\bullet)$ be a smooth genus-zero curve, and let $(W,M^1,\ds, M^m)$ be irreducible $V$-modules attached to $(P,Q_1,\ds, Q_m)$, respectively. Assume that $\S$ is closed under the fusion tensor product, $M^\bullet\subset \S$, and $W\notin \S$. Then the space of coinvariants
	\[
	\V(V;W\o M^\bullet)_{(C,P\sqcup Q_\bullet, t\sqcup s_\bullet)}=[W\o M^\bullet]_{\mathcal{L}_{C\bs P\sqcup Q_\bullet}(V)}=0.
	\]
\end{lemma}
\begin{proof}
This follows from the proof of \cite[Lemma 3.4]{ChakParafermion} and Lemma~\ref{lm:rigidity}, using the factorization and vector-bundle properties of the sheaves of coinvariants \cite{DGT24}.
\end{proof}

\begin{lemma}\label{lm:rigidity}
	If $\S\ssq \mathscr{W}$ is closed under taking fusion tensor products, then it is also closed under taking contragredient modules, i.e., if $W\in \S$, then $W'\in \S$. In particular, $V\in \S$ if $\S\neq\emptyset$.
\end{lemma}
\begin{proof}
For $W\in \mathscr{W}$, let $N_W$ be the matrix of the multiplication by $[W]$ on $\mathscr{F}(V)$ with respect to the basis $\{[W^0],\ds,[W^N]\}$. By the symmetry of fusion rules $N_{WM}^{K}=N_{W'K}^{M}$ \cite{FHL93,H08}, we have $N_{W'}=N_W^{T}$. Since $\mathscr{F}(V)$ is commutative, $\{N_W\}_{W\in \mathscr{W}}$ is a commuting family of real normal matrices. Hence they are simultaneously diagonalizable, and $\mathscr{F}(V)\o_\Z\C\cong \C^{N+1}$ as algebras, under which $[W']$ corresponds to the complex conjugate of $[W]$. Let $R$ be the $\C$-span of $\S$ in $\mathscr{F}(V)\o_\Z \C$. Since $\S$ is closed under the fusion tensor product, $R$ is a subalgebra of $\C^{N+1}$ without nonzero nilpotent elements, and so $R$ is spanned by idempotents, which are vectors with entries in $\{0,1\}$. Hence $R$ is closed under complex conjugation, and $[W']\in R$ for all $W\in \S$. Since $R$ is spanned by the basis elements $[M]$ with $M\in \S$, we have $W'\in \S$. Finally, $V\subset W\bt_V W'$ implies $V\in \S$.
\end{proof}

\subsubsection{Mode transition algebra}
We briefly recall the definition of the {\em mode transition algebra} from \cite{DGK25}. Recall that $\A=A(V)$ denotes the Zhu algebra \cite{Z96} of $V$. Let $T$ be a left $\A$-module, and let $\mathscr{U}=\mathscr{U}(V)$ be the universal enveloping algebra of $V$, which is a graded topological associative algebra \cite{FZ92,FBZ04}. The {\em generalized Verma module associated with $T$} is defined by
\[
\Phi^\L(T)=(\mathscr{U}/\mathsf{N}^1_\L\mathscr{U})\otimes_{\mathscr{U}_0}T,
\]
where $\mathsf{N}^1_\L\mathscr{U}=\mathscr{U}\cdot \mathscr{U}_{\leq -1}\subset \mathscr{U}$ is the first fundamental neighborhood.
The space $\Phi^\L(T)$ is an admissible $V$-module such that $\Phi^\L(T)(0)=T$ \cite{DLM98,DGK25}. Let $T=\A$ be the left regular $\A$-module. Since $V$ is a strongly rational VOA, $\A$ is semisimple and
\begin{equation}\label{eq:isoofA}
	\A\cong  \prod_{W\in \mathscr{W}} \End_\C(W(0))\cong  \bigoplus_{W\in \mathscr{W}} W(0)\otimes W(0)^\ast
\end{equation}
as an associative algebra and as a bimodule over itself.
Moreover, for any $W\in \mathscr{W}$, it follows from the rationality of $V$ that $W\cong \Phi^\L(W(0))$ \cite{DLM98}.
Therefore, as a $(V,\A)$-bimodule,
\begin{equation}\label{eq:PhiLA}
	\Phi^\L(\A)\cong \bigoplus_{W\in \mathscr{W}} W\o W(0)^\ast.
\end{equation}
Under the identification \eqref{eq:isoofA}, the identity $1_\A$ of $\A$ can be written as
\[
1_\A=\sum_{W\in  \mathscr{W}} \Id_{W(0)}=\sum_{W\in \mathscr{W}} \sum_i e(W)_i\o e(W)^i,
\]
where $\{e(W)_i\}$ and $\{e(W)^i\}$ are dual bases of $W(0)$ and $W(0)^\ast$, respectively. The element $1_\A$ can be viewed as an element of $\Phi^\L(\A)$, since the bottom degree is $\Phi^\L(\A)(0)=\A$.

Let $R$ be a right $\A$-module. One can similarly define the right induced module $\Phi^\sR(R)=R\o_{\scrU_0}(\scrU/\mathsf{N}^1_\sR\scrU)$, which is naturally a right $\scrU$-module. Since there is an anti-involution of associative algebras $\theta: \A\ra \A, [a]\mapsto [e^{L(1)}(-1)^{L(0)}a]$, which generalizes to an anti-involution $\theta: \scrU\ra \scrU$, we write $\Phi^\sR(R)^\theta$ for the left $\scrU$-module structure induced by $\theta$, i.e.,
\[
\scrU\times \Phi^\sR(R)^\theta\ra \Phi^\sR(R)^\theta,\quad a_{[n]}.u:=u\theta(a_{[n]}),\quad \forall a_{[n]}\in L(V).
\]
Note that $L(V)$ generates $\scrU$ as a topological algebra.
Then $\Phi^\sR(R)^\theta$ is an admissible $V$-module whose bottom degree is the left $\A$-module $R^\theta$. In particular, for any $M\in \mathscr{W}$, the dual space $M(0)^\ast$ is naturally a right $\A$-module, and $(M(0)^\ast)^\theta$ is isomorphic to the bottom degree of the contragredient module $M'$. We have $\Phi^\sR(M(0)^\ast)^\theta\cong M'$, and
\begin{equation}\label{eq:PhiRA}
	\Phi^\sR(\A)\cong \bigoplus_{M\in \mathscr{W}} M(0)\o \Phi^\sR(M(0)^\ast), \quad \Phi^\sR(\A)^\theta\cong \bigoplus_{M\in \mathscr{W}} M(0)\o M'.
\end{equation}

The {\em mode transition algebra} of the VOA $V$ is the induced $\scrU-\scrU$-bimodule \cite{DGK25}:
\[
\mathfrak{A}=\Phi(\A)=\Phi^\L\Phi^\sR(\A)\cong \Phi^\L(\A)\o_\A \Phi^\sR(\A).
\]

\subsubsection{Mode transition algebra associated with an ideal of $\A$} For a subset $\S\ssq \mathscr{W}$ that is closed under the fusion product, we set
\begin{equation}\label{eq:decofAS}
	\A_\S:=\prod_{W\in \S} \End_\C(W(0))\cong \bigoplus_{W\in \S} W(0)\o W(0)^\ast.
\end{equation}
Then $\A_\S\subset \A$ is a direct factor of $\A$ as an algebra, and a two-sided ideal of $\A$, with $1_{\A_\S}=\sum_{W\in \S} \sum_i e(W)_i\o e(W)^i$.
As in \eqref{eq:PhiLA} and \eqref{eq:PhiRA}, the decomposition \eqref{eq:decofAS} gives
\begin{equation}\label{eq:MTAforAS}
	\Phi^\L(\A_\S)\cong \bigoplus_{W\in \S} W\o W(0)^\ast,\quad \Phi^\sR(\A_\S)^\theta\cong \bigoplus_{M\in \S} M(0)\o M'
\end{equation}
as $V$-modules, where $M'\in \S$ by Lemma~\ref{lm:rigidity}.

Note that each direct summand $W\o W(0)^\ast$ of $\Phi^\L(\A_\S)$ is a right $\A$-module, while each summand $M(0)\o M'$ of $\Phi^\sR(\A_\S)^\theta$ is a left $\A$-module. By Schur's lemma, we have $W(0)^\ast \o_\A M(0)\cong  \delta_{W,M}\C=0$ if $W\in \S$ and $M\notin \S$. Therefore,
\begin{equation}\label{eq:AStensorA}
	\Phi^\L(\A_\S)\o _\A \Phi^\sR(\A)\cong \Phi^\L(\A_\S)\o _\A \Phi^\sR(\A_\S).
\end{equation}
Furthermore, the bottom degrees are given by
\[
\Phi^\L(\A_\S)(0)\cong \bigoplus_{W\in \S} W(0)\o W(0)^\ast\cong \A_\S\cong \bigoplus_{M\in \S} M(0)\o M(0)^\ast \cong \Phi^\sR(\A_\S)^\theta(0).
\]
Observe that there are two left actions of $\A$ on the bottom degree $\A_\S=\Phi^\sR(\A_\S)^\theta(0)$:
\begin{enumerate}
	\item The first action is the natural auxiliary left $\A$-multiplication action, induced by the ideal structure of $\A_\S\subset \A$ and the fact that $\Phi^\sR(\A_\S)^\theta=\A_\S\o _{\scrU_0} (\scrU/\rN^1_\sR\scrU)$ as a vector space. We denote it by $[a]\ast s$, for $[a]\in \A$ and $s\in \A_\S$.
	\item The second action is given by the usual Zhu algebra $\A$-action on the bottom degree of an admissible $V$-module, which, in this case, is defined by $\theta$ and the right $\scrU_0$-action on $\A_\S$. We denote it by $[a].s=s\o \theta(a_{(\wt a-1)})$, for all $[a]\in \A$ and $s\in \A_\S$.
\end{enumerate}

\subsubsection{Coinvariants for nodal curves}

\begin{figure}[htbp]
	\centering
	\begin{tikzpicture}[thick, line cap=round, line join=round, scale=0.8]

		\begin{scope}[xshift=0cm]
			\coordinate (A) at (-1,  3);
			\coordinate (B) at ( 2, -3);
			\coordinate (C) at ( 1,  3);
			\coordinate (D) at (-2, -3);

			\draw (A) -- (B);
			\draw (C) -- (D);

			\fill (0,1) circle (2pt);
			\node[left] at (0.8,1) {$\qq$};

			\coordinate (Pw1)   at ($(C)!0.45!(D)$);
			\coordinate (Pw2)   at ($(C)!0.55!(D)$);
			\coordinate (Pw3)   at ($(C)!0.78!(D)$);
			\coordinate (PinfL) at ($(C)!0.9!(D)$);
			\fill (Pw1)   circle (2pt); \node[left]  at (Pw1)   {$P_m$};
			\fill (Pw2)   circle (2pt); \node[left]  at (Pw2)   {$P_{m-1} $};
			\fill (Pw3)   circle (2pt); \node[left]  at (Pw3)   {$P_{2}$};
			\fill (PinfL) circle (2pt); \node[left]  at (PinfL) {$P_1$};

			\coordinate (Pw1) at ($(A)!0.45!(B)$);
			\coordinate (Pw2) at ($(A)!0.55!(B)$);
			\coordinate (Pw3) at ($(A)!0.78!(B)$);
			\coordinate (P0)  at ($(A)!0.9!(B)$);
			\fill (Pw1) circle (2pt); \node[right] at (Pw1) {$Q_n$};
			\fill (Pw2) circle (2pt); \node[right] at (Pw2) {$Q_{n-1}$};
			\fill (Pw3)   circle (2pt); \node[right]  at (Pw3)   {$Q_{2}$};
			\fill (P0)  circle (2pt); \node[right] at (P0)  {$Q_1$};

			\node at (0,-3.6) { $(X,P_\bullet\sqcup Q_\bullet)$};

		\end{scope}

		\draw[->, thick]
		(4.2,0) -- (2.5,0)
		node[midway, above] {$\eta$};

		\begin{scope}[xshift=7cm]

			\coordinate (Ltop) at (-0.5,  3);
			\coordinate (Lbot) at (-1.5, -3);

			\coordinate (Rtop) at ( 1,  3);
			\coordinate (Rbot) at ( 2, -3);

			\draw (Lbot) -- (Ltop);
			\draw (Rtop) -- (Rbot);

			\coordinate (L0)    at ($(Ltop)!0.32!(Lbot)$);
			\coordinate (LQ1)   at ($(Ltop)!0.45!(Lbot)$);
			\coordinate (LQ2)   at ($(Ltop)!0.55!(Lbot)$);
			\coordinate (LQ3)   at ($(Ltop)!0.8!(Lbot)$);
			\coordinate (Linf)  at ($(Ltop)!0.9!(Lbot)$);
			\fill (L0)   circle (2pt); \node[left] at (L0)   {$\qq_+$};
			\fill (LQ1)  circle (2pt); \node[left] at (LQ1)  {$P_m$};
			\fill (LQ2)  circle (2pt); \node[left] at (LQ2)  {$P_{m-1}$};
			\fill (LQ3)  circle (2pt); \node[left] at (LQ3)  {$P_{2}$};
			\fill (Linf) circle (2pt); \node[left] at (Linf) {$P_1$};

			\coordinate (Rinf) at ($(Rtop)!0.32!(Rbot)$);
			\coordinate (RR1)  at ($(Rtop)!0.45!(Rbot)$);
			\coordinate (RR2)  at ($(Rtop)!0.55!(Rbot)$);
			\coordinate (RR3)  at ($(Rtop)!0.8!(Rbot)$);
			\coordinate (R0)   at ($(Rtop)!0.9!(Rbot)$);
			\fill (Rinf) circle (2pt); \node[right] at (Rinf) {$\qq_-$};
			\fill (RR1)  circle (2pt); \node[right] at (RR1)  {$Q_n$};
			\fill (RR2)  circle (2pt); \node[right] at (RR2)  {$Q_{n-1}$};
			\fill (RR3)  circle (2pt); \node[right] at (RR3)  {$Q_{2}$};
			\fill (R0)   circle (2pt); \node[right] at (R0)   {$Q_1$};

			\node at (0,-3.6) {$(\tilde{X}=X_+\sqcup X_-,P_\bullet\sqcup Q_\bullet \sqcup \qq_\pm)$};

		\end{scope}
	\end{tikzpicture}
	\caption{Genus-zero curve with one node}
	\label{fig1}
\end{figure}
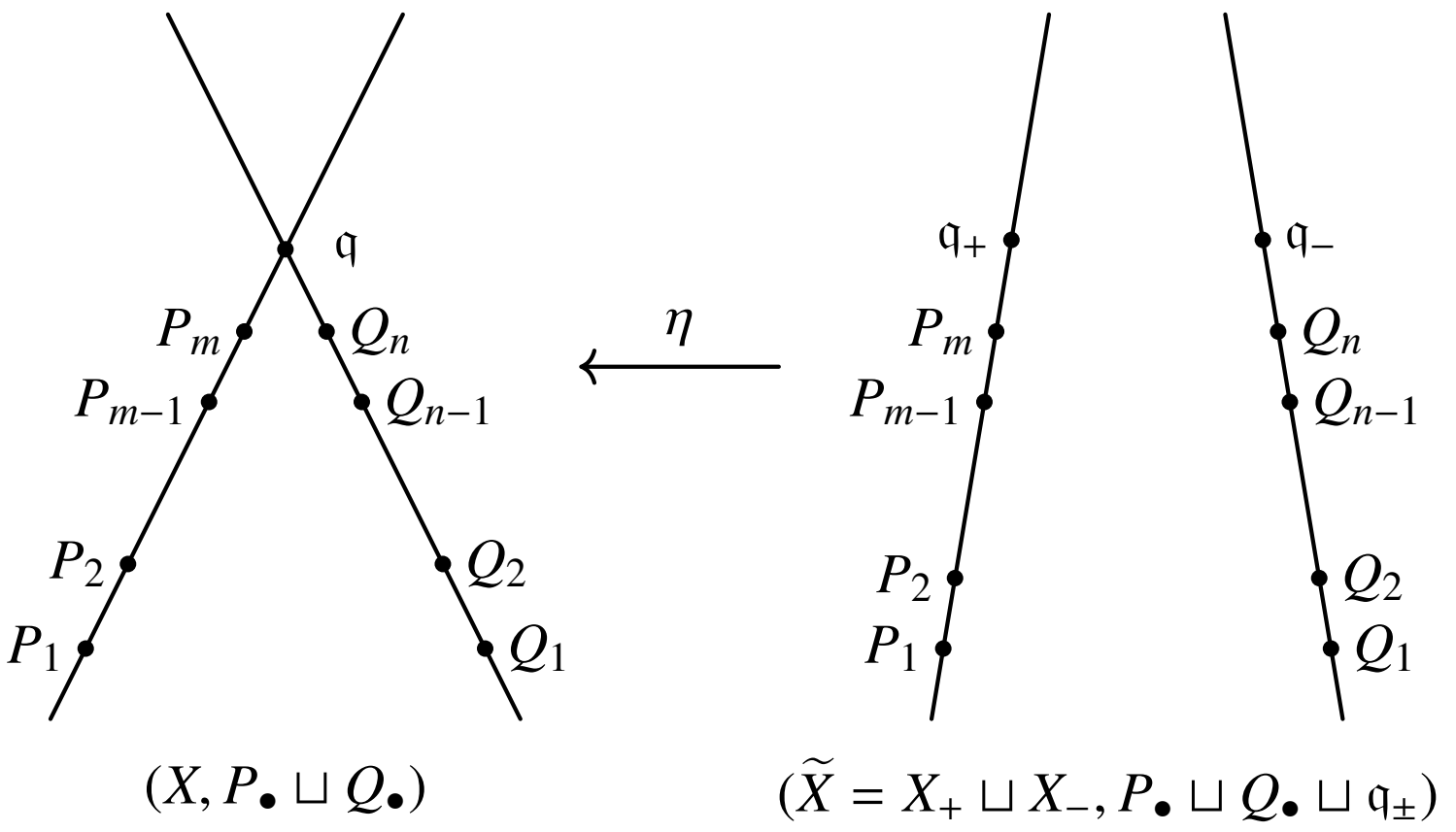

We first recall the definition of the chiral Lie algebra associated with a nodal curve \cite{DG23,DGT24}.
Let $(X,P_\bullet\sqcup Q_\bullet)$ be a genus-zero nodal curve with two sets of smooth marked points labeled by $P_\bullet=(P_1,\ds, P_m)$ and $Q_\bullet=(Q_1,\ds, Q_{n})$, and one simple node $\qq$. Let $(\widetilde{X},P_\bullet\sqcup Q_\bullet\sqcup \qq_\pm)$ be the resolution of the singularity at the node $\qq$, as illustrated in Figure~\ref{fig1}.

Note that both $X_+$ and $X_-$ are smooth genus-zero curves.
For simplicity, denote $\mathcal{L}_+=\mathcal{L}_{X_+\bs P_\bullet\sqcup \qq_+}(V)$ and $\mathcal{L}_-=\mathcal{L}_{X_-\bs Q_\bullet\sqcup \qq_-}(V)$.
Then the chiral Lie algebra for $(\tilde{X},P_\bullet\sqcup Q_\bullet \sqcup \qq_\pm)$ decomposes into subalgebras:
\[
\mathcal{L}_{\tilde{X}\bs P_\bullet\sqcup Q_\bullet \sqcup \qq_\pm}(V)\cong \mathcal{L}_+\op \mathcal{L}_-.
\]
It follows from \cite[Proposition 3.3.1]{DGT24} that the chiral Lie algebra associated with the nodal curve $(X,P_\bullet\sqcup Q_\bullet)$ can be described as follows:
\begin{equation}\label{eq:nodalLie}
	\mathcal{L}_{X\bs P_\bullet\sqcup Q_\bullet}(V)=\left\{
	\si=(\si_+,\si_-)\in \mathcal{L}_+\op\mathcal{L}_-: \begin{aligned}
		&(\si_+)_{\qq_+},(\si_{-})_{\qq_-}\in L_{\qq_\pm}(V)_{\leq 0}\\
		& [(\si_+)_{\qq_+}]_0=\theta([(\si_{-})_{\qq_-}]_0)
	\end{aligned}
	\right\}.
\end{equation}

Furthermore, let $M^\bullet,W^\bullet\subset \mathscr{W}$ be tuples of irreducible $V$-modules attached at $P_\bullet$ and $Q_\bullet$, respectively.
Then
the space of nodal coinvariants has the following identification; see \cite[eq.(40)]{DGT24} and \cite[Lemma 4.4.4]{DGK25}:

\begin{equation}\label{eq:onenodecoinviso1}
	\begin{aligned}
		h: 	[M^\bullet\o W^\bullet]_{\mathcal{L}_{X\bs P_\bullet\sqcup Q_\bullet}(V)}&\cong [M^\bullet\o \Phi(\A)\o W^\bullet]_{\mathcal{L}_{\tilde{X}\bs P_\bullet\sqcup Q_\bullet \sqcup \qq_\pm}(V)}\\
		& \cong [M^\bullet \o \Phi^\L(\A)]_{\mathcal{L}_{X_+\bs P_\bullet\sqcup \qq_+}(V)}\o_\A [\Phi^\sR(\A)^\theta\o W^\bullet]_{\mathcal{L}_{X_-\bs  \qq_-\sqcup Q_\bullet}(V)},\\
		[u^\bullet\o w^\bullet]&\mapsto [u^\bullet\o 1_\A]\o_\A [1_\A\o w^\bullet].
	\end{aligned}
\end{equation}
The space of coinvariants $[\Phi^\sR(\A)^\theta\o W^\bullet]_{\mathcal{L}_{X_-\bs \qq_-\sqcup Q_\bullet}(V)}$ is viewed as $\bigoplus_{M\in \mathscr{W}}m_M [M'\o W^\bullet]_{\mathcal{L}_{X_-\bs \qq_-\sqcup Q_\bullet}(V)}$, where $m_M=\dim M(0)$ and $M'$ is the contragredient module of $M$; see \eqref{eq:PhiRA}.
The identification \eqref{eq:onenodecoinviso1} also holds when $\qq$ is one of several nodes of a genus-zero curve, with $X_\pm$ replaced by the two connected components of the normalization at $\qq$. Indeed, after adding a marked point with the module $V$ on each irreducible component without smooth marked points, which does not change the coinvariants by Theorem~\ref{thm: POV}, this follows from \cite[Lemma 4.4.4]{DGK25} as in the one-node case.

\begin{lemma}\label{lm:PhiAandAS}
	Let $(X,P_\bullet\sqcup Q_\bullet)$ be a nodal curve with one simple node $\qq$, with normalization $(\widetilde{X},P_\bullet\sqcup Q_\bullet\sqcup \qq_\pm)$; see Figure~\ref{fig1}.
	Let $\S\ssq \mathscr{W}$ be a subset of modules that is closed under the fusion tensor product, and let $M^\bullet, W^\bullet\subset \S$. Then
	\begin{align*}
		[M^\bullet \o \Phi^\L(\A)]_{\mathcal{L}_{X_+\bs P_\bullet\sqcup \qq_+}(V)}&\cong [M^\bullet \o \Phi^\L(\A_\S)]_{\mathcal{L}_{X_+\bs P_\bullet\sqcup \qq_+}(V)},\\
		[\Phi^\sR(\A)^\theta\o W^\bullet]_{\mathcal{L}_{X_-\bs Q_\bullet\sqcup \qq_-}(V)}&\cong [\Phi^\sR(\A_\S)^\theta\o W^\bullet]_{\mathcal{L}_{X_-\bs Q_\bullet\sqcup \qq_-}(V)},
	\end{align*}
	where $\Phi^\L(\A_\S)$ and $\Phi^\sR(\A_\S)^\theta$ are given by \eqref{eq:MTAforAS}, and the isomorphisms are right and left $\A$-module isomorphisms, respectively.
	Moreover, the isomorphism \eqref{eq:onenodecoinviso1} becomes
	\begin{equation}\label{eq:4.10}
		\begin{aligned}
			[M^\bullet\o W^\bullet]_{\mathcal{L}_{X\bs P_\bullet\sqcup Q_\bullet}(V)}
			& \cong [M^\bullet \o \Phi^\L(\A_\S)]_{\mathcal{L}_{+}}\o_\A [\Phi^\sR(\A_\S)^\theta\o W^\bullet]_{\mathcal{L}_{-}},\\
			[u^\bullet\o w^\bullet]&\mapsto [u^\bullet\o 1_{\A_\S}]\o_\A [1_{\A_\S}\o w^\bullet].
		\end{aligned}
	\end{equation}
\end{lemma}
\begin{proof}
	By \eqref{eq:PhiLA} and Lemma~\ref{lm:subring}, we have an identification
	\begin{align*}
		[M^\bullet \o \Phi^\L(\A)]_{\mathcal{L}_{X_+\bs P_\bullet\sqcup \qq_+}(V)}&=\bigoplus_{W\in \mathscr{W}} [M^\bullet\o W]_{\mathcal{L}_{X_+\bs P_\bullet\sqcup \qq_+}(V)}\o W(0)^\ast\\
		&=\bigoplus_{W\in \S} [M^\bullet\o W]_{\mathcal{L}_{X_+\bs P_\bullet\sqcup \qq_+}(V)}\o W(0)^\ast\\
		&\cong  [M^\bullet \o \Phi^\L(\A_\S)]_{\mathcal{L}_{X_+\bs P_\bullet\sqcup \qq_+}(V)},
	\end{align*}
	under which $[u^\bullet \o 1_\A]$ corresponds to $[u^\bullet\o 1_{\A_\S}]$.
	By Lemma~\ref{lm:rigidity}, $\S\ra \S,\ M\mapsto M'$ is a one-to-one correspondence. It follows from \eqref{eq:PhiRA} and \eqref{eq:MTAforAS} that
	\begin{align*}
		[\Phi^\sR(\A)^\theta\o W^\bullet]_{\mathcal{L}_{X_-\bs Q_\bullet\sqcup \qq_-}(V)}&\cong \bigoplus_{M'\in \S} M(0)\o[M'\o W^\bullet]_{\mathcal{L}_{X_-\bs Q_\bullet\sqcup \qq_-}(V)}=[\Phi^\sR(\A_\S)^\theta\o W^\bullet]_{\mathcal{L}_{X_-\bs Q_\bullet\sqcup \qq_-}(V)},
	\end{align*}
	under which $[1_{\A}\o w^\bullet]$ corresponds to $[1_{\A_\S}\o w^\bullet]$. Now \eqref{eq:4.10} is clear.
\end{proof}

\subsubsection{Reduction of the nodal coinvariants} We show that, under certain conditions on the subset $\S$, for any nodal curve $(C,P_\bullet)$ lying on the boundary of the moduli space $\overline{\mathrm{M}}_{0,n}$, the space of coinvariants $[M^\bullet]_{\mathcal{L}_{C\bs P_\bullet}(V)}$ admits a surjection from $M^\bullet(0)$. For simplicity of notation, write
\begin{equation}\label{eq:Phiast}
	\Phi^\ast(\A_\S) =  \Phi^\L(\A_\S)\ \  \mathrm{or}\ \   \Phi^\sR(\A_\S)^\theta.
\end{equation}
Note that $\Phi^\ast(\A_\S)$ is an admissible $V$-module.

\begin{definition}\label{df:stronglyreducible}
	Assume that $V$ is a strongly rational VOA and $V=U+C_1(V)$.
	Let $\S\ssq \mathscr{W}$ be a subset satisfying the following conditions:
	\begin{enumerate}[($\S$1)]
		\item $\S$ is closed under the fusion tensor product;
		\item Each $M\in \S$ is strongly generated by $U_+$ from its bottom degree $M(0)$ \eqref{eq:Wspanning};
		\item Each $M\in \S$ satisfies \eqref{eq:assumptioncriterion1}: $b_{(n)} M(0)\ssq \mathscr{U}(\widehat{V_1}).M(0)$, for all $b\in U_+\cap V_k$, $k\geq 1$, and $n\geq 0$.
	\end{enumerate}
	We call $\S$ a {\bf (genus-zero) strongly geometrically reducible} set of $V$-modules.
\end{definition}

By Theorem~\ref{thm:criterion1}, if \(\S\) satisfies conditions \((\S2)\) and \((\S3)\), then \(\S\) is geometrically reducible in the sense of Definition~\ref{df:geored}. Examples of subsets \(\S\subset\mathscr{W}\) satisfying these conditions over various VOAs are given in Sections~\ref{sec:3.3}, \ref{sec:3.4}, and \ref{sec:3.5}. Some of the examples in these sections also satisfy condition \((\S1)\). Note that \((\S1)\) forces \(V\in\S\), and that \(V\) always satisfies \((\S2)\) and \((\S3)\).
Recall that if \(V\) is generated by \(V_1\), then the full set \(\mathscr{W}\) is strongly geometrically reducible; see Corollary~\ref{cor:degreeone}.

\subsubsection{One-node case} From now on, unless stated otherwise, $\S\ssq \mathscr{W}$ is a strongly geometrically reducible subset.

\begin{lemma}\label{lm:redAS}
	Let $(C,P_\bullet\sqcup Q)$ be a smooth genus-zero curve, with $P_\bullet=(P_1,\ds, P_n)$ such that $n\geq 2$. Let $M^1,\ds, M^n\in \S$ (possibly with repetitions). Then
	\begin{equation}\label{eq:1Areduction}
		[M^\bullet\o \Phi^\ast(\A_\S)]_{\mathcal{L}_{C\bs P_\bullet\sqcup Q}(V)}=[M^\bullet(0)\o 1_{\A_\S}]=\spn\{[v^\bullet\o 1_{\A_\S}]: v^\bullet\in M^\bullet(0)\}.
	\end{equation}
	In particular, the canonical quotient map is surjective:
	\[
	\psi: M^\bullet(0)\twoheadrightarrow [M^\bullet\o \Phi^\ast(\A_\S)],\quad v^\bullet\mapsto [v^\bullet\o 1_{\A_\S}].
	\]
\end{lemma}
\begin{proof}
	Since $\S$ is closed under taking contragredient modules by Lemma~\ref{lm:rigidity}, by \eqref{eq:MTAforAS} and \eqref{eq:Phiast}, we may write
	$
	\Phi^\ast(\A_\S)\cong \bigoplus_{W\in \S} m_W W,
	$
	where $m_W=\dim W(0)$. By Definition~\ref{df:stronglyreducible} and Theorem~\ref{thm:criterion1}, together with the fact that $\A_\S\cong \bigoplus_{W\in \S} m_W W(0)$ as a left $\A$-module, we have
	\begin{align*}
		[M^\bullet\o \Phi^\ast(\A_\S)]_{\mathcal{L}_{C\bs P_\bullet\sqcup Q}(V)}&\cong \bigoplus_{W\in \S} m_W [M^\bullet\o W]_{\mathcal{L}_{C\bs P_\bullet\sqcup Q}(V)}=\bigoplus_{W\in \S} m_W [M^\bullet(0)\o W(0)]\\
		&\cong  [M^\bullet(0)\o \A_\S].
	\end{align*}
	In other words, $[M^\bullet\o \Phi^\ast(\A_\S)]_{\mathcal{L}_{C\bs P_\bullet\sqcup Q}(V)}=\spn\{[v^\bullet\o \al]: v^\bullet\in M^\bullet(0),\al\in \A_{\S}\}.$ To further reduce $\al\in \A_\S$ to $1_{\A_\S}$, we note that $\al=[a]\cdot 1_{\A_\S}$ for some $[a]\in \A$ (for $\Phi^\ast(\A_\S)=\Phi^\sR(\A_\S)^\theta$ we use instead $\al=1_{\A_\S}\cdot\theta([a])$; since $\theta$ is an anti-automorphism of $\A$, the argument below is the same). Since the left regular $\A$-module action on $\A$ is defined by $o(a).[b]=[a].[b]=[a]\cdot [b]$, and $o(a)=a_{(k-1)}$ for $a\in V_k$, we have
	$\al=a_{(k-1)}. 1_{\A_\S}$. Moreover, by \eqref{eq:Aspanning}, we may assume that $[a]=[a^1]\ds [a^r]$, where $a^i\in U_+\cap V_{k_i}$ for all $i$. Then $\al=a^1_{(k_1-1)}\ds a^r_{(k_r-1)}.1_{\A_\S}$. To show \eqref{eq:1Areduction}, we need to show
	\begin{equation}\label{eq:ntsonenode}
		[v^\bullet\o \al]=[v^\bullet \o a^1_{(k_1-1)}\ds a^r_{(k_r-1)}.1_{\A_\S}]\in [M^\bullet(0)\o 1_{\A_\S}],\quad \forall v^\bullet\in M^\bullet(0).
	\end{equation}
	We use induction on $r\geq 1$. First, consider $\al=a_{(k-1)}.1_{\A_\S}$, where $a\in U_+\cap V_k$. Since $k\geq 1$, there exists a nonzero section $\mu\in H^0(C, \om_C^{1-k}(-(k-1)Q))\bs H^0(C, \om_C^{1-k}(-kQ))$ that is regular at all $P_\bullet$. Let $\si=a\o \mu+\mathrm{Im}\nabla\in \mathcal{L}_{C\bs P_\bullet\sqcup Q}(V)$. Then we may write $\si_Q=a_{(k-1)}+\sum_{i\geq 1} \la_i \cdot a_{(k-1+i)}$, with $\si_Q.1_{\A_\S}=a_{(k-1)}.1_{\A_\S}$. Moreover, since $\si_{P_j}=\sum_{i\geq 0} \la'_i \cdot a_{(i)}$, we have $\si_{P_j}.v^j\in \scrU(\widehat{V_1}).M^j(0)$ by condition $(\S3)$ for all $1\leq j\leq n$. Therefore,
	\begin{align*}
		[v^\bullet\o a_{(k-1)}.1_{\A_\S}]&=-\sum_{j=1}^n[v^1\o\ds \o \si_{P_j}.v^j\o \ds \o v^n\o 1_{\A_\S}]\\
		&\equiv 0\pmod{ [\scrU(\widehat{V_1}).M^\bullet(0)\o 1_{\A_\S}]}
	\end{align*}
	Observe that $1_{\A_\S}\in \Om (\Phi^\ast (\A_\S))$. Then, by Proposition~\ref{prop:reductiononeptleaveout}, we have
	$[\scrU(\widehat{V_1}).M^\bullet(0)\o 1_{\A_\S}]=[M^\bullet(0)\o 1_{\A_\S}]$, and so $[v^\bullet\o a_{(k-1)}.1_{\A_\S}]\in [M^\bullet(0)\o 1_{\A_\S}]$.

	Now let $r\geq 2$. We may write $\al=a^1_{(k_1-1)}\b$, where $\b=a^2_{(k_2-1)}\ds a^r_{(k_r-1)}.1_{\A_\S}\in \Om (\Phi^\ast (\A_\S))$. Similarly, we may choose a nonzero section $\mu'\in H^0(C, \om_C^{1-k_1}(-(k_1-1)Q))\bs H^0(C, \om_C^{1-k_1}(-k_1Q))$, which gives rise to an element $\si'=a^1\o \mu'+\mathrm{Im}\nabla\in  \mathcal{L}_{C\bs P_\bullet\sqcup Q}(V)$. Note that $\si'_Q. \b=a^1_{(k_1-1)}\b=\al$. Hence,
	\[
	[v^\bullet\o \al]=-\sum_{j=1}^n [v^1\o\ds \o \si'_{P_j}.v^j\o \ds \o \b],
	\]
	where $\si'_{P_j}.v^j\in \scrU(\widehat{V_1}).M^j(0)$ for all $1\leq j\leq n$ by $(\S3).$ By Lemma~\ref{lm:UV1spanning}, we may write $\si'_{P_j}.v^j=\sum b^1_{(-n_1)}\ds b^s_{(-n_s)}v'$, where $n_1,\ds, n_s\geq 1$, $b^1,\ds, b^s\in V_1$, and $v'\in M^j(0)$. As in the proof of Proposition~\ref{prop:reductiononeptleaveout}, we may choose a nonzero section $\eta\in H^0(C,\om_C^{1-1}(n_1P_j-Q))\bs H^0(C,\om_C^{1-1}((n_1-1)P_j-Q))$, which gives rise to an element $\tau=b^1\o \eta\in \mathcal{L}_{C\bs P_\bullet\sqcup Q}(V)$. Note that $\tau_Q.\beta=0$, $\tau_{P_j}=b^1_{(-n_1)}+\sum_{i\geq 1}\la_i''\cdot b^1_{(-n_1+i)}$, and $\tau_{P_\ell}.v^\ell=b^1_{(0)}v^\ell\in M^\ell(0)$ for all $\ell\neq j$. We have
	\begin{equation}\label{eq:intermediate}
		\begin{aligned}
			&[v^1\o\ds \o b^1_{(-n_1)}\ds b^s_{(-n_s)}v'\o \ds  \o v^n\o \b]\\
			&=-\sum_{i\geq 1} \la_i''\cdot [v^1\o\ds \o b^1_{(-n_1+i)}b^2_{(-n_2)}\ds b^s_{(-n_s)}v'\o\ds \o v^n\o \b]\\
			&\ \ -\sum _{\ell \neq j} [v^1\o\ds \o b^2_{(-n_2)}\ds b^s_{(-n_s)}v'\o \ds b^1_{(0)}v^\ell  \ds  \o v^n\o \b].
		\end{aligned}
	\end{equation}
	Using induction on the degree of $b^1_{(-n_1)}\ds b^s_{(-n_s)}v'$, we can show that the right-hand side of \eqref{eq:intermediate} is in $[M^\bullet(0)\o \b]$. Now, applying the induction hypothesis of \eqref{eq:ntsonenode} to the shorter-length element $\b$, we see that the right-hand side of \eqref{eq:intermediate} is in $[M^\bullet(0)\o 1_{\A_\S}]$. This completes the induction step for \eqref{eq:ntsonenode}.
\end{proof}

\begin{proposition}\label{prop:onenode}
	Let $(X,P_\bullet\sqcup Q_\bullet)$ be a genus-zero stable curve with one node $\qq$, and let $M^\bullet$ and $W^\bullet$ be irreducible $V$-modules in $\S$ attached to the marked points $P_\bullet$ and $Q_\bullet$, respectively. Then the quotient map
	\begin{equation}\label{coinvbottomnode}
		\pi:	M^\bullet(0)\o W^\bullet(0)\to [M^\bullet\o W^\bullet]_{\mathcal{L}_{X\bs P_\bullet\sqcup Q_\bullet}(V)}\quad u^\bullet\o v^\bullet\mapsto [u^\bullet\o v^\bullet]
	\end{equation}
	is surjective.
\end{proposition}
\begin{proof}
	Since $\S\ssq \mathscr{W}$ is strongly geometrically reducible and the connected components $X_\pm$ of the normalization $\widetilde{X}$ are smooth of genus zero, the conditions in Lemma~\ref{lm:PhiAandAS} and Lemma~\ref{lm:redAS} are satisfied. Therefore, we have the following commutative diagram:
	\[
	\begin{tikzcd}
		&  	{[}M^\bullet\o W^\bullet{]}_{\mathcal{L}_{X\bs P_\bullet\sqcup Q_\bullet}(V)} \arrow[d,"\cong\ h\ \eqref{eq:onenodecoinviso1}"]\\
		M^\bullet(0)\o W^\bullet(0)\arrow[r,dashed,"\pi'"] \arrow[ur,"\pi"]\arrow[dr,"\psi\otimes \psi\ \mathrm{Lemma}\ \ref{lm:redAS}"',two heads]& 	{[}M^\bullet\o \Phi^\L(\A){]}_{\mathcal{L}_+}\o_\A [\Phi^\sR(\A)^\theta\o W^\bullet]_{\mathcal{L}_-} \arrow[d,"\cong\ \mathrm{Lemma}\ \ref{lm:PhiAandAS}"]\\
		&		{[}M^\bullet\o \Phi^\L(\A_\S){]}_{\mathcal{L}_+}\o_\A [\Phi^\sR(\A_\S)^\theta\o W^\bullet]_{\mathcal{L}_-},
	\end{tikzcd}
	\]
	with all maps canonically defined.
	By Lemmas~\ref{lm:PhiAandAS} and \ref{lm:redAS}, the horizontal map $\pi'$ is surjective. Since $h$ \eqref{eq:onenodecoinviso1} is an isomorphism \cite{DGT24}, the map $\pi$ is also surjective.
\end{proof}

\subsubsection{One-node case with a special marked point} Next, we prove a stronger version of Proposition~\ref{prop:onenode}, which serves as a base case for the induction argument later.

\begin{lemma}\label{lm:onenodespecialred}
	Let $(C,\qq_1\sqcup Q_\diamond\sqcup \qq_2)$ be a smooth genus-zero curve, with $Q_\diamond=(Q_1,\ds, Q_m)$ such that $m\geq 1$. Let $W^1,\ds, W^m\subset \S$ (possibly with repetitions). Then
	\[
	\left[\Phi^\ast(\A_\S)\o W^\diamond \o \Phi^\ast(\A_\S)\right]_{\mathcal{L}_{C\bs \qq_1\sqcup Q_\diamond \sqcup \qq_2}(V)}=[\A_\S\o W^\diamond(0)\o 1_{\A_\S}].
	\]
\end{lemma}
\begin{proof}
	Decompose $\Phi^\ast(\A_\S)$ attached at $\qq_1$ as $\bigoplus_{M\in \S} m_M M$. Then the bottom degree is $\Phi^\ast(\A_\S)(0)=\A_\S=\bigoplus_{M\in \S} m_M M(0)$.
	By Lemma~\ref{lm:redAS}, we have
	\begin{align*}
		&\left[\Phi^\ast(\A_\S)\o W^\diamond \o \Phi^\ast(\A_\S)\right]_{\mathcal{L}_{C\bs \qq_1\sqcup Q_\diamond \sqcup \qq_2}(V)}\cong \bigoplus_{M\in \S} m_M \left[M\o W^\diamond \o \Phi^\ast(\A_\S)\right]_{\mathcal{L}_{C\bs \qq_1\sqcup Q_\diamond \sqcup \qq_2}(V)}\\
		&=\bigoplus_{M\in \S} m_M \left[M(0)\o W^\diamond(0) \o 1_{\A_\S}\right]=[\A_\S\o W^\diamond(0)\o 1_{\A_\S}],
	\end{align*}
	where $[M\o W^\diamond \o \Phi^\ast(\A_\S)]_{\mathcal{L}_{C\bs \qq_1\sqcup Q_\diamond \sqcup \qq_2}(V)}= [M(0)\o W^\diamond(0) \o 1_{\A_\S}]$ follows from \eqref{eq:1Areduction}.
\end{proof}

\begin{proposition}\label{prop:onenodeAS}
	Let $(X,P_\bullet\sqcup Q_\diamond\sqcup \qq_2)$ be a genus-zero stable curve with one node $\qq$, such that $P_\bullet$ is contained in one irreducible component of $X$, and $(Q_\diamond,\qq_2)$ are contained in the other component. Let $M^\bullet, W^\diamond \subset \S$, and let $\Phi^\ast(\A_\S)$ be the admissible $V$-module attached at $\qq_2$. Then
	\begin{equation}\label{eq:nodalonepointred}
		[M^\bullet\o W^\diamond\o \Phi^\ast(\A_\S)]_{\mathcal{L}_{X\bs P_\bullet\sqcup Q_\diamond\sqcup \qq_2}(V)}=[M^\bullet(0)\o W^\diamond(0)\o 1_{\A_\S}].
	\end{equation}
\end{proposition}
\begin{proof}
	See the following picture for an illustration:
	\begin{figure}[htbp]
		\centering
		\begin{tikzpicture}[thick, line cap=round, line join=round, scale=0.8]

			\begin{scope}[xshift=0cm]
				\coordinate (A) at (-1,  3);
				\coordinate (B) at ( 2, -3);
				\coordinate (C) at ( 1,  3);
				\coordinate (D) at (-2, -3);

				\draw (A) -- (B);
				\draw (C) -- (D);

				\fill (0,1) circle (2pt);
				\node[left] at (0.8,1) {$\qq$};

				\coordinate (Pw1)   at ($(C)!0.45!(D)$);
				\coordinate (Pw2)   at ($(C)!0.55!(D)$);
				\coordinate (Pw3)   at ($(C)!0.78!(D)$);
				\coordinate (PinfL) at ($(C)!0.9!(D)$);
				\fill (Pw1)   circle (2pt); \node[left]  at (Pw1)   {$P_m$};
				\fill (Pw3)   circle (2pt); \node[left]  at (Pw3)   {$P_{2}$};
				\fill (PinfL) circle (2pt); \node[left]  at (PinfL) {$P_1$};

				\coordinate (Pw1) at ($(A)!0.45!(B)$);
				\coordinate (Pw2) at ($(A)!0.55!(B)$);
				\coordinate (Pw3) at ($(A)!0.78!(B)$);
				\coordinate (P0)  at ($(A)!0.9!(B)$);
				\fill (Pw1) circle (2pt); \node[right] at (Pw1) {$Q_n$};
				\fill (Pw3)   circle (2pt); \node[right]  at (Pw3)   {$Q_{1}$};
				\fill (P0)  circle (2pt); \node[right] at (P0)  {$\qq_2$};

				\node at (0,-3.6) { $(X,P_\bullet\sqcup Q_\diamond)$};

			\end{scope}

			\draw[->, thick]
			(4.2,0) -- (2.5,0)
			node[midway, above] {$\eta$};

			\begin{scope}[xshift=7cm]

				\coordinate (Ltop) at (-0.5,  3);
				\coordinate (Lbot) at (-1.5, -3);

				\coordinate (Rtop) at ( 1,  3);
				\coordinate (Rbot) at ( 2, -3);

				\draw (Lbot) -- (Ltop);
				\draw (Rtop) -- (Rbot);

				\coordinate (L0)    at ($(Ltop)!0.32!(Lbot)$);
				\coordinate (LQ1)   at ($(Ltop)!0.45!(Lbot)$);
				\coordinate (LQ2)   at ($(Ltop)!0.55!(Lbot)$);
				\coordinate (LQ3)   at ($(Ltop)!0.8!(Lbot)$);
				\coordinate (Linf)  at ($(Ltop)!0.9!(Lbot)$);
				\fill (L0)   circle (2pt); \node[left] at (L0)   {$\qq_+$};
				\fill (LQ1)  circle (2pt); \node[left] at (LQ1)  {$P_m$};
				\fill (LQ3)  circle (2pt); \node[left] at (LQ3)  {$P_{2}$};
				\fill (Linf) circle (2pt); \node[left] at (Linf) {$P_1$};

				\coordinate (Rinf) at ($(Rtop)!0.32!(Rbot)$);
				\coordinate (RR1)  at ($(Rtop)!0.45!(Rbot)$);
				\coordinate (RR2)  at ($(Rtop)!0.55!(Rbot)$);
				\coordinate (RR3)  at ($(Rtop)!0.8!(Rbot)$);
				\coordinate (R0)   at ($(Rtop)!0.9!(Rbot)$);
				\fill (Rinf) circle (2pt); \node[right] at (Rinf) {$\qq_-$};
				\fill (RR1)  circle (2pt); \node[right] at (RR1)  {$Q_n$};
				\fill (RR3)  circle (2pt); \node[right] at (RR3)  {$Q_{1}$};
				\fill (R0)   circle (2pt); \node[right] at (R0)   {$\qq_2$};

				\node at (0,-3.6) {$(X_+,P_\bullet\sqcup \qq_+)\quad (X_-,\qq_-\sqcup Q_\diamond\sqcup \qq_2)$};

			\end{scope}
		\end{tikzpicture}
		\caption{One-node curve with a special marked point $\qq_2$}
		\label{fig:specialmarkedpoint}
	\end{figure}
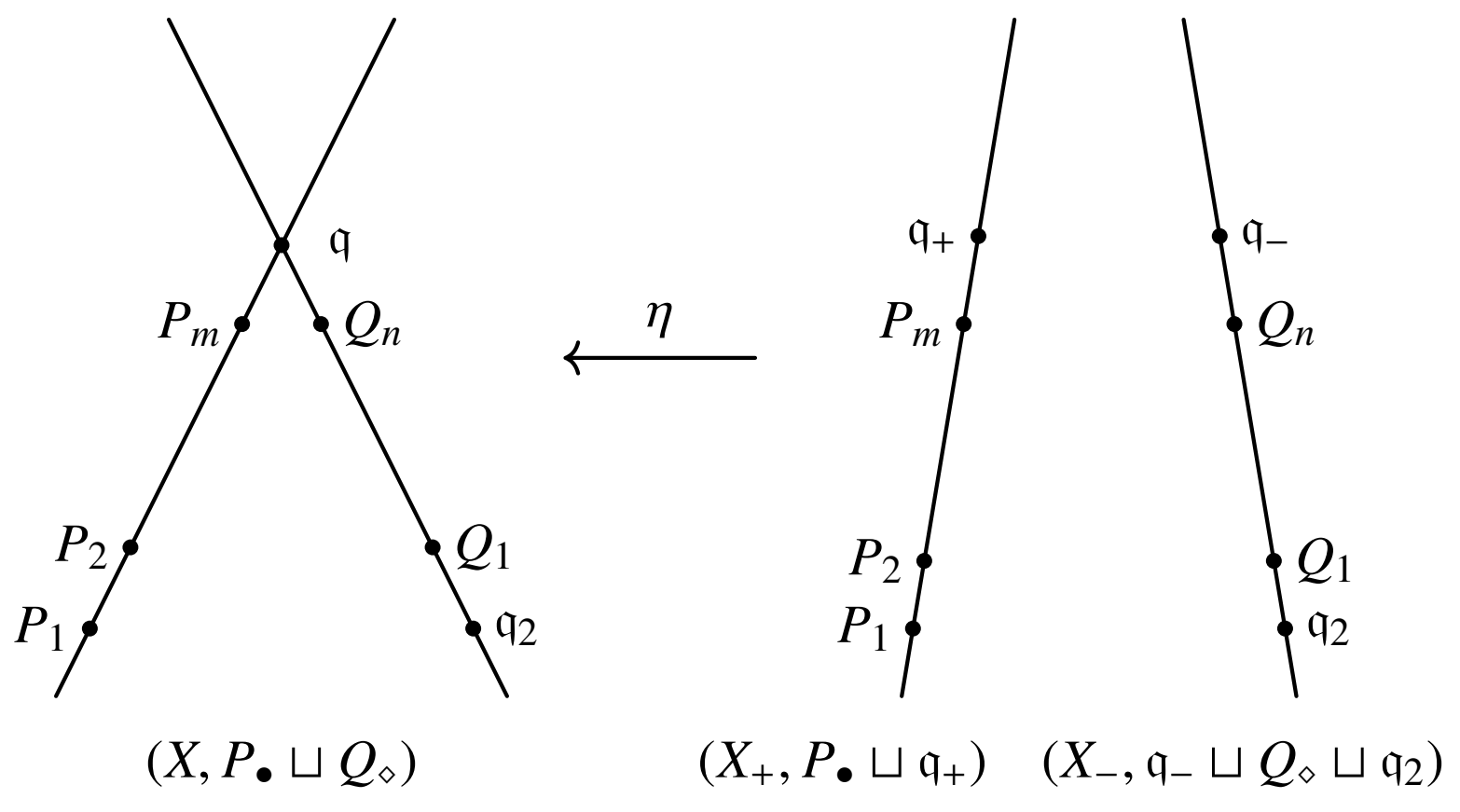

	Let $\mathcal{L}_+=\mathcal{L}_{X_+\bs P_\bullet\sqcup \qq_+}(V)$ and $\mathcal{L}_-=\mathcal{L}_{X_-\bs \qq_-\sqcup Q_\diamond\sqcup \qq_2}(V)$.
	By \eqref{eq:4.10}, we have an isomorphism:
	\begin{align*}
		h:[M^\bullet\o W^\diamond\o \Phi^\ast(\A_\S)]_{\mathcal{L}_{X\bs P_\bullet\sqcup Q_\diamond\sqcup \qq_2}(V)}&\cong [M^\bullet\o \Phi^\L(\A_\S)]_{\mathcal{L}_+}\o_\A [\Phi^\sR(\A_\S)^\theta\o W^\diamond\o \Phi^\ast(\A_\S)]_{\mathcal{L}_-}\\
		[x^\bullet\o w^\diamond\o \xi]&\mapsto [x^\bullet\o 1_{\A_\S}]\o_\A [1_{\A_\S}\o w^\diamond\o \xi].
	\end{align*}
	Therefore, to show \eqref{eq:nodalonepointred}, it suffices to show
	\begin{equation}\label{eq:4.18}
		\begin{aligned}
			&[M^\bullet\o \Phi^\L(\A_\S)]_{\mathcal{L}_+}\o_\A [\Phi^\sR(\A_\S)^\theta\o W^\diamond\o \Phi^\ast(\A_\S)]_{\mathcal{L}_-}\\
			&=\spn\{[u^\bullet\o 1_{\A_\S}]\o_\A [1_{\A_\S}\o v^\diamond\o 1_{\A_\S}]: u^\bullet\in M^\bullet(0), v^\diamond\in W^\diamond(0) \},
		\end{aligned}
	\end{equation}
	since the preimage of $[u^\bullet\o 1_{\A_\S}]\o_\A [1_{\A_\S}\o v^\diamond\o 1_{\A_\S}]$ under $h$ is $[u^\bullet\o v^\diamond\o 1_{\A_\S}]\in [M^\bullet(0)\o W^\diamond(0)\o 1_{\A_\S}]$.
	Note that $\A_\S=\spn\{\al\ast 1_{\A_\S}:\al=a^1_{(k_1-1)}\ds a^r_{(k_r-1)}\in \A \ \}$, where $\al\ast 1_{\A_\S}$ is the auxiliary left multiplication action of $\A$ on $\A_\S$, as explained in Section~\ref{sec:4.1}.
	By Lemma~\ref{lm:onenodespecialred}, we have $$[\Phi^\sR(\A_\S)^\theta\o W^\diamond\o \Phi^\ast(\A_\S)]_{\mathcal{L}_-}=\spn\{\al\ast 1_{\A_\S}\o v^\diamond\o 1_{\A_\S}: \al\in \A,v^\diamond\in W^\diamond(0) \}.$$
	On the other hand, by \eqref{eq:1Areduction}, we have $[M^\bullet\o \Phi^\L(\A_\S)]_{\mathcal{L}_+}=[M^\bullet(0)\o 1_{\A_\S}]$. Therefore, a general spanning element in $[M^\bullet\o \Phi^\L(\A_\S)]_{\mathcal{L}_+}\o_\A [\Phi^\sR(\A_\S)^\theta\o W^\diamond\o \Phi^\ast(\A_\S)]_{\mathcal{L}_-}$ has the form
	\[
	[u^\bullet \o 1_{\A_\S}]\o_\A [\al\ast 1_{\A_\S}\o v^\diamond\o 1_{\A_\S}]=[u^\bullet \o \al.1_{\A_\S}]\o_\A [1_{\A_\S}\o v^\diamond\o 1_{\A_\S}],
	\]
	where we used the fact that $\al\in \A$ can be moved from right to left within the tensor product over $\A$, and $1_{\A_\S}\ast \al=\al.1_{\A_\S}$ in $\Phi^\L(\A_\S)(0)$.
	Applying \eqref{eq:1Areduction} again, $[u^\bullet \o \al.1_{\A_\S}]$ is a sum of elements of the form $[\tilde{u}^\bullet\o 1_{\A_\S}]$, where $\tilde{u}^\bullet\in M^\bullet(0)$. This proves \eqref{eq:4.18}.
\end{proof}



\subsubsection{The multiple-node case and global generation}
Note that the boundary $\partial \overline{\mathrm{M}}_{0,n}$ of the moduli space $\overline{\mathrm{M}}_{0,n}$ consists of stable rational curves with $n$ marked points and at least one node. Let $[n]=\{1,2,\ds, n\}.$ An element of $\partial \overline{\mathrm{M}}_{0,n}$ can be represented by
\[
C=C_1\cup \ds \cup C_t,\quad [n]=\bigsqcup_{j=1}^t I_j,
\]
where $C_1,\ds, C_t$ are the irreducible components of $C$, and $I_j$ is the set of labels of the marked points $P^j_\bullet$ lying on $C_j$. Each component satisfies
\[
C_j\simeq \P^1,\quad  |I_j|+|C_j^{\mathrm{node}}|\geq 3, \quad \forall 1\leq j\leq t,
\]
where $C_j^{\mathrm{node}}$ denotes the set of nodes lying on $C_j$. Moreover, the dual graph of $C$ is a tree.
Let $\S\ssq \mathscr{W}$ be a strongly geometrically reducible subset. We will use induction on the number of nodes of $C$, with the results in the previous subsection as the base case, to show that
\[
[M^\bullet]_{\mathcal{L}_{C\bs P_\bullet}(V)}=[M^\bullet(0)],
\]
where $M^\bullet \subset \S$, and $P_\bullet$ are smooth marked points on $C$ \eqref{eq:classincoinv}. Therefore, the quotient map is surjective: $M^\bullet(0)\twoheadrightarrow [M^\bullet]_{\mathcal{L}_{C\bs P_\bullet}(V)}$.

\subsubsection{Closure under fusion and mode transition algebra}
Let $\qq$ be a node on some connected component $C_j$ of $C$. Let $\widetilde{C}=C_+\sqcup C_-$ be the normalization at $\qq$. Let $P_\bullet$ and $Q_\diamond$ be the sets of smooth marked points on $C_+$ and $C_-$, respectively.

We first generalize Lemma~\ref{lm:PhiAandAS} from the one-node case to the general case. In other words, given $M^\bullet,W^\diamond\subset \S$, we want to show
\begin{equation}\label{eq:multinodeAS}
	\begin{aligned}
		[M^\bullet\o W^\diamond]_{\mathcal{L}_{C\bs P_\bullet\sqcup Q_\diamond}(V)}
		&\cong [M^\bullet \o \Phi^\L(\A_\S)]_{\mathcal{L}_{C_+\bs P_\bullet\sqcup \qq_+}(V)}\o_\A [\Phi^\sR(\A_\S)^\theta\o W^\diamond]_{\mathcal{L}_{C_-\bs Q_\diamond\sqcup \qq_-}(V)}\\
		[u^\bullet\o w^\diamond]&\mapsto [u^\bullet\o 1_{\A_\S}]\o_\A [1_{\A_\S}\o w^\diamond].
	\end{aligned}
\end{equation}
Note that \eqref{eq:multinodeAS} is not immediate, since, a priori, the admissible modules attached at $\qq_+$ and $\qq_-$ are $\Phi^\L(\A)$ and $\Phi^\sR(\A)^\theta$, respectively, by \eqref{eq:onenodecoinviso1}.

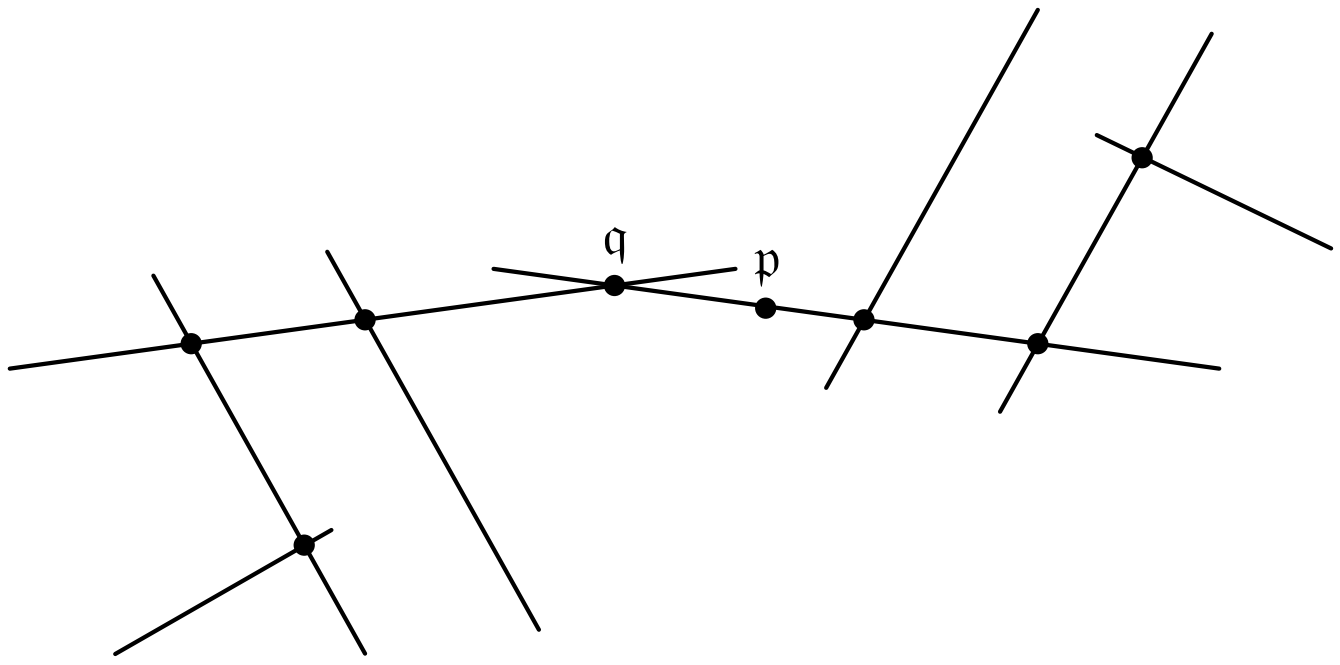
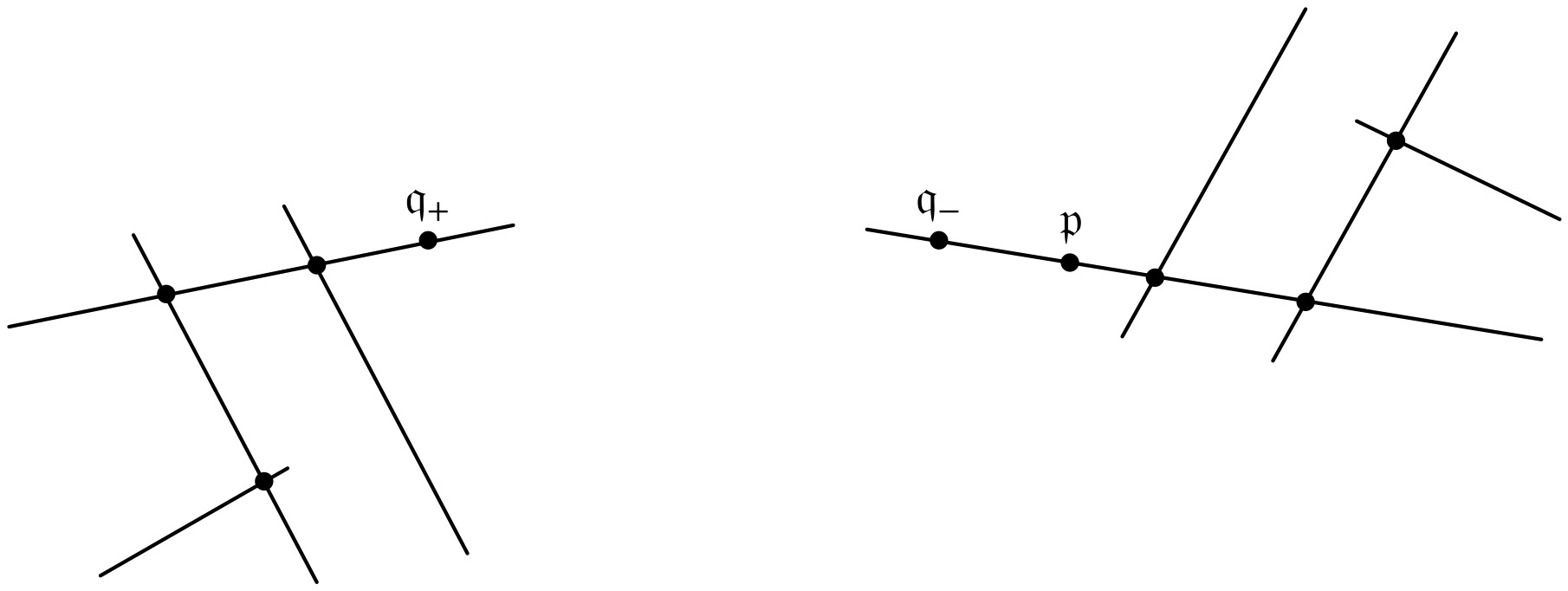
\begin{figure}[ht]
	\centering
	\begin{tikzpicture}[
		thick,
		line cap=round,
		line join=round,
		scale=1
		]

		\coordinate (q) at (0,0);
		\coordinate (p) at (1,-0.15);

		\draw (-4,-0.55) -- (0.8,0.11);

		\draw (-0.8,0.11) -- (4,-0.55);

		\fill (q) circle (2pt);
		\node[above] at (q) {$\qq$};

		\fill (p) circle (2pt);
		\node[above] at (p) {$\pp$};

		\coordinate (q1) at (-2.8,-0.385);

		\coordinate (q2) at (-1.65,-0.227); 

		\coordinate (q3) at ( 1.65,-0.227); 

		\coordinate (q4) at ( 2.8,-0.385);

		\fill (q1) circle (2pt);
		\fill (q2) circle (2pt);
		\fill (q3) circle (2pt);
		\fill (q4) circle (2pt);


		\coordinate (I1end) at ($(q1)+(1.15,-2.05)$);
		\coordinate (I2end) at ($(q2)+(1.15,-2.05)$);

		\draw ($(q1)+(-0.25,0.45)$) -- (I1end);
		\draw ($(q2)+(-0.25,0.45)$) -- (I2end);


		\coordinate (r1) at ($(q1)!0.65!(I1end)$);
		\fill (r1) circle (2pt);

		\draw
		($(r1)+(0.18,0.10)$)
		--
		($(r1)+(-1.25,-0.72)$);


		\coordinate (I3end) at ($(q3)+(1.15,2.05)$);
		\coordinate (I4end) at ($(q4)+(1.15,2.05)$);

		\draw ($(q3)+(-0.25,-0.45)$) -- (I3end);
		\draw ($(q4)+(-0.25,-0.45)$) -- (I4end);


		\coordinate (r4) at ($(q4)!0.6!(I4end)$);
		\fill (r4) circle (2pt);

		\draw
		($(r4)+(-0.3,0.15)$)
		--
		($(r4)+(1.25,-0.6)$);

		\node at (0,-3) { $(C,P_\bullet\sqcup \pp)$};
	\end{tikzpicture}
	\caption{A general nodal curve $(C,P_\bullet)\in \partial \overline{\mathrm{M}}_{0,n}$}
	\label{fig:Fcurve-general}
\end{figure}

\begin{lemma}\label{lm:notinSvanishing}
	Given a boundary point $(C,P_\bullet\sqcup \pp)\in  \partial \overline{\mathrm{M}}_{0,n}$, with a particular smooth marked point $\pp$ lying on an irreducible component, see Figure~\ref{fig:Fcurve-general}. Let $M^\bullet\subset \S$, but let $\tilde{M}\in \mathscr{W}\bs \S$. Then
	\[
	[M^\bullet\o \tilde{M}]_{\mathcal{L}_{C\bs P_\bullet\sqcup \pp}(V)}=0.
	\]
\end{lemma}
\begin{proof}
	We use induction on the total number of nodes $\mathcal{N}=|C^{\mathrm{node}} |$. If $\mathcal{N}=1$, let $\qq\in C$ be the node, and let $\tilde{C}=C_+\sqcup C_-$ be the normalization at $\qq$. Assume that $\pp\in C_-$. Then, by \eqref{eq:onenodecoinviso1}, we have
	\[
	[M^\bullet\o \tilde{M}]_{\mathcal{L}_{C\bs P_\bullet\sqcup \pp}(V)}\cong [M^\star\o \Phi^\L(\A)]_{\mathcal{L}_+}\o_\A [\Phi^\sR(\A)^\theta\o M^\circ \o \tilde{M}]_{\mathcal{L}_-},
	\]
	where $\bullet=\star\sqcup \circ$ are the labels of the marked points.
	Moreover, since $M^\star\subset \S$, it follows from the proof of Lemma~\ref{lm:PhiAandAS} that
	$
	[M^\star\o \Phi^\L(\A)]_{\mathcal{L}_+}=[M^\star\o \Phi^\L(\A_\S)]_{\mathcal{L}_+}
	$, and, by \eqref{eq:AStensorA},
	\[
	[M^\bullet\o \tilde{M}]_{\mathcal{L}_{C\bs P_\bullet\sqcup \pp}(V)}\cong [M^\star\o \Phi^\L(\A_\S)]_{\mathcal{L}_+}\o_\A [\Phi^\sR(\A_\S)^\theta\o M^\circ \o \tilde{M}]_{\mathcal{L}_-}.
	\]
	By Lemma~\ref{lm:rigidity}, $\S$ is closed under taking contragredient modules, and $\dim W(0)=\dim W(0)^\ast=m_W$. Hence, we can write $\Phi^\sR(\A_\S)^\theta=\bigoplus_{W\in\S}m_W W$, in view of \eqref{eq:MTAforAS}.
	Since $\mathcal{L}_-$ is the chiral Lie algebra associated with a smooth genus-zero curve $(C_-,\qq_-\sqcup P_\circ\sqcup \pp)$, by Lemma~\ref{lm:subring}, we have $[\Phi^\sR(\A_\S)^\theta\o M^\circ \o \tilde{M}]_{\mathcal{L}_-}=0$. Hence, $[M^\bullet\o \tilde{M}]_{\mathcal{L}_{C\bs P_\bullet\sqcup \pp}(V)}=0$.

	Now assume that $C$ has $\mathcal{N}>1$ nodes. Choose a node $\qq$ of $C$, and let
	$\tilde{C}=C_+\sqcup C_- $ be the normalization at $\qq$; see Figure~\ref{fig:Fcurve-normalization}. Then, by \eqref{eq:onenodecoinviso1}, we have
	\[
	[M^\bullet\o \tilde{M}]_{\mathcal{L}_{C\bs P_\bullet\sqcup \pp}(V)}\cong [M^\star\o \Phi^\L(\A)]_{\mathcal{L}_+}\o_\A [\Phi^\sR(\A)^\theta\o M^\circ \o \tilde{M}]_{\mathcal{L}_-},
	\]
	where $\mathcal{L}_+=\mathcal{L}_{C_+\bs P_\star\sqcup \qq_+}(V)$ and $\mathcal{L}_-=\mathcal{L}_{C_-\bs  \qq_-\sqcup P_\circ\sqcup \pp}(V)$.
	Since $|C_+^{\mathrm{node}}|<\mathcal{N}$, by the induction hypothesis,
	we have $[M^\star\o W]_{\mathcal{L}_+}=0$ for all $W\in \mathscr{W}\bs \S$. Note also that if $C_+^{\mathrm{node}}=\emptyset$, then $[M^\star\o W]_{\mathcal{L}_+}=0$ by Lemma~\ref{lm:subring}.
	Hence,
	\[
	[M^\star\o \Phi^\L(\A)]_{\mathcal{L}_+}\cong \bigoplus_{W\in \mathscr{W}} m_W[M^\star\o W]_{\mathcal{L}_+}\cong \bigoplus_{W\in \S} m_W[M^\star\o W]_{\mathcal{L}_+}\cong [M^\star\o \Phi^\L(\A_\S)]_{\mathcal{L}_+}.
	\]
	Now it follows from \eqref{eq:AStensorA} that
	\[
	[M^\bullet\o \tilde{M}]_{\mathcal{L}_{C\bs P_\bullet\sqcup \pp}(V)}\cong [M^\star\o \Phi^\L(\A_\S)]_{\mathcal{L}_+}\o_\A [\Phi^\sR(\A_\S)^\theta\o M^\circ \o \tilde{M}]_{\mathcal{L}_-}.
	\]
	Since we also have $|C_-^{\mathrm{node}}|<\mathcal{N}$ and $\Phi^\sR(\A_\S)^\theta$ is a direct sum of modules in $\S$, by the induction hypothesis, $[\Phi^\sR(\A_\S)^\theta\o M^\circ \o \tilde{M}]_{\mathcal{L}_-}=0$; if $C_-$ is smooth, this follows from Lemma~\ref{lm:subring}.
\end{proof}

\begin{figure}[ht]
	\centering
	\begin{tikzpicture}[
		thick,
		line cap=round,
		line join=round,
		scale=1
		]


		\coordinate (qL) at (-0.8,0.11);

		\draw (-4,-0.55) -- (-0.15,0.225);

		\fill (qL) circle (2pt);
		\node[above] at (qL) {$\qq_+$};

		\coordinate (q1) at (-2.8,-0.3);
		\coordinate (q2) at (-1.65,-0.08);

		\fill (q1) circle (2pt);
		\fill (q2) circle (2pt);

		\coordinate (I1end) at ($(q1)+(1.15,-2.2)$);
		\coordinate (I2end) at ($(q2)+(1.15,-2.2)$);

		\draw
		($(q1)+(-0.25,0.45)$)
		--
		(I1end);

		\draw
		($(q2)+(-0.25,0.45)$)
		--
		(I2end);

		\coordinate (r1) at ($(q1)!0.65!(I1end)$);
		\fill (r1) circle (2pt);

		\draw
		($(r1)+(0.18,0.10)$)
		--
		($(r1)+(-1.25,-0.72)$);
		\node at (-1.5,-3) { $(C_+,P_\star\sqcup \qq_+)$};

		\begin{scope}[xshift=2.3cm]

			\coordinate (qR) at (0.8,0.11);
			\coordinate (p) at (1.8,-0.06);

			\draw (0.25,0.193) -- (5.4,-0.647);

			\fill (qR) circle (2pt);
			\node[above] at (qR) {$\qq_-$};

			\fill (p) circle (2pt);
			\node[above] at (p) {$\pp$};

			\coordinate (q3) at (2.45,-0.175);
			\coordinate (q4) at (3.60,-0.36);

			\fill (q3) circle (2pt);
			\fill (q4) circle (2pt);

			\coordinate (I3end) at ($(q3)+(1.15,2.05)$);
			\coordinate (I4end) at ($(q4)+(1.15,2.05)$);

			\draw
			($(q3)+(-0.25,-0.45)$)
			--
			(I3end);

			\draw
			($(q4)+(-0.25,-0.45)$)
			--
			(I4end);

			\coordinate (r4) at ($(q4)!0.6!(I4end)$);
			\fill (r4) circle (2pt);

			\draw
			($(r4)+(-0.3,0.15)$)
			--
			($(r4)+(1.25,-0.6)$);
			\node at (2.5,-3) { $(C_-, \qq_-\sqcup \pp\sqcup P_\circ)$};
		\end{scope}

	\end{tikzpicture}
	\caption{Normalization of $(C,P_\bullet\sqcup \pp)$ at the node $\qq$}
	\label{fig:Fcurve-normalization}
\end{figure}

Now \eqref{eq:multinodeAS} is an easy corollary of Lemma~\ref{lm:notinSvanishing}. Indeed, by \eqref{eq:onenodecoinviso1}, we have
\[
[M^\bullet\o W^\diamond]_{\mathcal{L}_{C\bs P_\bullet\sqcup Q_\diamond}(V)}\cong [M^\bullet \o \Phi^\L(\A)]_{\mathcal{L}_{C_+\bs P_\bullet\sqcup \qq_+}(V)}\o_\A [\Phi^\sR(\A)^\theta\o W^\diamond]_{\mathcal{L}_{C_-\bs Q_\diamond\sqcup \qq_-}(V)}.
\]
If $C_+$ has at least one node, then, by Lemma~\ref{lm:notinSvanishing},
$$[M^\bullet \o \Phi^\L(\A)]_{\mathcal{L}_{+}}=\bigoplus_{W\in \mathscr{W}} m_W [M^\bullet \o W]_{\mathcal{L}_{+}}=\bigoplus_{W\in S} m_W [M^\bullet \o W]_{\mathcal{L}_{+}}\cong[M^\bullet \o \Phi^\L(\A_\S)]_{\mathcal{L}_{+}}.$$
The same equality holds if $C_+$ is smooth, in view of Lemma~\ref{lm:subring}.
Similarly, we have $[\Phi^\sR(\A)^\theta\o W^\diamond]_{\mathcal{L}_{C_-\bs Q_\diamond\sqcup \qq_-}(V)}= [\Phi^\sR(\A_\S)^\theta\o W^\diamond]_{\mathcal{L}_{C_-\bs Q_\diamond\sqcup \qq_-}(V)}.$

\subsubsection{Reduction of nodal coinvariants with $\Phi^\ast(\A_\S)$} Now we generalize Proposition~\ref{prop:onenodeAS} to the multiple-node case.
\begin{proposition}\label{prop:arbitrarynodal}
	Given a boundary point $(C,P_\bullet\sqcup \pp)\in  \partial \overline{\mathrm{M}}_{0,n}$, with a particular smooth marked point $\pp$ lying on an irreducible component, see Figure~\ref{fig:Fcurve-general}. Let $M^\bullet\subset \S$. With the notation \eqref{eq:Phiast}, we have
	\begin{equation}\label{eq:PhiASreductionmultiode}
		[M^\bullet\o \Phi^\ast(\A_\S)]_{\mathcal{L}_{C\bs P_\bullet\sqcup \pp}(V)}=[M^\bullet(0)\o 1_{\A_\S}].
	\end{equation}
\end{proposition}
\begin{proof}
	We again use induction on the number of nodes $\mathcal{N}=|C^{\mathrm{node}}|$. If $\mathcal{N}=1$, \eqref{eq:PhiASreductionmultiode} follows from Proposition~\ref{prop:onenodeAS}.

	Assume that $C$ has more than one node. Pick a node $\qq\in C$, and let $\tilde{C}=C_+\sqcup C_-$ be the normalization at $\qq$; see Figure~\ref{fig:Fcurve-general}. Assume that the smooth marked points on $C_\pm$ are given by $(C_+,P_\star\sqcup \qq_+)$ and $(C_-,\qq_-\sqcup P_\circ\sqcup \pp)$, where $\bullet=\star\sqcup \circ$; see Figure~\ref{fig:Fcurve-normalization}.

	Then, by \eqref{eq:multinodeAS}, we have an isomorphism
	\begin{equation}\label{eq:multinodeASAS}
		\begin{aligned}
			[M^\star \o M^\circ \o \Phi^\ast(\A_\S)]_{\mathcal{L}_{C\bs P_\bullet\sqcup \pp} (V)}
			&\cong [M^\star \o \Phi^\L(\A_\S)]_{\mathcal{L}_+}\o_\A [\Phi^\sR(\A_\S)^\theta\o M^\circ \o \Phi^\ast(\A_\S)]_{\mathcal{L}_{-}}\\
			[u^\star\o u^\circ\o \xi]&\mapsto [u^\star\o 1_{\A_\S}]\o_\A [1_{\A_\S}\o u^\circ\o \xi],
		\end{aligned}
	\end{equation}
	where $\mathcal{L}_+=\mathcal{L}_{C_+\bs P_\star\sqcup \qq_+}(V)$ and $\mathcal{L}_-=\mathcal{L}_{C_-\bs  \qq_-\sqcup P_\circ\sqcup \pp}(V)$, respectively.

	First, consider the curve $(C_-,\qq_-\sqcup P_\circ\sqcup \pp)$. Since $\Phi^\sR(\A_\S)^\theta\cong \bigoplus_{W\in \S}m_W W$, we have
	\[
	[\Phi^\sR(\A_\S)^\theta\o M^\circ \o \Phi^\ast(\A_\S)]_{\mathcal{L}_{-}}\cong \bigoplus_{W\in \S} m_W [ W\o M^\circ \o \Phi^\ast(\A_\S)]_{\mathcal{L}_-}.
	\]
	Since $C_-$ has fewer nodes than $C$, by the induction hypothesis for \eqref{eq:PhiASreductionmultiode}, we have
	$ [ W\o M^\circ \o \Phi^\ast(\A_\S)]_{\mathcal{L}_-}=[W(0)\o M^\circ(0)\o 1_{\A_\S}]$ for each $W\in \S$. Therefore,
	$$[\Phi^\sR(\A_\S)^\theta\o M^\circ \o \Phi^\ast(\A_\S)]_{\mathcal{L}_{-}}=[\A_\S\o M^\circ(0)\o 1_{\A_\S}].$$
	Note that if $C_-$ is smooth, then this equality follows from Lemma~\ref{lm:onenodespecialred}.
	In particular, we have
	\begin{align*}
		&[M^\star \o \Phi^\L(\A_\S)]_{\mathcal{L}_+}\o_\A [\Phi^\sR(\A_\S)^\theta\o M^\circ \o \Phi^\ast(\A_\S)]_{\mathcal{L}_{-}}\\
		&=\spn\{[u^\star\o \eta]\o_\A [\al.1_{\A_\S}\o v^\circ\o 1_{\A_\S}]: \eta\in \Phi^\L(\A_\S),\al=a^1_{(k_1-1)}\ds a^r_{(k_r-1)}\in \A,v^\circ\in M^\circ(0) \}.
	\end{align*}
	Now applying the induction hypothesis for \eqref{eq:PhiASreductionmultiode} to the curve $(C_+,P_\star\sqcup \qq_+)$, we have
	\begin{align*}
		&[u^\star\o \eta]\o_\A [\al.1_{\A_\S}\o v^\circ\o 1_{\A_\S}]=[u^\star\o \eta \cdot \al ]\o_\A [1_{\A_\S}\o v^\circ\o 1_{\A_\S}]\\
		&\equiv 0\mod{\spn\{[v^\star\o 1_{\A_\S}]\o_\A [1_{\A_\S}\o v^\circ\o 1_{\A_\S}]: v^\star \in  M^\star(0)\} }.
	\end{align*}
	Note that if $C_+$ is smooth, then the congruence above still holds by Lemma~\ref{lm:redAS}.
	Therefore, elements of $ [M^\bullet \o \Phi^\L(\A_\S)]_{\mathcal{L}_+}\o_\A [\Phi^\sR(\A_\S)^\theta\o M^\circ \o \Phi^\ast(\A_\S)]_{\mathcal{L}_{-}}$ can be written as linear combinations of elements $[v^\star\o 1_{\A_\S}]\o_\A [1_{\A_\S}\o v^\circ\o 1_{\A_\S}]$, where $v^\star\in M^\star(0)$ and $v^\circ \in M^\circ (0)$. Then, by the isomorphism \eqref{eq:multinodeASAS}, we have
	\[
	[M^\star \o M^\circ \o \Phi^\ast(\A_\S)]_{\mathcal{L}_{C\bs P_\bullet\sqcup \pp} (V)}=\spn\{[v^\star\o v^\circ\o 1_{\A_\S}]: v^\star\in M^\star(0),\ v^\circ \in M^\circ (0)\}.
	\]
	This proves \eqref{eq:PhiASreductionmultiode} for a general nodal curve $C$ and completes the induction step.
\end{proof}
\subsubsection{Final reduction of nodal coinvariants}
\begin{theorem}\label{thm:main1}
	Let $\S\subset \mathscr{W}$ be a strongly geometrically reducible subset of irreducible $V$-modules, and let $M^\bullet\subset \S$ (possibly with repetitions). Then the vector bundle of VOA-coinvariants $ \V(V;M^\bullet)$ on the moduli space $\overline{\mathrm{M}}_{0,n}$ is globally generated.
\end{theorem}
\begin{proof}
Recall that $\widetriangle{\mathrm{M}}_{0,n}$ denotes the moduli space of coordinatized pointed stable rational curves $(C,P_\bullet,t_\bullet)$, where $t_\bullet$ is an $n$-tuple of formal coordinates at the marked points, that the forgetful map $\pi: \widetriangle{\mathrm{M}}_{0,n}\ra \overline{\mathrm{M}}_{0,n}$ is a torsor under the group scheme $(\operatorname{Aut}\O)^n$ of changes of formal coordinates, and that $\V(V;M^\bullet)$ is obtained from $\widetriangle{\V}(V;M^\bullet)$ by descent along $\pi$; see \cite[Sections 2.2.2 and 8.7]{DGT24}. By \cite[Lemma 3.2.1]{DG23}, the bottom degrees of $M^\bullet$ give rise to a constant sheaf:
	\[
	\pi_\ast\left(M^\bullet(0)\o \mathcal{O}_{\widetriangle{\mathrm{M}}_{0,n}}\right)^{(\mathrm{Aut}\ \O)^n} =M^\bullet(0)\o \left(\pi_\ast \mathcal{O}_{\widetriangle{\mathrm{M}}_{0,n}}\right)^{(\mathrm{Aut}\ \O)^n}\cong M^\bullet(0)\o \O_{\overline{\mathrm{M}}_{0,n}},
	\]
	where $\pi$ is as above. This gives rise to a natural morphism
	\[
	M^\bullet(0)\o \O_{\overline{\mathrm{M}}_{0,n}}\ra  \V(V;M^\bullet)
	\]
	defined by the quotient map from $M^\bullet(0)$ to the coinvariants.

	Since all modules in $\S$ satisfy conditions \eqref{eq:Wspanning} and \eqref{eq:assumptioncriterion1}, by Theorem~\ref{thm:criterion1} and Definition~\ref{df:geored}, given a smooth curve $(C,P_\bullet)$ corresponding to a point in the interior of $\overline{\mathrm{M}}_{0,n}$, the quotient map is surjective:
	\begin{equation}\label{eq:smoothsurj}
		M^\bullet(0)\ra   \V(V;M^\bullet)_{(C,P_\bullet)}\quad v^\bullet\mapsto [v^\bullet].
	\end{equation}
	Here we also use the fact that the space of coinvariants does not depend on the choice of local coordinates $t_\bullet$ around $P_\bullet$ \cite[Section 8]{DGT24}, i.e., $ \V(V;M^\bullet)_{(C,P_\bullet)}\cong \V(V;M^\bullet)_{(C,P_\bullet,t_\bullet)}$.

	Now let $(C,P_\bullet)\in  \partial \overline{\mathrm{M}}_{0,n}$ correspond to a nodal curve with $n$ smooth marked points $P_\bullet$, and let $M^\bullet\subset \S$. We claim that the canonical quotient map from the bottom degrees to the space of nodal coinvariants is surjective:
	\begin{equation}\label{eq:generalsurj}
		M^\bullet(0)\twoheadrightarrow	[M^\bullet]_{\mathcal{L}_{C\bs P_\bullet}(V)}=\V(V;M^\bullet)_{(C,P_\bullet)},\quad v^\bullet\mapsto [v^\bullet].
	\end{equation}

	Indeed, if $C$ has only one node, then \eqref{eq:generalsurj} follows from Proposition~\ref{prop:onenode}. Otherwise, choose a node $\qq$ of $C$, and let $\tilde{C}=C_+\sqcup C_-$ be the normalization at $\qq$, such that the smooth marked points of $C_+$ and $C_-$ are labeled by $P_\star$ and $P_\circ$, respectively, where $\bullet=\star\sqcup \circ$.
	By \eqref{eq:multinodeAS}, we have an isomorphism:
	\begin{equation}\label{eq:mutinodegeneral}
		\begin{aligned}
			[M^\star\o M^\circ]_{\mathcal{L}_{C\bs P_\bullet}(V)}
			&\cong [M^\star \o \Phi^\L(\A_\S)]_{\mathcal{L}_{C_+\bs P_\star\sqcup \qq_+}(V)}\o_\A [\Phi^\sR(\A_\S)^\theta\o M^\circ]_{\mathcal{L}_{C_-\bs  \qq_-\sqcup P_\circ}(V)}\\
			[u^\star\o u^\circ]&\mapsto [u^\star\o 1_{\A_\S}]\o_\A [1_{\A_\S}\o u^\circ].
		\end{aligned}
	\end{equation}
	Assume that $|C_+^{\mathrm{node}}|\leq |C_-^{\mathrm{node}}|$. If $|C_+^{\mathrm{node}}|\geq 1$, then, by Proposition~\ref{prop:arbitrarynodal}, we have $$[M^\star \o \Phi^\L(\A_\S)]_{\mathcal{L}_{C_+\bs P_\star\sqcup \qq_+}(V)}=[M^\star(0)\o 1_{\A_\S}].$$
	If $|C_+^{\mathrm{node}}|=0$, then, by Lemma~\ref{lm:redAS}, the reduction above still holds. Since $|C_-^{\mathrm{node}}|\geq 1$, we also have $[\Phi^\sR(\A_\S)^\theta\o M^\circ]_{\mathcal{L}_{C_-\bs  \qq_-\sqcup P_\circ}(V)}=[1_{\A_\S}\o M^\circ(0)]$. Therefore,
	\begin{align*}
		&[M^\star \o \Phi^\L(\A_\S)]_{\mathcal{L}_{C_+\bs P_\star\sqcup \qq_+}(V)}\o_\A [\Phi^\sR(\A_\S)^\theta\o M^\circ]_{\mathcal{L}_{C_-\bs  \qq_-\sqcup P_\circ}(V)}\\
		&=\spn\{[v^\star \o 1_{\A_\S}]\o_\A [1_{\A_\S}\o v^\circ]: v^\star\in M^\star(0), v^\circ\in M^\circ(0)\}.
	\end{align*}
	It then follows from \eqref{eq:mutinodegeneral} that $[M^\star\o M^\circ]_{\mathcal{L}_{C\bs P_\bullet}(V)}=[M^\star(0)\o M^\circ (0)]$. This proves \eqref{eq:generalsurj}.

	By \eqref{eq:smoothsurj} and \eqref{eq:generalsurj}, the morphism $M^\bullet(0)\o \O_{\overline{\mathrm{M}}_{0,n}}\ra  \V(V;M^\bullet)$ between sheaves on $\overline{\mathrm{M}}_{0,n}$ is surjective on each fiber. Hence, it is an epimorphism of sheaves. Therefore, $ \V(V;M^\bullet)$ is globally generated.
\end{proof}
As a straightforward corollary, we obtain the global generation of coinvariant bundles associated with strongly rational VOAs generated in degree $1$ \cite{DG23}, which we record below.

\begin{corollary}\label{cor:degreeone}
	Let \(V\) be a strongly rational VOA such that \(V=U+C_1(V)\), where \(U=\C\vac+V_1\) (note that $L(1)V_1=0$, since $V\cong V'$). Then \(\S=\mathscr{W}\) is a strongly geometrically reducible set. Therefore, the coinvariant bundle \(\V(V;M^\bullet)\) on \(\overline{\mathrm{M}}_{0,n}\) is globally generated for every \(M^\bullet\subseteq\mathscr{W}^{\oplus n}\).
\end{corollary}

\section{Globally generated coinvariant bundles defined by general VOAs}
\label{sec: examples}
In all examples of this section, the subspace $U$ with $V=U+C_1(V)$ can be chosen to be spanned by quasi-primary vectors, i.e., $L(1)U=0$, as required in Theorem~\ref{thm:criterion1}.

\subsection{Lattice VOAs}\label{sec:3.3}
In this subsection, we show that a subset $\S$ of irreducible modules over various small-rank lattice VOAs is strongly geometrically reducible as in Definition~\ref{df:stronglyreducible}.  We refer to the Appendix Section~\ref{sec:appendix} for the basics of lattice VOAs. 

\subsubsection{Rank-one lattice VOAs}
Let $L=\Z\al$, with $(\al|\al)=2k$ for some fixed $k\in \Z_{\geq 1}$. Let $V_L$ be the associated lattice VOA. 
We can take
$$V_L=U+C_1(V_L),\quad U=\spn\{\vac,\ \al(-1)\vac,\ e^{\pm \al}\},$$
see \cite{KL99}. 
Note that $V_L$ is not generated by $(V_L)_1$ unless $k=1$. The irreducible modules are given by
\[
\mathscr{W}=\left\{V_L,V_{L+\frac{\al}{2k}},\ds, V_{L+\frac{\al}{2}},\ds, V_{L+\frac{(2k-1)\al}{2k}} \right\}.
\]
The fusion tensor products of these modules can be found in \cite{DL93} and are given as follows:
\begin{equation}\label{eq:rankonelatticefusion}
	V_{L+\frac{i\al}{2k}}\bt_{V_L} V_{L+\frac{j\al}{2k}}\cong V_{L+\frac{(i+j)\al}{2k}},\quad 0\leq i,j\leq 2k-1.
\end{equation}

\begin{proposition}\label{prop:rankonelattice}
	For the rank-one lattice VOA $V_L$, the subset of irreducible modules $\mathcal{S}=\{V_L, V_{L+\frac{1}{2}\al} \} $ is strongly geometrically reducible.
\end{proposition}
\begin{proof}
	For $M=V_{L}$, we have $M(0)=\C\vac$. Since $V_L$ is strongly generated by $U$ from $\C\vac$ and $b_{(n)}\vac=0$ for all $n\geq 0$, it is clear that $(\S2)$ \eqref{eq:Wspanning} and $(\S3)$ \eqref{eq:assumptioncriterion1} are satisfied.

	Now let $M=V_{L+\frac{1}{2}\al}$. Note that $M(0)=\spn\{e^{\pm \frac{\al}{2}}\}$. As a Heisenberg module, we can write
	\[
	V_{L+\frac{1}{2}\al}=\bigoplus_{n=0} ^\infty\mathscr{U}(M_{\hat{\h}}(1,0)).e^{\pm\left(n\al+\frac{\al}{2}\right)}=\bigoplus_{n=0}^\infty M_{\hat{\h}}\left(1,\pm(n\al+\al/2)\right).
	\]
	Therefore, in order to show that $V_{L+\frac{1}{2}\al}$ is strongly generated from its bottom degree by the VOA generating set $U_+=\spn\{\al(-1)\vac, e^{\pm \al}\}$, we only need to show that $e^{\pm\left(n\al+\frac{\al}{2}\right)}$, $n\geq 0$, are contained in
	\[
	N=\spn\{ b^1_{(-n_1)}\ds b^r_{(-n_r)} e^{\pm \frac{\al}{2}}: b^i\in \{\al(-1)\vac, e^{\pm \al}\},\ n_i\geq 0\}.
	\]
	This follows by induction on $n$ and from the fact that $e^{\pm\left((n+1)\al+\frac{\al}{2}\right)}=e^{\pm\al}_{-(2n+1)k-1} e^{\pm\left(n\al+\frac{\al}{2}\right)}$, in view of \eqref{latticeealpha}.
	To show \eqref{eq:assumptioncriterion1} for $b\in \{\al(-1)\vac, e^{\pm \al}\}$, we first note that, for any $n\geq 0$,
	$$(\al(-1)\vac)_{(n)} e^{\pm \frac{\al}{2}}=\al(n) e^{\pm\frac{\al}{2}}=\pm\delta_{n,0} k e^{\pm \frac{\al}{2}}\in M(0),$$
	in view of \eqref{latticeh}.
	Moreover, it follows from \eqref{latticeealpha} that $e^\al_{(n)} e^{\frac{\al}{2}}=0$ for any $n\geq 0$, and
	\[
	e^\al_{(k-1)} e^{-\frac{\al}{2}}= e^{\frac{\al}{2}},\quad e^\al_{(n)} e^{-\frac{\al}{2}}=0,\quad \forall n\geq k.
	\]
	For $0\leq m\leq k-1$, since the generator $\al(-1)\vac$ of the subVOA $M_{\hat{\h}}(1,0)$ is contained in $(V_L)_1$, it follows from \eqref{latticeealpha} that
	\[
	e^\al_{(m)} e^{-\frac{\al}{2}}=\Res_z  \exp\left(\sum_{n<0} \frac{-\al(n)}{n} z^{-n}\right) z^{-k+m} e^{\frac{\al}{2}}\in \mathscr{U}(M_{\hat{\h}}(1,0)).e^{\frac{\al}{2}}\ssq \mathscr{U}(\widehat{V_1}).e^{\frac{\al}{2}}.
	\]
	Therefore, $e^\al_{(n)} M(0)\ssq \mathscr{U}(\widehat{V_1}).M(0)$ for all $n\geq 0$. Similarly, we can show \eqref{eq:assumptioncriterion1} for $b=e^{-\al}$. Finally, it follows from \eqref{eq:rankonelatticefusion} that $\S$ is closed under fusion tensor product. 
\end{proof}

\begin{remark}
When $k\geq 2$,	the other irreducible $V_L$-modules $V_{L+\frac{i\al}{2k}}$, $i\neq k$, do not satisfy condition $(\S3)$.  For instance, let $W=V_{L+\frac{\al}{2k}}$. Then we have 
    $$W(0)=\C e^{\frac{\al}{2k}}, \quad  \mathscr{U}(\widehat{V_1}).W(0)=M_{\hat{\h}}(1,\al/2k);\quad \mathrm{But}\quad e^{-\al}_{(0)} e^{\frac{\al}{2k}}=e^{-\frac{(2k-1)\al}{2k}}\notin M_{\hat{\h}}(1,\al/2k).$$ Thus, the irreducible module $W$ does not satisfy \eqref{eq:assumptioncriterion1}.

In fact, conditions $(\S2)$ and $(\S3)$, given in \eqref{eq:Wspanning} and \eqref{eq:assumptioncriterion1}, respectively, were extracted from this lattice VOA example. In other words, they generalize the properties we used in the proof of Proposition~\ref{prop:rankonelattice}.
\end{remark}


\subsubsection{Rank-two lattice VOAs}
In general, the defining relations \eqref{latticeh}--\eqref{latticegeneral} of a lattice VOA are difficult to analyze for higher-rank lattices. Here we give two additional examples of $\S$ in rank-two lattice VOAs that are strongly geometrically reducible. 

\begin{example}\label{ex:latticesqrtKA2}
	Let $L=\sqrt{k}A_2=\Z\al\op \Z \b$, for some $k\in \Z_{\geq 2}$, where $A_2$ is the root lattice. Then
	\[
	(\al|\al)=(\b|\b)=2k,\quad (\al|\b)=-k,
	\]
	and $V_L=U+C_1(V_L)$, where $U=\C\vac+ (V_L)_1+\spn\{e^{\pm \al}, e^{\pm \b}, e^{\pm (\al+\b)} \}$, with $\wt (e^{\pm \al})=\wt (e^{\pm \b})=\wt (e^{\pm (\al+\b)})=k\geq 2$, i.e., $V_L$ is not generated by $(V_L)_1$. Consider the following irreducible modules:
	\[
	W^+=V_{L+\mu},\quad W^-=V_{L-\mu},\quad \mu=\frac{2\al+\b}{3}\in L^\circ.
	\]
	As in the previous example and Proposition~\ref{prop:rankonelattice}, we can show that
$W^\pm$ satisfy conditions \eqref{eq:Wspanning} and \eqref{eq:assumptioncriterion1}.  Moreover, by the simple current property of lattice VOA modules \cite{DL93}, we have
	\[
	W^{\pm}\bt_{V_L} W^{\mp}\cong V_L,\quad W^{+}\bt_{V_L} W^+\cong W^-,\quad W^{-}\bt_{V_L} W^-\cong W^+.
	\]
Thus, $\S=\{V_L, W^\pm \}$ is strongly geometrically reducible. 
\end{example}

\begin{example}\label{ex:latticeA1+A1}
	Let $L=\sqrt{k}(A_1\op A_1)=\Z\al\op \Z\b$, for some $k\in \Z_{\geq 2}$, where
	\[
	(\al|\al)=(\b|\b)=2k,\quad (\al|\b)=0.
	\]
	In this case, we may take $U=\C\vac+ (V_L)_1+\spn\{e^{\pm \al}, e^{\pm \b} \}$, and we let $\S=\{V_L, V_{L+\frac{\al}{2}}, V_{L+\frac{\b}{2}}, V_{L+\frac{\al+\b}{2}}\}.$
	Again, as in the proof of Proposition~\ref{prop:rankonelattice}, we can show that
all modules in $\S$ satisfy conditions \eqref{eq:Wspanning} and \eqref{eq:assumptioncriterion1}. Moreover, $\S$ is fusion closed since it forms the group $\Z/2\Z\times \Z/2\Z$ under the fusion product. 
Thus, $\S$ is strongly geometrically reducible.
\end{example}

\subsubsection{An example of $\S$ that is geometrically reducible but not strongly geometrically reducible}
\begin{example}\label{ex:latticenontrivial}
	Let $L=\Z\al \op \Z\b$, with
	\[
	(\al|\al)=2,\quad (\b|\b)=4,\quad (\al|\b)=1.
	\]
	Then we may take $U=\C\vac+ (V_L)_1+ \spn\{e^{\pm \b}, e^{\pm (\b-\al)}\}$, where $(V_L)_1=\h(-1)\vac\op \C e^{\pm \al}$ and $e^{\pm \b}, e^{\pm (\b-\al)}\in (V_L)_2$. In particular, $V_L$ is not generated by $(V_L)_1$. Consider the following two irreducible modules:
	\[
	M^+=V_{L+\la},\quad M^-=V_{L-\la},\quad \la=\frac{2\al+3\b}{7}\in L^\circ.
	\]
	As in the proof of Proposition~\ref{prop:rankonelattice}, using the lattice VOA relations \eqref{latticeh}--\eqref{latticegeneral}, we can show that $M^\pm$ satisfy conditions \eqref{eq:Wspanning} and \eqref{eq:assumptioncriterion1}. Hence, $\S=\{V_L, M^\pm\}$ is geometrically reducible by Theorem~\ref{thm:criterion1}.
However, the fusion tensor products satisfy 
	\[
	M^\pm \bt_{V_L} M^\mp \cong V_L,\quad M^\pm \bt_{V_L} M^\pm\cong V_{L\pm 2\la} \notin \S.
	\]
In other words, $\S=\{V_L,M^\pm\}$ does not satisfy fusion closedness $(\S1)$, hence it is NOT strongly geometrically reducible. 
\end{example}

\subsection{Affine simple-current extension VOAs }\label{sec:3.4}
We refer to the Appendix Section~\ref{sec:appendix} for the basics of affine VOAs and simple current extension. 

\subsubsection{The simple-current extension of $L_8(\sl_2)$}
Consider the type-$A_1$ level-$8$ affine VOA $V=L_8(\sl_2)$. It has $9$ irreducible modules
$L_8(i)=L_8\left(\sl_2,i\al/2\right), 0\leq i\leq 8,$
where $\al$ is the unique positive root, with $(\al|\al)=2$. 
In particular, $J=L_8(8)$ is the unique simple-current irreducible $V$-module, with
\begin{equation}\label{L8tensorsec4}
	L_8(8)\boxtimes L_8(i)= L_8(8-i),\quad 0\leq i\leq 8.
\end{equation}
By \eqref{eq:affineconformalwt}, the conformal weight of $J=L_8(8)$ is
$
\frac{(8\al/2+\al|8\al/2)}{2(h^\vee+k)}=\frac{(5\al|4\al)}{2(2+8)}=2\in \Z.
$
It follows from Lemma~\ref{thm:DLM} that the following direct sum space is a VOA:
\begin{equation}\label{eq:Ve}
	V_e=V\op J=L_8(\sl_2)\op L_8(8),
\end{equation}
which is the simple-current extension of the affine VOA $V=L_8(\sl_2)$.

Let $I\in I_V\fusion{J}{J}{V}$ be the unique nonzero intertwining operator such that
$
v^8_{(3)}u^8=\vac\in V_0,
$
where $u^8,v^8$ denote lowest-weight vectors of $J$.
Then the vertex operator $Y_{V_e}$ is given by
\begin{equation}\label{eq:YVsec4}
	Y_{V_e}((a,u),z)(b,v)=\left(Y_V(a,z)b+I(u,z)v, Y_J(a,z)v+e^{zL(-1)}Y_J(b,-z)u \right),
\end{equation}
where $(a,u), (b,v)\in V_e=V\op J$. Note that $V_e$ is not generated by $(V_e)_1=V_1=\sl_2$.

\begin{lemma}[{\cite[Corollary 5.8]{Li97}}]\label{prop:irrmodVsec4}
	The simple-current extension VOA $V_e=L_8(\sl_2)\op L_8(8)$ has the following irreducible modules up to isomorphism:
	\begin{align*}
		&V_e=L_8(\sl_2)\op L_8(8),&& X=L_8(2)\op L_8(6),\\
		&W^+=L_8(4),\ Y_{W^+}=Y_{L_8(4)}+\Y,&& W^-= L_8(4),\ Y_{W^-}=Y_{L_8(4)}-\Y ,
	\end{align*}
with conformal weights $h_{X}=1/5$ and $h_{W^\pm}=3/5$. 
	where $\Y\in I_V\fusion{J}{L_8(4)}{L_8(4)}$ is a normalized intertwining operator. In other words, $\mathscr{W}=\{V_e,X,W^+,W^-\}$. 
\end{lemma}

Moreover, the fusion tensor products among the irreducible $V_e$-modules are given by 
\begin{equation}\label{eq:Vefusion}
\begin{aligned}
		&	X\boxtimes_{V_e} X \cong V_e\oplus X\oplus W^+\oplus W^-,&& 	W^+\boxtimes_{V_e} W^- \cong X,\\
		&	X\boxtimes_{V_e} W^+ \cong X\oplus W^-, && X\boxtimes_{V_e} W^- \cong X\oplus W^+,\\
		&	W^+\boxtimes_{V_e} W^+ \cong V_e\oplus W^+, && 	W^-\boxtimes_{V_e} W^- \cong V_e\oplus W^-,
	\end{aligned}  
\end{equation}	
see Proposition~\ref{prop:fusionsimplecurrent} for the proof. 
\begin{proposition}\label{prop:4.7}
	Let $V_e$ be the affine simple current extension \eqref{eq:Ve}. Then the subset $\S=\{W^+,W^-\}$ satisfies conditions $(\S2)$ \eqref{eq:Wspanning} and $(\S3)$ \eqref{eq:assumptioncriterion1}. In particular, $\S$ is geometrically reducible.
\end{proposition}
\begin{proof}
	By the definition of $Y_{V_e}$ \eqref{eq:YVsec4}, it is clear that $V_e=U+C_1(V_e)$, where $U=\C\vac+V_1+J_2=\C\vac+\sl_2+J_2$, and $J_2=J(0)$ is the bottom degree of $J$. By the definition of $Y_{W^\pm}$ in Lemma~\ref{prop:irrmodVsec4}, $W^\pm$ are the irreducible modules $L_8(4)$ when viewed as $V=L_8(\sl_2)$-modules. Therefore, $W^\pm$ are strongly generated by $\sl_2\subset U_+$ from $W^\pm(0)$ \eqref{eq:Wspanning}.

	Moreover, since the lowest weight of $J$ is $2$, we have $(V_e)_1=V_1\op J_1=L_8(\sl_2)_1\cong \sl_2$. Therefore, $\scrU(\widehat{(V_e)_1}).W^\pm(0)=\scrU(\widehat{\sl}_2).W^\pm(0)=W^\pm$. In particular, $b_{(n)}.W^\pm(0)\ssq \scrU(\widehat{(V_e)_1}).W^\pm(0)$ for any $b\in U_+\cap (V_e)_k$, $k\geq 1$, and $n\geq 0$. This proves \eqref{eq:assumptioncriterion1}.
\end{proof}

In view of \eqref{eq:Vefusion},  $\S=\{W^+,W^-\}$ does not satisfy the fusion closedness $(\S1)$. Therefore, it is {\bf not} strongly geometrically reducible. We will give an example in subsection~\ref{subsubsec: non-example simple-current extension} to show that the coinvariants bundle $\V(V_e,M^\bullet)$ for $M^\bullet\subseteq \{W^+,W^-\}^n$ is not globally generated.  

\begin{example}\label{ex:VeGG}
Let $V_e$ be the affine simple current extension \eqref{eq:Ve}. Let $\S^+=\{V_e,W^+\}$ and $\S^-=\{V_e,W^-\}$. Then by Proposition~\ref{prop:4.7}, both $\S^\pm$ satisfy conditions $(\S2)$ and $(\S3)$. Moreover, by \eqref{eq:Vefusion}, we have 
\[
V_e\bt_{V_e}W^\pm=W^\pm\quad \mathrm{and} \quad W^\pm \bt_{V_e} W^\pm=V_e\op W^\pm.
\]
Hence both $\S^+$ and $\S^-$ are strongly geometrically reducible subsets of $V_e$-modules. 
\end{example}

\subsection{Tensor product VOAs and diagonal simple current extensions}\label{sec:3.5}

We refer to the Appendix Section~\ref{sec:appendix} for the basics of tensor product VOAs.
\begin{proposition}\label{prop:tensorproduct}
	Let $V^1$ and $V^2$ be strongly rational VOAs, and let $V=V^1\o V^2$. Let $\S^1\subset \mathscr{W}^1$ and $\S^2\subset \mathscr{W}^2$ be two subsets of strongly geometrically reducible  $V^1$ and $V^2$-modules, respectively. Then
	\[
	\S^1\bt \S^2:=\{M^{1}\o M^{2}: M^{1}\in \S^1,\ M^{2}\in \S^2 \}
	\]
	is a subset of irreducible $V^1\o V^2$-modules that is strongly geometrically reducible. 
\end{proposition}
\begin{proof}
	Let $V^i=U^i+C_1(V^i)$ for $i=1,2$. It follows from \eqref{eq:tensor} that
	\begin{equation}\label{eq:modetsensor}
		(a^1\o a^2)_{(n)} (b^1\o b^2)=\sum_{j\in \Z} a^1_{(j)}b^1\o a^2_{(n-j-1)} b^2,\quad \forall n\in \Z.
	\end{equation}
	In particular, since $\vac_{(n)}c=0$ for all $n\geq 0$, we have
$
	a^1\otimes a^2
	=(a^1\otimes\vac_2)_{(-1)}(\vac_1\otimes a^2)
	\in C_1(V^1\otimes V^2),
$ for $a^i\in V^i_+$. 
	Moreover, $
	C_1(V^1)\otimes\vac_2 $ and $
	\vac_1\otimes C_1(V^2)
	$ are contained in $C_1(V^1\otimes V^2)$.
	Hence, $$V^1\o V^2= (U^1\o \C \vac_2+ \C\vac_1\o U^2)+C_1(V^1\o V^2),$$
	and we can take $U=U^1\o \C \vac_2+ \C\vac_1\o U^2$ in Definition~\ref{df:stronglygen}. Let $M^i\in \S^i$ for $i=1,2$. It follows from \eqref{eq:modetsensor} that
	\[
	(U^1\o \C\vac_2)_{(-n)} (M^1(0)\o M^2(0))=\left((U^1)_{(-n)} M^1(0)\right)\o M^2(0),\quad \forall n\geq 0.
	\]
	Hence, $M^1\o M^2$ satisfies condition $(\S2)$  \eqref{eq:Wspanning} with respect to the subspace $U$.

	Furthermore, since $L(0)=L^1(0)\o \Id+\Id\o L^2(0)$, we have $V_1=V^1_1\o \C\vac_2+\C\vac_1\o V^2_1$. Then
	\[
	\scrU(\widehat{V_1}).(M^1(0)\o M^2(0))=\scrU(\widehat{V^1_1}).M^1(0)\o \scrU(\widehat{V^2_1}).M^2(0).
	\]
	Let $b^1\in U^1\cap V^1_k$ for some $k\geq 1$, and let $n\geq 0$. Then, by \eqref{eq:modetsensor} and the fact that $M^1$ satisfies condition \eqref{eq:assumptioncriterion1}, we have
	\[
	(b^1\o \vac_2)_{(n)}. (M^1(0)\o M^2(0))=(b^1_{(n)}M^1(0))\o M^2(0)\ssq \scrU(\widehat{V_1}).(M^1(0)\o M^2(0)),
	\]
	and similarly for $(\vac_1\o b^2)_{(n)}. (M^1(0)\o M^2(0))$. Therefore, for any $b=b^1\o \vac_2+\vac_1\o b^2\in U\cap V_k$ with $k\geq 1$, we have
	\[
	b_{(n)}.  (M^1(0)\o M^2(0))\ssq \scrU(\widehat{V_1}).(M^1(0)\o M^2(0)),\quad \forall n\geq 0.
	\]
	Thus, $M^1\o M^2$ satisfies condition $(\S3)$ \eqref{eq:assumptioncriterion1}.
    
The fusion tensor product of modules over a tensor product VOA $V^1\o V^2$ can also be calculated componentwise. More precisely, let $M^i,N^i\in \mathscr{W}^i$ for $i=1,2$. Then
\begin{equation}\label{eq:fusionten}
	\left(M^1\o M^2\right)\bt_{V^1\o V^2} \left(N^1\o N^2\right) \cong \left(M^1\bt_{V^1} N^1\right)\o \left(M^2\bt_{V^2} N^2\right),
\end{equation}
see \cite[Theorem 2.10]{ADL05} for more details. It follows that if $\S^i\subset \mathscr{W}^i$ for $i=1,2$ are both closed under fusion, then $\S^1\bt \S^2$ is also closed under fusion.
\end{proof}

\subsubsection{A diagonal simple-current extension VOA}
Let $V=L_k(\sl_2)$, where $k=2m$ for some $m\geq 1$. By \eqref{eq:sl2fusion}, $V$ has a simple-current module $J_m=L_k(2m)$.
Now let $\tilde{V}:=V\o V$ and $\tilde{J}_m:=J_m\o J_m$. Then $\tilde{J}_m$ is a module over the strongly rational VOA $\tilde{V}$ by Lemma~\ref{lm:tensor}, with $\tilde{J}_m\bt_{\tilde{V}} \tilde{J}_m\cong \tilde{V}$.

Hence, $\tilde{J}_m$ is a simple-current module over $\tilde{V}$, and the conformal weight $h_{\tilde{J}_m}=2 h_{J_m}=m$ is a positive integer. Then
\[
V^{(m)}:= \tilde{V}\op \tilde{J}_m=(V\o V)\op (J_m\o J_m)
\]
is a simple-current extension of the VOA $\tilde{V}$, and it is again strongly rational \cite{DLM97,HKL15}. We call $V^{(m)}$ the {\em diagonal simple-current extension of the affine VOA $L_{2m}(\sl_2)$}.

Now we construct an irreducible module over $V^{(m)}$. Consider the irreducible $\tilde{V}$-module $N_m:=V\o J_m$. By \eqref{eq:fusionten}, we have
$
\tilde{J}_m\bt_{\tilde{V}}N_m\cong (J_m\bt_V V)\o (J_m\bt_V J_m)\cong J_m\o V.$ Then, with respect to the embedding $\tilde{V}\hookrightarrow V^{(m)}$, we have
\begin{equation}\label{eq:defXm}
X^{(m)}:=\mathrm{Ind}^{V^{(m)}}_{\tilde{V}}(N_m)=\left(\tilde{V}\op \tilde{J}_m\right)\bt_{\tilde{V}} N_m\cong (V\o J_m)\op (J_m\o V),
\end{equation}
which is an ordinary module over the VOA $V^{(m)}$, with conformal weight $h_{X^{(m)}}=m/2$.

Moreover, $X^{(m)}$ is irreducible since both $(V\o J_m)$ and $(J_m\o V)$ are irreducible $\tilde{V}=V\o V$-modules, and $\tilde{J}_m.(V\o J_m)=J_m\o V$ and $\tilde{J}_m. (J_m\o V)=(V\o J_m)$ \cite{HKL15}.

\begin{proposition}\label{prop:digonalsimplecurrect}
	Let $V^{(m)}$ be the diagonal simple-current extension VOA, with $m\geq 2$, and let $\S=\{V^{(m)},X^{(m)}\}$ be a subset of irreducible $V^{(m)}$-modules. Then $\S$ is strongly geometrically reducible. 
\end{proposition}
\begin{proof}
	Since the conformal weight of $\tilde{J}_m$ is $m\geq 2$, we have $V^{(m)}_1=\tilde{V}_1\op (\tilde{J}_m)_1=\tilde{V}_1\cong \sl_2\o \C\vac_2+\C\vac_1\o \sl_2$. Moreover, by \eqref{eq:defXm}, we have
	\[
	X^{(m)}(0)\cong \left((V\o J_m)\op (J_m\o V)\right)(0)=(\C\vac_1\o J_m(0))\op (J_m(0)\o\C\vac_2).
	\]
	Since $J_m=\scrU(\widehat{\sl_2}).J_m(0)=\scrU(\widehat{V_1}).J_m(0)$, by \eqref{eq:modetsensor}, we have
	\begin{align*}
		X^{(m)}&=\left(\scrU(\widehat{\C\vac_1 \o \sl_2}).(\C\vac_1\o J_m(0))\right)\op\left( \scrU(\widehat{\sl_2\o \C \vac_2}). (J_m(0)\o \C\vac_2)\right)=\scrU(\widehat{V^{(m)}_1}). X^{(m)}(0).
	\end{align*}
	Therefore, \eqref{eq:Wspanning} and \eqref{eq:assumptioncriterion1} are clearly satisfied.
	Furthermore, as a $\tilde{V}$-module, we have
	\[
	N_m\bt_{\tilde{V}} N_m=(V\o J_m)\bt_{\tilde{V}} (V\o J_m)\cong (V\bt_V V)\o (J_m\bt_V J_m)\cong V\o V=\tilde{V}.
	\]
	Applying the induction functor $\mathrm{Ind}^{V^{(m)}}_{\tilde{V}}(-)$ to both sides, we have
	\begin{align*}
		X^{(m)}\bt_{V^{(m)}} X^{(m)}&=\mathrm{Ind}^{V^{(m)}}_{\tilde{V}}(N_m)\bt_{V^{(m)}} \mathrm{Ind}^{V^{(m)}}_{\tilde{V}}(N_m)\cong \mathrm{Ind}^{V^{(m)}}_{\tilde{V}}(N_m\bt_{\tilde{V}}N_m)\\
		&\cong \mathrm{Ind}^{V^{(m)}}_{\tilde{V}}(\tilde{V})=V^{(m)}.
	\end{align*}
	Therefore, $\S=\{V^{(m)}, X^{(m)}\}$ is closed under the fusion tensor product.
\end{proof}



\subsection{Examples of globally generated coinvariant bundles over general VOAs}
\label{exmp: list of examples}
We list below new examples of globally generated coinvariant bundles arising from strongly rational VOAs $V$ that are not generated by $V_1$.
	\begin{enumerate}
		\item Let $V=V_L$ be a rank-one lattice VOA, with $L=\Z\al$ and $(\al|\al)=2K$ for some $K\in \Z_{\geq 2}$, and let $\S=\{V_{L},V_{L+\frac{\al}{2}}\}$. Then, by Proposition~\ref{prop:rankonelattice} and Theorem~\ref{thm:main1}, $\V(V_L;M^\bullet)$ is globally generated on $\overline{\mathrm{M}}_{0,n}$ for all $n\geq 3$ and any choice of $M^\bullet\subset\S$.
		\item Let 	$L=\sqrt{K}A_2=\Z\al\op \Z \b$, for some $K\in \Z_{\geq 2}$, with
		\[
		(\al|\al)=(\b|\b)=2K,\quad (\al|\b)=-K,
		\]
		and let $V=V_L$ be the associated lattice VOA. Let $\S=\{V_{L}, V_{L\pm \mu} \}$, where $\mu=\frac{2\al+\b}{3}$.
		Then, by Example~\ref{ex:latticesqrtKA2} and Theorem~\ref{thm:main1}, the coinvariant bundle $\V(V_L;M^\bullet)$ is globally generated on $\overline{\mathrm{M}}_{0,n}$ for all $n\geq 3$ and any choice of $M^\bullet\subset \S$.
		\item Let $L=\sqrt{K}(A_1\op A_1)=\Z\al\op \Z\b$, for some $K\in \Z_{\geq 2}$, with
		\[
		(\al|\al)=(\b|\b)=2K,\quad (\al|\b)=0,
		\]
		and let $V=V_L$ be the associated lattice VOA.
		Let $\S=\left\{V_L, V_{L+\frac{\al}{2}}, V_{L+\frac{\b}{2}}, V_{L+\frac{\al+\b}{2}} \right\}.$ Then, by Example~\ref{ex:latticeA1+A1} and Theorem~\ref{thm:main1}, the coinvariant bundle $\V(V_L;M^\bullet)$ is globally generated on $\overline{\mathrm{M}}_{0,n}$ for all $n\geq 3$ and any choice of $M^\bullet\subset \S$.
    \item Let $V_e=L_8(\sl_2)\op L_8(8)$ be the simple current extension of the affine VOA $L_8(\sl_2)$, and let $\S^+=\{V_e,W^+\}$ and $\S^-=\{V_e,W^-\}$. Then, by Example~\ref{ex:VeGG} and Theorem~\ref{thm:main1}, the bundles $\V(V_e;(W^+)^\bullet)$ and $\V(V_e;(W^-)^\bullet)$ are both globally generated over $\overline{M}_{0,n}$.

Moreover, it is easy to show using the fusion tensor products and the factorization theorem (Theorem~\ref{thm: factorization}) that
\[
\mathrm{rank}\,\V(V_e;(W^\pm)^{\otimes n})=F_{n-1},
\]
where the Fibonacci numbers are defined by $F_0=0$, $F_1=1$, and $F_{m+1}=F_m+F_{m-1}$.
In particular, these globally generated bundles are not line bundles for $n\geq 4$.
		\item Let $V^{(m)}$ be the diagonal simple-current extension VOA $\tilde{V}\op \tilde{J}_m$ associated with the affine VOA $V=L_{2m}(\sl_2)$ and $J_m=L_{2m}(2m)$, where $m\geq 2$, $\tilde{V}=V\o V$, and $\tilde{J}_m=J_m\o J_m$. Let $\S=\{ V^{(m)}, X^{(m)}\}$, where $X^{(m)}=\mathrm{Ind}^{V^{(m)}}_{\tilde{V}}(V\o J_m)$.

		Then, by Proposition~\ref{prop:digonalsimplecurrect} and Theorem~\ref{thm:main1}, the coinvariant bundle $\V(V^{(m)}; M^\bullet)$ is globally generated on $\overline{\mathrm{M}}_{0,n}$ for all $n\geq 3$ and any choice of $M^\bullet\subset\S$.

		\item Let $V$ be the tensor product VOA $V^1\o V^2$, where $V^1$ and $V^2$ are chosen from the VOAs in the examples above or are affine VOAs of positive integral level $k\in \Z_{\geq 1}$:
		\[
		V_{\Z\al},\quad V_{\sqrt{K}A_2},\quad V_{\sqrt{K}(A_1\op A_1)},\quad V^{(m)},\quad L_k(\g),
		\]
		and let $\S^1$ and $\S^2$ be the corresponding strongly geometrically reducible subsets as in items $(1)$--$(5)$ above. When $V^i=L_k(\g)$, let $\S^i=\mathscr{W}$ be the full set of irreducible modules. Then, by Proposition~\ref{prop:tensorproduct} and Theorem~\ref{thm:main1}, the coinvariant bundle $\V(V^1\o V^2; M^\bullet)$ is globally generated on $\overline{\mathrm{M}}_{0,n}$ for all $n\geq 3$ and any choice of $M^\bullet\subset \S=\S^1\bt \S^2$.
	\end{enumerate}

\subsection{Fusion closure is essential for global generation}
\label{subsubsec: non-example simple-current extension}
By Definition~\ref{df:stronglyreducible}, a strongly geometrically reducible subset $\S$ must be closed under the fusion tensor product, as required by $(\S1)$. We now use the simple-current extension VOA $V_e=L_8(\sl_2)\op L_8(8)$ of Section~\ref{sec:3.4} to show that this condition cannot be dropped from our criterion for the global generation of coinvariant bundles.

By Theorem~\ref{thm:criterion1} and Proposition~\ref{prop:fusionsimplecurrent}, $\S=\{W^+,W^- \}$ is geometrically reducible in the sense of Definition~\ref{df:geored}, but it is {\bf not} closed under the fusion tensor product $(\S1)$. We will show that $\mathbb{V}:=\mathbb{V}_{0,8}(V_e,(W^+)^{\otimes4}\otimes(W^-)^{\otimes 4})$ is not globally generated on $\overline{\mathrm M}_{0,8}$. If $\mathbb{V}$ is globally generated, then its first Chern class $c_1\mathbb{V}$ is semi-ample and therefore intersects all curves on $\overline{\mathrm M}_{0,8}$ non-negatively. We produce an $F$-curve that $\mathbb{V}$ intersects strictly negatively, thereby proving that $\mathbb{V}$ is not globally generated. Let us recall the fusion rules needed for this computation:
\[
W^+\boxtimes W^- = X;\qquad X\boxtimes X = V_e\oplus X\oplus W^+\oplus W^-.
\]

The specific $F$-curve we use is given by the partition $\bigsqcup_{i=1}^4\{i,i+4\}$ of the set $\{1,2,\ldots,8\}$. The intersection number between the $F$-curve and the divisor $c_1\mathbb{V}$ is the sum
\[
\sum_{S_1,\ldots,S_4\in \mathscr{W}}\left(\deg\mathbb{V}_{0,4}(V_e;S_1\otimes\cdots\otimes S_4) \;\prod_{i=1}^4\text{rank} \mathbb{V}_{0,3}(V_e;W^+\otimes W^-\otimes S'_i)\right),
\]
where the modules $S_i$ run over the set of all simple $V_e$-modules $\mathscr{W}:=\{V_e,X,W^{+},W^-\}$. Since $W^+\boxtimes W^-=X$, and $X$ is self-dual, we see that 
\[
\text{rank}\mathbb{V}_{0,3}(V_e;W^+\otimes W^-\otimes S'_i) = 1 \text{ if and only if } S_i=X\text{ and it is zero for all other $S_i\in\mathscr{W}$}.
\]
Therefore, the sum above reduces to a single summand $\deg\mathbb{V}_{0,4}(V_e;X^{\otimes4})$, which is equal to
\[
4\;\text{cw}(X)\;\text{rank}\mathbb{V}_{0,4}(V_e;X^{\otimes4}) - 3\sum_{S\in\mathscr{W}}\text{cw}(S)\;\text{rank}\mathbb{V}_{0,3}(V_e;X^{\otimes2}\otimes S)\;\text{rank}\mathbb{V}_{0,3}(V_e;X^{\otimes2}\otimes S').
\]
As $X\boxtimes X=\bigoplus_{S\in\mathscr{W}}S$, we have $\text{rank}\mathbb{V}_{0,3}(V_e; X^{\otimes2}\otimes S)=1$ for each $S\in\mathscr{W}$. Therefore,
\begin{align*}
\text{rank}\mathbb{V}_{0,4}(V_e;X^{\otimes 4}) = \sum_{S\in\mathscr{W}}\text{rank}\mathbb{V}_{0,3}(V_e;X^{\otimes 2}\otimes S)\;\text{rank}\mathbb{V}_{0,3}(V_e;X^{\otimes 2}\otimes S') = 4,
\end{align*}
and we see that $(c_1\mathbb{V})\cdot F =  \deg\mathbb{V}_{0,4}(V_e;X^{\otimes4})$ is equal to
\begin{align*}
16\;\text{cw}(X) - 3\;\sum_{S\in\mathscr{W}}\text{cw}(S)= 13\; \text{cw}(X)-3\left(\text{cw}(W^+)+\text{cw}(W^-)\right)= \frac{13}{5} - 6\cdot\frac{3}{5} = -1.
\end{align*}

\section{Cohomological approach to global generation}
\label{sec: cohom approach}

\subsection{Proof of Theorem A}
\label{subsec: proof of Theorem A}

In this section, we assume that all vertex operator algebras are strongly rational and of CohFT-type in the sense of \cite{DGT3, DG23}.
\begin{notation}
Let $\mathcal{S}^1$ and $\mathcal{S}^2$ be subrings of the fusion rings of the vertex operator algebras $V^1$ and $V^2$, respectively. Assume that $f:\mathcal{S}^1\hookrightarrow\mathcal{S}^2$ defines a proportional pairing with proportionality $k\in\mathbb{Z}_{\geq1}$. We denote the coinvariant bundle on $\overline{\mathrm{M}}_{0,n}$ with simple $V^1$-modules $M^1,\ldots,M^n\in\mathcal{S}^1$ by $\mathbb{V}_{0,n}(V^1,\otimes_{i=1}^n M^i)$ and its determinant bundle is denoted by $\mathbb{D}_{0,n}(V^1,\otimes_{i=1}^n M^i)$. Moreover, the image $f(M^i)\in\mathcal{S}^2$ is a simple $V^2$-module and we denote the corresponding coinvariant bundle by $\mathbb{V}_{0,n}(V^2,\otimes_{i=1}^n f(M^i))$. Finally, denoting  $\cw(M^i):=a_i$, we have $\cw(f(M^i))=ka_i$.
\end{notation}
\begin{theorem}
\label{thm: K-class result}
With the notation as above,
\(
[\mathbb{V}_{0,n}(V^2,\otimes_{i=1}^n f(M^i))] = [\psi^k\mathbb{V}_{0,n}(V^1,\otimes_{i=1}^n M^i)]
\)
in $K_0(\overline{\mathrm{M}}_{0,n})\otimes\mathbb{Q}$, where $\psi^k:K_0(\overline{\mathrm{M}}_{0,n})\otimes\mathbb{Q}\to K_0(\overline{\mathrm{M}}_{0,n})\otimes\mathbb{Q}$ is the $k$-th Adams ring endomorphism.
\end{theorem}
\begin{proof}
By \cite{DG23}, the Chern character of $\mathbb{V}_2:=\mathbb{V}_{0,n}(V^2,\otimes_{i=1}^n f(M^i))$ is 
\[
\sum_{\Gamma,\mu}\frac{1}{|\mathrm{Aut}\Gamma|} {(\xi_{\Gamma})}_* C^k_{\Gamma,\mu}\ ,\quad \text{ where }
\quad
C^k_{\Gamma,\mu} = \prod_{i=1}^n e^{ka_i\psi_i} \prod_{v\in V(\Gamma)}\mathrm{cont}_{\mu}(v) \prod_{e\in E(\Gamma)}\frac{1-e^{ka_{\mu}\psi_e}}{\psi_e},
\]
$\psi_e=\psi_h+\psi_{h'}$, and the sum runs over all stable graphs $\Gamma$ and module arrangements $\mu$.  A priori, the modules $\mu$ runs over all simple $V^2$-modules. However, we can restrict to simple modules $\mu$ in $\mathcal{S}^2$ by Proposition~\ref{prop: factorization respects subrings}, as $\mathrm{cont}_{\mu}(v):=\rank\ \mathbb{V}(V^2,\otimes_{i\in L(v)}f(M^i)\otimes M^\mu)$ is zero if $\mu\notin\mathcal{S}^2$. In fact, we can further restrict to $\mu\in f((\mathcal{S}^1)^{\mathrm{sim}})$ by definition of proposition pairing \ref{def: pairing}. Here, $L(v)$ is the set of legs attached to the vertex $v$ of the stable graph $\Gamma$, and $M^\mu$ denotes the tensor product of the modules attached by $\mu$ to the half-edges at $v$. Since ranks of bundles of coinvariants on $\overline{\mathrm{M}}_{0,m}$ are determined by the fusion rules via the factorization theorem, and $f$ is an injective ring homomorphism with $f((\mathcal{S}^1)^{\mathrm{sim}})\subseteq(\mathcal{S}^2)^{\mathrm{sim}}$, the vertex contributions $\mathrm{cont}_\mu(v)$ are the same for $\mathbb{V}_2$ and for $\mathbb{V}_1:=\mathbb{V}_{0,n}(V^1,\otimes_{i=1}^n M^i)$. The $t$-th Chern character $\ch_t(\mathbb{V}_2)\in H^{2t}(\overline{\mathrm{M}}_{0,n})$ is defined as
\[
\ch_t(\mathbb{V}_2) = \sum_{q=0}^t \sum_{\substack{\Gamma,\mu\\ |E(\Gamma)|=q}}\frac{(-1)^q}{|\mathrm{Aut}\Gamma|} B^k_{\Gamma,\mu},
\]
where 
{\footnotesize  \[ 
B^k_{\Gamma,\mu} = (\xi_{\Gamma})_*\left(\prod_{v\in V(\Gamma)}\mathrm{cont}_{\mu}(v)\left(\sum_{\substack{x,y\\ x+y=t-q}}\left(\sum_{\substack{s_1,\ldots,s_n\in\nn\\ s_1+\cdots+s_n=x}}\ \prod_{p=1}^n\frac{(ka_p)^{s_p}}{{s_p}!}\psi_p^{s_p}\right)
\left(\sum_{\substack{z_1,\ldots,z_q\in\nn\\ z_1+\cdots+z_q=y}}\ \prod_{p=1}^q(ka_{\mu_p})^{z_p+1}\frac{\psi_{e_p}^{z_p}}{{(z_p+1)}!}\right)
\right)\right).
\]}

Upon some calculation, we find that $B^k_{\Gamma,\mu}=k^t B^1_{\Gamma,\mu}$. Therefore, $\ch_t(\mathbb{V}_2)= k^t \ch_t(\mathbb{V}_1)$ in $H^{2t}(\overline{\mathrm{M}}_{0,n};\mathbb{Q})$. On the other hand, $\ch_t(\psi^k\mathbb{V}_1)=k^t\,\ch_t(\mathbb{V}_1)$ by the defining property of the Adams operations. Since the Chern character $K_0(\overline{\mathrm{M}}_{0,n})\otimes\mathbb{Q}\to H^{2\ast}(\overline{\mathrm{M}}_{0,n};\mathbb{Q})$ is an isomorphism, the conclusion follows.
\end{proof}
\begin{proposition}
\label{prop: tensor VB}
If $\mathbb{V}_{0,n}(V^1,\otimes_{i=1}^nM^i)$ has rank $1$, then so does $\mathbb{V}_{0,n}(V^2,\otimes_{i=1}^nf(M^i))$. Moreover,
\[
\mathbb{V}_{0,n}(V^2,\otimes_{i=1}^nf(M^i))\cong \mathbb{V}_{0,n}(V^1,\otimes_{i=1}^nM^i)^{\otimes k}
\]
as vector bundles on $\overline{\mathrm{M}}_{0,n}$.
\end{proposition}
\begin{proof}
Write $\mathbb{L}_1:=\mathbb{V}_{0,n}(V^1,\otimes_{i=1}^nM^{i})$ and $\mathbb{L}_2:=\mathbb{V}_{0,n}(V^2,\otimes_{i=1}^nf(M^{i}))$. The rank of a bundle of coinvariants on $\overline{\mathrm{M}}_{0,n}$ is determined by the fusion rules through the factorization theorem (Theorem~\ref{thm: factorization}) and Proposition~\ref{prop: factorization respects subrings}. Since the proportional pairing $f$ is an injective ring homomorphism of fusion subrings taking simple modules to simple modules, the two computations agree, and $\mathbb{L}_2$ has rank $1$. By Theorem~\ref{thm: K-class result}, $c_1(\mathbb{L}_2)=\ch_1(\mathbb{L}_2)=k\ \ch_1(\mathbb{L}_1)=k\ c_1(\mathbb{L}_1)$ in $H^2(\overline{\mathrm{M}}_{0,n},\mathbb{Q})$. Since $\overline{\mathrm{M}}_{0,n}$ is a rational smooth projective variety, we have $\mathrm{Pic}^0\left(\overline{\mathrm{M}}_{0,n}\right)=0$. Therefore, $c_1:\mathrm{Pic}\left(\overline{\mathrm{M}}_{0,n}\right)\to \mathrm{NS}\left(\overline{\mathrm{M}}_{0,n}\right)$ is a group  isomorphism and the conclusion follows.
\end{proof}
\begin{corollary}
\label{cor: main corollary}
With the notation as above, $\mathbb{D}_{0,n}(V^2,\otimes_{i=1}^n f(M^i))\cong \mathbb{D}_{0,n}(V^1,\otimes_{i=1}^n M^i)^{\otimes k}$ as line bundles on $\overline{\mathrm{M}}_{0,n}$. Moreover, if $\mathbb{D}_{0,n}(V^1,\otimes_{i=1}^n M^i)$ is globally generated, so is $\mathbb{D}_{0,n}(V^2,\otimes_{i=1}^n f(M^i))$.
\end{corollary}
\begin{proof}
It follows from Proposition~\ref{prop: tensor VB} since $c_1\mathbb{D}_{0,n}(V;M^{\bullet}) = c_1\mathbb{V}_{0,n}(V;M^{\bullet})$ for $n$-modules  $M^{\bullet}$ of $V$.
\end{proof}


\subsection{Globally generated bundles from type $A$ even lattice VOAs}
\label{subsec: typeA}
In this section, we study global generation of coinvariant bundles arising from representation theory of certain type $A$ even lattice VOAs. The lattices under consideration form an infinite family indexed by pairs $(k,r)\in \mathbb{Z}_{\geq 1}^2$. Let us denote the family by $\{L(k,r)\mid k,r\in\zz_{\geq1}\}$. For each fixed pair $(k,r)$, the corresponding lattice VOA $V_{L(k,r)}$ is generated in degrees $1$ and $k$. In Proposition~\ref{prop: subring str}, we identify a subring $\mathcal{S}_r$ of the fusion ring so that any coinvariant vector bundle with modules in $\mathcal{S}_r$ is globally generated (cf. Theorem~\ref{thm: global generation}). 

Fix integers $r,k\geq1$. Let $B:=k\cdot C(A_r)$, where $C(A_r)$ is the Cartan matrix corresponding to the Lie algebra $A_r$. Consider the rank $r$ lattice vertex operator algebra $V_L$ corresponding the lattice $L:=\oplus_{i=1}^r\zz\gb_i$ with the bi-linear pairing $(\cdot|\cdot):L\times L\to\zz$ given by the matrix $B$. Since $\det(B)=k^r(r+1)$ and for any integers $m,n\in\zz$ and $1\leq i,j\leq r$, we have
\[
(m\gb_i+n\gb_j \mid m\gb_i +n\gb_j)
= m^2 (\gb_i | \gb_i) +n^2(\gb_j| \gb_j)+2mn(\gb_i | \gb_j)
= 2\left(km^2+kn^2+mn(\gb_i | \gb_j)\right),
\]
$L$ is a non-degenerate even lattice. Therefore, $V_L$ is a strongly rational VOA. Let $L^{\circ}$ be the dual lattice of $L$ and let $L+\varphi_1,\ldots,L+\varphi_r\in L^{\circ}/L$ be defined by
\begin{equation}
\label{eqn: varphi}
\varphi_i := \frac{1}{k(r+1)}\left((r+1-i)\sum_{j=1}^i j\gb_j+i\sum_{j=i+1}^r(r+1-j)\gb_j\right).
\end{equation}
\begin{proposition}\label{prop: subring str}
The set $S_r:=\{V_L, V_{L+k\varphi_1},\ldots,V_{L+k\varphi_r}\}$
generates a subring $\mathcal{S}_r$ of the fusion ring of $V_L$.
\end{proposition}
The proof of Proposition~\ref{prop: subring str} is given after Theorem~\ref{thm: global generation}, the main result of this section. First, we discuss two important results that are needed to prove Theorem~\ref {thm: global generation}.
\begin{lemma}
The conformal weight of the irreducible $V_L$-module $V_{L+k\varphi_i}$ is $\frac{ki(r+1-i)}{2(r+1)}$ for each $1\leq i\leq r$.
\end{lemma}
\begin{proof}
The conformal weight of $V_{L+k\varphi_i}$ is given by
\begin{align*}
\frac{1}{2}\min_{a_1,\ldots,a_r\in\zz}\left(k\varphi_i + \sum_{j=1}^ra_j\gb_j \mid k\varphi_i + \sum_{j=1}^ra_j\gb_j\right)
&=\frac{1}{2}(k\varphi_i | k\varphi_i) = \frac{ki(r+1-i)}{2(r+1)}.
\end{align*}
\end{proof}
\begin{proposition}\label{prop: main prop}
Let $\mathcal{R}$ be the fusion ring of the affine VOA $L_{\hat{\mathfrak{sl}}_{r+1}}(1,0)$. Then $(\mathcal{R},\mathcal{S}_r)$ forms a proportional pairing with proportionality $k$.
\end{proposition}
\begin{proof}
We define a map $f:\mathcal{R}\to\mathcal{S}_r$ by sending $L_{\hat{\mathfrak{sl}}_{r+1}}(1,\varpi_i)\mapsto V_{L+k\varphi_i}$ for each $1\leq i\leq r$ and $L_{\hat{\mathfrak{sl}}_{r+1}}(1,0)\mapsto V_L$. Since both fusion rings are the group ring of $\Z/(r+1)\Z$ (see the proof of Proposition~\ref{prop: subring str} below), $f$ is a ring isomorphism. Moreover, for each $1\leq i\leq r$, we have
\[
\frac{\cw(V_{L+k\varphi_i})}{\cw(L_{\hat{\mathfrak{sl}}_{r+1}}(1,\varpi_i))} = \left(\frac{ki(r+1-i)}{2(r+1)}\right)/\left(\frac{i(r+1-i)}{2(r+1)}\right) = k.
\]
Finally, $\cw(V_L)=0=k\cdot 0 = k\cdot \cw(L_{\hat{\mathfrak{sl}}_{r+1}}(1,0))$. This verifies the conditions of Definition~\ref{def: pairing}.
\end{proof}
\begin{theorem}\label{thm: global generation}
Any coinvariant bundle on $\overline{\mathrm{M}}_{0,n}$ associated to irreducible $V_L$-representations in the subring $\mathcal{S}_r$ is globally generated.
\end{theorem}
\begin{proof}
Pick any integers $t_i\geq0$ with $t_1+\cdots+t_r=n$, and let
\[
\mathbb{V}_1:=\mathbb{V}_{0,n}(L_{\hat{\mathfrak{sl}}_{r+1}}(1,0),\otimes_{i=1}^r L_{\hat{\mathfrak{sl}}_{r+1}}(1,\varpi_i)^{\otimes t_i}) \text{ and }\mathbb{V}_2:=\mathbb{V}_{0,n}(V_L,\otimes_{i=1}^r V_{L+k\varphi_i}^{\otimes t_i})
\]
be the coinvariant vector bundles corresponding to representations of $L_{\hat{\mathfrak{sl}}_{r+1}}(1,0)$ and of $V_L$, respectively; by Propagation of Vacua (Theorem~\ref{thm: POV}) we may assume that none of the marked modules is the VOA itself. Note that the rank of $\mathbb{V}_2$ is either $0$ or $1$ and is equal to $1$ if and only if $k\sum_{i=1}^r t_i\varphi_i\in L$ (by Proposition~\ref{prop: non-triviality condition}). A similar result holds for $\mathbb{V}_1$ by \cite{Fak12}.
By Proposition~\ref{prop: main prop}, $\mathbb{V}_2$ is obtained from $\mathbb{V}_1$ by applying the proportional pairing $f$ with proportionality $k$. Hence, by Proposition~\ref{prop: tensor VB}, we have
\(
c_1\mathbb{V}_2 = k\ \cdot\  c_1\mathbb{V}_{1}
\)
and $\mathbb{V}_2\cong\mathbb{V}_1^{\otimes k}$. Since $\mathbb{V}_1$ is globally generated by Theorem~\ref{thm: Fakh}, so is $\mathbb{V}_2$.
\end{proof}

\begin{proof}[Proof of Proposition~\ref{prop: subring str}]
The elements in the set $S_r$ are pairwise inequivalent irreducible representations. Indeed, for any $1\leq i\leq r$, $V_{L}\cong V_{L+k\varphi_i}$ as $V_L$-representations if and only if $k\varphi_i\in L$. But $k\varphi_i\notin L$ since the coefficient of $\gb_1$ in $k\varphi_i$ is $\frac{r+1-i}{r+1}\notin\zz$. Therefore, $V_L$ and $V_{L+k\varphi_i}$ are inequivalent irreducible $V_L$-representations. Similarly, for any $1\leq i<j\leq r$, $V_{L+k\varphi_i}\cong V_{L+k\varphi_j}$ if and only if $k\varphi_i-k\varphi_j\in L$. Again the coefficient of $\gb_1$ in $k(\varphi_i-\varphi_j)$ is $\frac{r+1-i}{r+1}-\frac{r+1-j}{r+1}=\frac{j-i}{r+1}\notin\zz$, as $0<j-i<r+1$. Therefore, $V_{L+k\varphi_i}$ and $V_{L+k\varphi_j}$ for $i\ne j$ are inequivalent irreducible $V_L$-representations. Define $\mathcal{S}_r:=\bigoplus_{i=1}^r\zz V_{L+k\varphi_i}\oplus\zz V_L$ to be a formal group. To prove that $\mathcal{S}_r$ is a subring of the fusion ring of $V_L$, it suffices to prove that the set $S_r$ is closed under the fusion product. We claim that
\[
V_{L+k\varphi_a}\boxtimes V_{L+k\varphi_b} = V_{L+k\varphi_{a+b}}\quad \text{ for any integers $1\leq a,b\leq r$},
\]
where we set $\varphi_{r+1}:=0$ and $\varphi_{s}:=\varphi_{s-(r+1)}$ for $s>r+1$. This amounts to proving the following three cases: 
\begin{itemize}
\item If $a+b\leq r$, we have $k(\varphi_a+\varphi_b-\varphi_{a+b})=0\mod L$,
\item If $a+b=r+1$, we have $k(\varphi_a+\varphi_b)=0\mod L$, and
\item If $r+1<a+b\leq 2r$, we have $k(\varphi_a+\varphi_b-\varphi_{a+b-(r+1)})=0\mod L$.
\end{itemize}
For the second case, one sees $k(\varphi_a+\varphi_b) = \sum_{j=1}^r\min(a,b,j,r+1-j)\gb_j\in L$ and we leave the proof for the other two cases to the reader.
Therefore, we see that $\mathcal{S}_r$ is a subring of simple currents, which completes the proof.
\end{proof}
\begin{proposition}
\label{prop: non-triviality condition}
Let $k\geq1$ be an integer. Then the coinvariant bundle $\mathbb{V}:=\mathbb{V}_{0,n}(V_L,\bigotimes_{i=1}^n V_{L+k\varphi_{j_i}})$, where $1\leq j_i\leq r$ are integers for each $i$, has non-zero rank if and only if $k\sum_{i=1}^n\varphi_{j_i}\in L$.
Moreover, $\mathbb{V}$ is non-trivial if and only if there is a partition $I_1\sqcup I_2\sqcup I_3\sqcup I_4=\{1,\ldots,n\}$ so that the following holds: for $\nu_p:=\sum_{i\in I_p}j_i\mod (r+1)$,
\begin{enumerate}[i.]
\item $\sum_{p=1}^4\nu_p=2(r+1)$.
\item If $\left(\max_p\nu_p\right)+\left(\min_p\nu_p\right)<r+1$, then $\min_p\nu_p\ne0$.
\end{enumerate}
\end{proposition}
\begin{proof}
By the Propagation of Vacua (Theorem~
\ref{thm: POV}), we only need to consider bundles where no $V_L$-modules are $V_L$ itself. Using Factorization Theorem (Theorem~\ref{thm: factorization}), we see that
\begin{align*}
\rank \mathbb{V} = \sum_{M} \rank \mathbb{V}_{0,n-1}(V_L,\bigotimes_{i=1}^{n-2}V_{L+k\varphi_{j_i}}\otimes M)\cdot\rank \mathbb{V}_{0,3}(V_L,V_{L+k\varphi_{j_{n-1}}}\otimes V_{L+k\varphi_{j_{n}}}\otimes M'),
\end{align*}
where $M$ is running over all simple $V_L$-representations. By Proposition~\ref{prop: factorization respects subrings}, it suffices to consider simple $V_L$-representations in the subring $\mathcal{S}_r$. Therefore, $M=V_{L+k\varphi_t}$ for some $1\leq t\leq r$ or $M=V_L$. However, $\rank \mathbb{V}_{0,3}(V_L,V_{L+k\varphi_{j_{n-1}}}\otimes V_{L+k\varphi_{j_{n}}}\otimes M')=1$ if and only if $k\varphi_t = k(\varphi_{j_{n-1}}+\varphi_{j_{n}})\mod L$. Therefore,
\[
\rank \mathbb{V} = \rank \mathbb{V}_{0,n-1}(V_L,\bigotimes_{i=1}^{n-2}V_{L+k\varphi_{j_i}}\otimes V_{L+k\varphi_{j_{n-1}}+k\varphi_{j_{n}}})
\]
After finitely many such steps, one sees that $\rank \mathbb{V}\in\{0,1\}$, and that it equals $1$ if and only if $k\sum_{i=1}^n\varphi_{j_i}\equiv 0\mod L$. For non-triviality, it suffices to find a curve $C$ on $\overline{\mathrm{M}}_{0,n}$ so that $c_1\mathbb{V}\cdot C\ne0$. We consider $C$ to be an $F$-curve $F_{I_1,\ldots,I_4}$ given by a partition $I_1\sqcup\cdots\sqcup I_4$ of $\{1,\ldots,n\}$. We have
\[
c_1\mathbb{V} = k \cdot c_1\mathbb{V}_{0,n}(L_{\hat{\mathfrak{sl}}_{r+1}}(1,0),\bigotimes_{i=1}^nL_{\hat{\mathfrak{sl}}_{r+1}}(1,\varpi_{j_i}))
\]
by Proposition~\ref{prop: main prop} and Proposition~\ref{prop: tensor VB}. The conclusion then follows from Proposition 5.2 in \cite{Fak12}.
\end{proof}
\begin{example}[Case for $r=2$]
We let $k\geq1$ be an arbitrary integer. The lattice $L$ is simply $\zz\gb_1\oplus\zz\gb_2$ and the bi-linear pairing $(\cdot|\cdot):L\times L\to\zz$ is given by the matrix $B=\begin{pmatrix}2k&-k\\-k&2k\end{pmatrix}$. The set $S_2=\{V_L,V_{L+k\varphi_1},V_{L+k\varphi_2}\}$ embeds into a fusion subring of the fusion ring of $V_L$ with multiplication given by
\[
V_L\boxtimes V_{L+k\varphi_i}=V_{L+k\varphi_i},\quad V_{L+k\varphi_1}\boxtimes V_{L+k\varphi_1}=V_{L+k\varphi_2},\quad V_{L+k\varphi_2}\boxtimes V_{L+k\varphi_2}=V_{L+k\varphi_1},\quad  V_{L+k\varphi_1}\boxtimes V_{L+k\varphi_2}=V_L.
\]
Here $\varphi_1:=\frac{2\gb_1+\gb_2}{3k}$ and $\varphi_2:=\frac{\gb_1+2\gb_2}{3k}$.
This subring is isomorphic to the fusion ring $\mathcal{R}$ of the affine VOA $L_{\hat{\mathfrak{sl}}_{3}}(1,0)$ generated by $L_{\hat{\mathfrak{sl}}_{3}}(1,0),L_{\hat{\mathfrak{sl}}_{3}}(1,\varpi_1),L_{\hat{\mathfrak{sl}}_{3}}(1,\varpi_2)$.
Under the isomorphism, $L_{\hat{\mathfrak{sl}}_{3}}(1,0)\mapsto V_L$ and $L_{\hat{\mathfrak{sl}}_{3}}(1,\varpi_i)\mapsto V_{L+k\varphi_i}$ for $i\in\{1,2\}$. The proportionality is $\eta=k$, since
\[
\cw(V_{L+k\varphi_i})/\cw(L_{\hat{\mathfrak{sl}}_{3}}(1,\varpi_i)) = \frac{k/3}{1/3} = k.
\]
Therefore, by Proposition~\ref{prop: tensor VB}, we have that \[c_1 \mathbb{V}_{0,n}(V_L, {V_{L+k\varphi_1}}^{\otimes i}\otimes{V_{L+k\varphi_2}}^{\otimes (n-i)}) = k\cdot c_1\mathbb{V}_{0,n}(L_{\hat{\mathfrak{sl}}_{3}}(1,0), {L_{\hat{\mathfrak{sl}}_{3}}(1,\varpi_1)}^{\otimes i}\otimes{L_{\hat{\mathfrak{sl}}_{3}}(1,\varpi_2)}^{\otimes (n-i)}).\]
 By Corollary~\ref{cor: main corollary},  $\mathbb{V}_{0,n}(V_L, V_{L+k\varphi_1}^{\otimes i}\otimes V_{L+k\varphi_2}^{\otimes (n-i)})
\cong {\left(\mathbb{V}_{0,n}(L_{\hat{\mathfrak{sl}}_{3}}(1,0), L_{\hat{\mathfrak{sl}}_{3}}(1,\varpi_1)^{\otimes i}\otimes{L_{\hat{\mathfrak{sl}}_{3}}(1,\varpi_2)}^{\otimes (n-i)})\right)}^{\otimes k}.$
Finally, the global generation of $\mathbb{V}_{0,n}(V_L, V_{L+k\varphi_1}^{\otimes i}\otimes V_{L+k\varphi_2}^{\otimes (n-i)})$ follows from Theorem~\ref{thm: Fakh}.
\end{example}
\subsection{Positivity of bundles arising from rank $1$ even lattice VOAs}
\label{subsec: positivity rank 1 lattice}
In this section, we study coinvariant bundles arising from rank $1$ lattice VOAs for arbitrary pairing $(\alpha | \alpha) = 2k$. All rank $1$ even lattice VOAs arise this way. The VOA $V_L$, with  $(\alpha | \alpha) = 2$, is isomorphic to $L_{\hat{\mathfrak{sl}}_2}(1,0)$ and therefore all corresponding vector bundles are globally generated. Hence, we restrict our attention to cases $k\geq2$. All results in this section apply to any rank $1$ even lattice VOA with pairing $(\alpha | \alpha) \geq4$. The set of simple $V_L$-modules is
\[
\{M^r:=V_{L+r\varphi}\mid 0\leq r\leq 2k-1, \varphi:=\frac{\al}{2k}\},\quad \text{with conformal weight } \cw(M^{r}) = \frac{\left(\min(r,2k-r)\right)^2}{4k}.
\]

In Proposition~\ref{prop: subring of size 2p is not Fpositive} we describe large classes of subrings of the fusion ring of $V_L$ that are not $F$-positive. A fusion subring $\mathcal{S}$ is called \emph{$F$-positive} if every bundle of coinvariants on $\overline{\mathrm{M}}_{0,4}$ with marked modules in $\mathcal{S}$ has non-negative degree.
This result indicates that the fusion ring is not $F$-positive and a more in-depth analysis is needed to understand exactly which bundles on $\overline{\mathrm M}_{0,4}$ have non-negative degree. This is the content of Lemma~\ref{lem: 4 points rank 1 arb k}. In particular, the lemma identifies two types of $4$-tuples that exhibit positivity; which we will call here as type (iii) and type (iv) tuples (matching with the order in which they appear in the statement of the Lemma). Of the two, the type (iii) case is special as follows: for any $n$-tuple $(a_1,\ldots,a_n)$ with $\sum_{i=1}^n a_i\in\{2k,6k\}$, the intersection numbers of $c_1\mathbb{V}_{0,n}(V_L, \otimes_{i=1}^n M^{a_i})$ with $F$-curves can be described as positive linear combinations of the degrees of bundles $c_1\mathbb{V}_{0,4}(V_L, \otimes_{i=1}^n M^{b_i})$ satisfying $\sum b_i\in\{2k,6k\}$. This property indicates that semi-ampleness of type (iii) bundles may be extendable from $n=4$ to arbitrary $n$, which we prove in Lemma~\ref{lem: semiample 2k 6k}. A similar property does not hold for type (iv) tuples and we do not have any analogue of Lemma~\ref{lem: semiample 2k 6k} in the type (iv) case.



\begin{proposition}
\label{prop: subring of size 2p is not Fpositive}
Let $k\geq2$ be an integer. Consider any integer $p\geq2$ so that $p$ divides $k$. Then the subring $\zz[\zz/2p\zz]$ of the fusion ring $\zz[\zz/2k\zz]$ is not $F$-positive.
\end{proposition}
\begin{proof}
We can write the group $\zz/2p\zz$ as a set of modules by
\(
\{ M^{r\frac{k}{p}}\mid 0\leq r\leq 2p-1 \}.
\)
Here $M^{r\frac{k}{p}}:=V_{\zz\al+\frac{rk}{p}\frac{\al}{2k}}$. The conformal weight of such modules is given by
\[
\cw(M^{r\frac{k}{p}}) = \begin{cases}
\frac{kr^2}{4p^2} &\text{ if } 0\leq r\leq p\\
\frac{k(2p-r)^2}{4p^2} &\text{ if } p<r< 2p,
\end{cases}
\]
since $r\frac{k}{p}\leq 2k-r\frac{k}{p}$ if and only if $r\leq p$. Consider integers $i_1,\ldots,i_{2p-2}\geq0$ so that $\sum_{j=1}^{2p-2}i_j\leq 4$. The bundle
\[
\mathbb{V}:=\mathbb{V}_{0,4}(V_L,\bigotimes_{j=1}^{2p-2}(M^{j\frac{k}{p}})^{\otimes i_j}\otimes (M^{(2p-1)\frac{k}{p}})^{\otimes(4-\sum_{j=1}^{2p-2}i_j)})
\]
has rank $\leq1$ and it is equal to $1$ if and only if
$\sum_{j=1}^{2p-2}j\frac{k}{p}i_j+(2p-1)\frac{k}{p}(4-\sum_{j=1}^{2p-2}i_j)\equiv_{2k}0,$
that is, if and only if
\(
\sum_{j=1}^{2p-2}(j+1)i_j\equiv_{2p}4.
\)
If $p\geq 3$, let $i_1=i_{p-1}=2$ and all other $i_t=0$; if $p=2$, let $i_1=4$. In both cases the rank is $1$ and
\(
\deg c_1\mathbb{V}
= -\frac{k}{p}.
\)
\end{proof}
\begin{corollary}
Let $k\geq2$ be an integer. The fusion ring of $V_L$ is not $F$-positive.
\end{corollary}
\begin{proof}
It follows from Proposition~\ref{prop: subring of size 2p is not Fpositive} since every subring of an $F$-positive ring is also $F$-positive.
\end{proof}
\begin{example}
For $k=3$, the simple $V_L$-modules are of the form $V_{L+t\al/2k}$, where $0\leq t\leq 5$ are integers. Note that the subring $\mathcal{S}:=\cyc{V_L,V_{L+\al/2}}$ determines globally generated line bundles by Theorem~\ref{thm:main1} and Proposition~\ref{prop:rankonelattice}. By propagation of Vacua, we see that all such bundles are symmetric, that is, they are of the form $\mathbb{V}:=\mathbb{V}_{0,n}(V_L;{V^{\otimes n}_{L+\al/2}})$. A quick calculation shows that $\mathbb{V}$ is non-trivial if and only if $n$ is even. Another proper subring $\mathcal{T}$ of the fusion ring is generated by the modules $V_L, V_{L+\al/3}, V_{L+2\al/3}$, with relations $V_{L+a\al/3}\boxtimes V_{L+b\al/3}=V_{L+c\al/3}$, where $c\equiv a+b\mod 3$. The conformal weights of these modules are $0, \frac{1}{3}, \frac{1}{3}$, respectively, which coincide with the conformal weights of the irreducible modules of $L_{\hat{\mathfrak{sl}}_3}(1,0)$. Hence the map $L_{\hat{\mathfrak{sl}}_3}(1,\varpi_i)\mapsto V_{L+i\al/3}$ ($i=0,1,2$, with $\varpi_0:=0$) defines a proportional pairing between the fusion ring of $L_{\hat{\mathfrak{sl}}_3}(1,0)$ and $\mathcal{T}$ with proportionality $1$, and by Proposition~\ref{prop: tensor VB} every bundle $\mathbb{V}_{0,n}(V_L,V_{L+\al/3}^{\otimes i}\otimes V_{L+2\al/3}^{\otimes j})$, with $i+j=n$, is isomorphic to the corresponding level-one $\mathfrak{sl}_3$ bundle $\mathbb{V}_{0,n}(L_{\hat{\mathfrak{sl}}_3}(1,0),L_{\hat{\mathfrak{sl}}_3}(1,\varpi_1)^{\otimes i}\otimes L_{\hat{\mathfrak{sl}}_3}(1,\varpi_2)^{\otimes j})$; in particular it is globally generated, and it is nonzero if and only if $i\equiv j\mod 3$. For instance, for $n=4$ and $i=j=2$ one gets a line bundle of degree $1$ on $\overline{\mathrm{M}}_{0,4}$.
\end{example}



\begin{lemma}
\label{lem: 4 points rank 1 arb k}
Let $0<a_1,\ldots,a_4<2k$ be any integers and let $\mathbb{L}:=\mathbb{V}_{0,4}(V_L,M^{a_1}\otimes\cdots\otimes M^{a_4})$ be the corresponding coinvariant vector bundle.
\begin{enumerate}[(i)]
\item The rank of $\mathbb{L}$ is either $0$ or $1$.
\item If $\sum a_i\not\equiv 0\mod 2k$, then its first Chern class $c_1\mathbb{L}$ is zero.
\item If $\sum a_i\in \{2k,6k\}$, then $\deg c_1\mathbb{L}\leq 0$. Moreover, $\deg c_1\mathbb{L}\ne0$ if and only if $a_i<k$ for each $i$ (when $\sum a_i=2k$), respectively $a_i>k$ for each $i$ (when $\sum a_i=6k$).
\item If $\sum a_i=4k$, then $\deg c_1\mathbb{L}\geq0$.
\end{enumerate}
\end{lemma}
\begin{proof}
Since the fusion rules are of simple-current type, statements (i) and (ii) are easy to conclude. For (iii), note that the dual to $M^a$ is $M^{2k-a}$ and therefore, if $a_1+\cdots+a_4=6k$, then $a'_1+\cdots+a'_4=2k$ and since $\mathbb{V}_{0,4}(V_L,\bigotimes_{i=1}^4 M^{a_i})\cong \mathbb{V}_{0,4}(V_L,\bigotimes_{i=1}^4 {(M^{a_i})}')\cong \mathbb{V}_{0,4}(V_L,\bigotimes_{i=1}^4 M^{a'_i})$, we can restrict to $\sum a_i=2k$.


Assume that $0<a_1\leq a_2\leq a_3\leq a_4<k$ and $\sum a_i=2k$. For each integer $i=2,3,4$, let $b_i:=\min(a_1+a_i,2k-a_1-a_i)=\min(a_1+a_i,\ a_l+a_m)$, where $\{l,m\}=\{2,3,4\}\setminus\{i\}$, so that $\cw(\overline{a_1+a_i})=b_i^2/4k$. Then, $4k\deg\mathbb{V}_{0,4}(V_L,\bigotimes_{i=1}^4 M^{a_i})$ is equal to
\begin{align*}
\sum_{i=1}^4 a_i^2-\sum_{j=2}^4b_j^2
=\sum_{i=1}^4a_i^2-(a_1+a_2)^2-(a_1+a_3)^2-b_4^2
=\begin{cases}
4k(a_4-k) & \text{ if } (a_2+a_3\leq a_1+a_4),\\
-4ka_1 &\text{ if } (a_2+a_3\geq a_1+a_4).
\end{cases}
\end{align*}

To complete the proof of (iii), we need to prove that $\max(a_1,\ldots,a_4)\geq k$ implies that the degree of  $\mathbb{V}_{0,4}(V_L,\bigotimes_{i=1}^4 M^{a_i})$ is zero. Now assume instead that $a_1\geq a_2\geq a_3\geq a_4$ and $a_1\geq k$. Then, for each $i=2,3,4$, we have $a_i<k$ and $a_1+a_i>k$. Then, the degree is
\[
\frac{1}{4k}\left((2k-a_1)^2+\sum_{q=2}^4(a_q^2-(2k-a_1-a_q)^2)\right)=\frac{1}{4k}(2k-a_1)(2k-\sum a_i) = 0.
\]
This completes the proof for (iii).

For (iv), let $0<a\leq b\leq c\leq d<2k$ integers with $a+b+c+d=4k$. Let $m(x):=x^2-4k(x-k)_+$, where $(x)_+:=\max(x,0)$. Then, the conformal weight of $M^a$ is $\cw(a):=m(a)/4k$. Letting $\bar{x}:=x\mod 2k$, $\deg:=\deg\mathbb{V}_{0,4}(V_L,M^{a}\otimes M^b\otimes M^c\otimes M^d)$ is 
\begin{align*}
\deg=\cw(a)+\cw(b)+\cw(c)+\cw(d)-\cw(\overline{a+b})-\cw(\overline{a+c})-\cw(\overline{a+d}).
\end{align*}

\textbf{Case (A): $b+c\geq a+d$:} Note that $a+d\leq 2k$ and therefore, $a+b\leq a+c\leq a+d\leq 2k$ and we have
\[
\deg=-2a-d_+;\quad d_+:=(a-k)_++(b-k)_++(c-k)_++(d-k)_+-(a+b-k)_+-(a+c-k)_+-(a+d-k)_+.
\]
If $d<k$, then $a+b+c+d\leq 4d<4k$, a contradiction. Since $d\geq k$ and $a+d\leq 2k$, we have $a\leq k$. Similarly, if $c<k$, then $b+c\leq 2c<2k$, a contradiction. So, $a\leq k$ and $c,d\geq k$. Therefore,
\[
d_+=0+(b-k)_++(c-k)+(d-k)-(a+b-k)_+-(a+c-k)-(a+d-k)
=-2a+(b-k)_+-(a+b-k)_+.
\]
Putting it all together, we have
\[
\deg=(a+b-k)_+-(b-k)_+=\begin{cases}
a &\text{ if } b\geq k\\
(a+b-k)_+ &\text{ if } b<k.
\end{cases}
\]

\textbf{Case (B): $a+d\geq b+c$:} Since $a\leq k$, we have $2k\leq a+d<3k$. Also, $a+b\leq b+c\leq 2k$ and $a+c\leq b+c\leq 2k$. Therefore,
\[
\deg = \cw(a)+\cw(b)+\cw(c)+\cw(d)-\cw(a+b)-\cw(a+c)-\cw(a+d-2k) = d_{\mathrm{quad}}-d_+,
\]
\[
4kd_{\mathrm{quad}}:=a^2+b^2+c^2+d^2-(a+b)^2-(a+c)^2-(a+d-2k)^2=4k(d-a-k)
\]
and since $b>k$ implies $b+c\geq 2b>2k$, a contradiction, we have
\begin{align*}
d_+&:=(a-k)_++(b-k)_++(c-k)_++(d-k)_+-(a+b-k)_+-(a+c-k)_+-(a+d-3k)_+\\
&=0+0+(c-k)_++(d-k)-(a+b-k)_+-(a+c-k)_+-0.
\end{align*}
We claim that $a+c\geq k$. Indeed, the two conditions (a) $a+b,a+c<k$ and (b) $a+b\geq k$, $a+c<k$ are impossible: for (a), $a+b<k$ implies $c+d>3k$ and $a+c<k$ implies $c<k$; hence condition (a) would imply $d>2k$, a contradiction. For (b), $a+b\geq k$ implies $a+c\geq a+b\geq k$, which contradicts the other condition $a+c<k$. Therefore,
\[
 \deg = (d-a-k)-d_+=(a+b-k)_+ +(a+c-k)-a-(c-k)_+=(a+b-k)_+-(k-c)_+ .
\]
If $c\geq k$, this gives $\deg=(a+b-k)_+\geq 0$. If $c<k$, then $a+b\geq k$ (otherwise $d=4k-(a+b)-c>2k$), and $\deg=a+b+c-2k=2k-d\geq 0$. In summary,
\[
 \deg
=\begin{cases}
(a+b-k)_+ &\text{ if } c\geq k,\\
2k-d&\text{ if } c< k.
\end{cases}
\]
This completes the proof for (iv).
\end{proof}
This lemma gives all that we need for the study of $F$-nef divisors.

\begin{lemma}\label{lem: semiample 2k 6k}
For any positive integer $k\geq1$, let $0<a_1,\ldots,a_n<2k$ be $n$ integers so that $\sum a_i$ is equal to $2k$ or $2k(n-1)$. The conformal blocks divisor $c_1\mathbb{V}_{0,n}(V_L,\bigotimes_{i=1}^n M^{a_i})^{\vee}$ is semi-ample.
\end{lemma}
\begin{proof}
This follows from the degree calculation in Lemma~\ref{lem: 4 points rank 1 arb k}(iii) and comparing it with the degree calculation in \cite[Lemma 5.1]{Fak12}. If $\sum a_i=2k$ and every $a_i<k$, the divisor $-c_1\mathbb{V}_{0,n}(V_L,\otimes_{i=1}^n M^{a_i})$ is equal to the semi-ample divisor $c_1\mathbb{V}_{0,n}(L_{\hat{\mathfrak{sl}}_{k}}(1,0),\otimes_{i=1}^n L_{\hat{\mathfrak{sl}}_{k}}(1,a_i))$ by comparing their intersections with $F$-curves, and otherwise it is $0$. The case $2k(n-1)$ reduces to $2k$ by $a_i\mapsto 2k-a_i$. 
\end{proof}

\subsection{Classification of positive divisors in rank $1$ and $k=2$}
\label{subsec: classification rank 1 lattice k is 2}
In this subsection, we classify all positive bundles for $k=2$. Note that fixing $k$ fixes the lattice VOA we are working with and we will denote it by $V_L$. Let us denote the coinvariant bundles $\mathbb{V}_{0,n}(V_L,\otimes_{i=1}^n M^{a_i})$ by $\mathbb{L}_t$, where $t=(a_1,\ldots,a_n)$. One can prove that these bundles have rank $1$. This fact implies that the intersection of such a bundle $\mathbb{L}_t$ with an $F$-curve in $\overline{\mathrm M}_{0,n}$ is equal to the degree of a certain bundle on $\overline{\mathrm M}_{0,4}$ built from $V_L$-representations. Therefore, if we understand the degree of all bundles on $\overline{\mathrm M}_{0,4}$ built from $V_L$-representations (cf. Lemma~\ref{lem: k=2 4pts lemma}), then we can describe the intersection number $c_1\mathbb{L}_t\cdot C$ for any $F$-curve $C\s\overline{\mathrm M}_{0,n}$ (cf. Lemma~\ref{lem: k2 Fcurves}). Using these calculations, we classify all positive divisors arising from $V_L$-representations in Proposition~\ref{prop: classification in k2 rank1}. In particular, we describe all divisors, either coming from the coinvariant bundles $\mathbb{L}_t$ or from their duals $\mathbb{L}^{\vee}_t$, the conformal blocks bundles, that are positive. We find that all positive divisors are semi-ample and lie in the cone $\mathcal{C}\mathcal{B}(A,1)$ generated by type $A$, level $1$ conformal block divisors \cite{NoahAngela}.

\begin{lemma}
\label{lem: k=2 4pts lemma}
Let $0<a_1\leq\ldots\leq a_4<4$ be any four integers so that $\sum a_i=0\mod 4$. Let $\mathbb{L}(a_1,\ldots,a_4):=\mathbb{V}_{0,4}(V_L;\bigotimes_{i=1}^4M^{a_i})$.
\begin{enumerate}[(i)]
\item $\mathbb{L}(1,1,1,1)^{\vee}\cong \mathbb{L}(3,3,3,3)^{\vee}\cong\mathbb{V}_{0,4}(L_{\hat{\mathfrak{sl}}_2}(1,0), L_{\hat{\mathfrak{sl}}_2}(1,\varpi_1)^{\otimes 4})$; 
their degree is $-1$. 
\item If $\sum a_i=8$, then $(a_1,\ldots,a_4)$ is one of $(1,1,3,3)$, $(1,2,2,3)$, $(2,2,2,2)$, and
\[
\deg\mathbb{L}(1,1,3,3)=0,\qquad \deg\mathbb{L}(1,2,2,3)=1,\qquad \deg\mathbb{L}(2,2,2,2)=2.
\]
\end{enumerate}
\end{lemma}
\begin{proof}
It follows from Lemma~\ref{lem: 4 points rank 1 arb k}, or by a direct computation.
\end{proof}

\begin{proposition}
\label{prop: classification in k2 rank1}
All $F$-positive line bundles arising from representations of the lattice VOA $V_{\mathbb{Z}\al}$ with $(\al \mid \al)=4$ are semi-ample. Moreover, the first Chern class of each such bundle is proportional to the first Chern class of a line bundle coming from an affine VOA.
\end{proposition}

\begin{proof}
Let $t=(a_1,\ldots,a_n)\in \{1,2,3\}^{\oplus n}$, and write $\mathbb{L}_t$ for the corresponding coinvariant bundle. In what follows, $[\varpi_i]$ denotes the irreducible module $L_{\hat{\mathfrak{sl}}_m}(1,\varpi_i)$ of the level-one affine VOA $L_{\hat{\mathfrak{sl}}_m}(1,0)$. Lemma~\ref{lem: k2 Fcurves} below determines the intersections of $\mathbb{L}_t$ with all $F$-curves in $\overline{\mathrm{M}}_{0,n}$. We prove semi-ampleness of every $F$-positive bundle by comparing these intersection numbers with those of suitable affine VOA bundles.

The key point is that, in each $F$-positive case, the divisor $c_1(\mathbb{L}_t)$ has the same intersection numbers with all $F$-curves as a known globally generated divisor coming from an affine VOA, up to a positive scalar. Since $F$-curves span $N_1(\overline{\mathrm{M}}_{0,n})_{\mathbb{R}}$, it follows that these divisors lie on the same ray. Consequently, $c_1(\mathbb{L}_t)$ is semi-ample. More precisely, the $F$-positive cases are as follows.
\begin{itemize}
\item The divisor $c_1\mathbb{L}_{(1,3,2,\ldots,2)}$ lies on the ray spanned by
$c_1\mathbb{V}_{0,n}(L_{\hat{\mathfrak{sl}}_4}(1,0),[\varpi_1]\otimes[\varpi_3]\otimes[\varpi_2]^{\otimes(n-2)})$.

\item The divisors $c_1\mathbb{L}_{(1,1,2,\ldots,2)}$ and $c_1\mathbb{L}_{(3,3,2,\ldots,2)}$ lie on the ray spanned by
$c_1\mathbb{V}_{0,n}(L_{\hat{\mathfrak{sl}}_4}(1,0),[\varpi_1]^{\otimes 2}\otimes[\varpi_2]^{\otimes(n-2)})$.

\item If $\sum_{i=1}^n a_i\in\{4,4(n-1)\}$, then $t\in\{(1,1,1,1),(3,3,3,3)\}$, and
$\mathbb{L}^{\vee}_{t}\cong\mathbb{V}_{0,4}(L_{\hat{\mathfrak{sl}}_2}(1,0),[\varpi_{1}]^{\otimes 4})$.

\item The bundle $\mathbb{L}_{(2,\ldots,2)}$ satisfies
$\mathbb{L}_{(2,\ldots,2)}\cong \mathbb{V}_{0,n}(L_{\hat{\mathfrak{sl}}_2}(1,0),[\varpi_1]^{\otimes n})^{\otimes 2}$,
and is therefore globally generated.

\end{itemize}
This proves the claim.
\end{proof}

\begin{lemma}\label{lem: k2 Fcurves}
Let $k=2$ and let $n\geq 4$. Fix integers $a_i$ with $1\leq a_i\leq 3$, and denote by
$\mathbb{L}_t:=\mathbb{V}_{0,n}(V_L,\bigotimes_{i=1}^n M^{a_i})$
the coinvariant bundle corresponding to the tuple $t=(a_1,\ldots,a_n)$.
Then:
\begin{enumerate}[(i)]
\item $\rank\mathbb{L}_t\in\{0,1\}$.

\item If $n$ is even, then $\mathbb{L}_{(1,3,2,\ldots,2)}$ is $F$-positive and $\mathbb{L}_{(2,\ldots,2)}$ is globally generated.

\item If $n$ is odd, then $\mathbb{L}_{(1,1,2,\ldots,2)}$ and $\mathbb{L}_{(3,3,2,\ldots,2)}$ are $F$-positive.

\item If $\sum_{i=1}^n a_i\in\{4,4(n-1)\}$, then $\mathbb{L}^{\vee}_{(a_1,\ldots,a_n)}$ is $F$-positive.

\item Let $t=(2^b,3^c)$ or $t=(2^b,1^c)$, and assume that $c\equiv_4 2$ together with either $c\geq 6$ and $b=1$, or $c=2$ and $b\geq 3$. Then $\mathbb{L}_t$ is $F$-positive.

\item For every other tuple $t$, both the coinvariant bundle $\mathbb{L}_t$ and the conformal blocks bundle $\mathbb{L}_t^{\vee}$ fail to be $F$-positive.
\end{enumerate}
\end{lemma}

\begin{proof}
For (i), the rank is $1$ exactly when $\sum_{i=1}^n a_i\equiv 0 \mod 4$, and is $0$ otherwise.
For (ii), when all $a_i=2$, the claim was proved in Section~\ref{subsec: typeA}. We therefore focus on the remaining cases.

\textbf{Case A: $t=(1,3,2,\ldots,2)$ in (ii).}
It is enough to consider the following two types of $F$-curves.

First,
$F_1:=F_{\{(1,2^a),(3,2^b),(2^c),(2^d)\}}$,
where $a+b+c+d=n-2$, with $a,b\geq 0$, $a\equiv b\pmod 2$, and $c,d\geq 1$ odd.
For such a curve,
$\mathbb{L}_t\cdot F_1=\deg \mathbb{L}_{(1,2,2,3)}=1$.
Second,
$F_2:=F_{\{(1,3,2^a),(2^b),(2^c),(2^d)\}}$,
where $a+b+c+d=n-2$, $a\geq 0$, $b,c,d\geq 1$, and $a,b,c,d$ are all odd.
For such a curve,
$\mathbb{L}_t\cdot F_2=\deg \mathbb{L}_{(2,2,2,2)}=2$.

Hence $\mathbb{L}_t$ intersects every $F$-curve non-negatively, so it is $F$-positive.

\textbf{Case B: $t=(1,1,2,\ldots,2)$ in (iii).}
Again it suffices to consider two types of $F$-curves.

First,
$F_1:=F_{\{(1,2^a),(1,2^b),(2^c),(2^d)\}}$,
where $a,b\geq 0$, $a+b$ is odd, and $c,d\geq 1$ are odd.
Then
$\mathbb{L}_t\cdot F_1=\deg \mathbb{L}_{(1,2,2,3)}=1$.
Second,
$F_2:=F_{\{(1,1,2^a),(2^b),(2^c),(2^d)\}}$,
where $a\geq 0$ is even and $b,c,d\geq 1$ are all odd.
Then
$\mathbb{L}_t\cdot F_2=\deg \mathbb{L}_{(2,2,2,2)}=2$.
Thus $\mathbb{L}_t$ is $F$-positive. The case $t=(3,3,2,\ldots,2)$ is identical by symmetry.

\textbf{Case C: item (iv).} This follows from Lemma~\ref{lem: semiample 2k 6k}.

\textbf{Case D: The remaining tuples that constitutes item (v).} Let $t=(1^a,2^b,3^c)$ be the tuples that has not appeared in the proof yet. In this last part, our task is to find two $F$-curves, let us call them $F_{\mathrm{positive}}$ and $F_{\mathrm{negative}}$ so that any choices of $t$ in consideration in item (v), $c_1\mathbb{L}_t\cdot F_{\mathrm{positive}}>0$ and $c_1\mathbb{L}_t\cdot F_{\mathrm{negative}}<0$. This suffices to conclude that neither of $\mathbb{L}_t$ nor its dual $\mathbb{L}_t^{\vee}$ is $F$-positive. We need to formulate exactly the set of $3$-tuples $(a,b,c)$ we need to consider:  
\begin{itemize}
    \item We must have $a,b,c\in\zz_{\geq0}$ and $a+b+c\geq4$ for clear reasons. 
    \item As $\mathrm{rank } \mathbb{L}_t \ne0$ if and only if $a+2b+3c=0\mod 4$, we require $a+2b+3c=0\mod 4$.
    \item As in the proof of Lemma~\ref{lem: 4 points rank 1 arb k}(iii), $\mathbb{L}_{(1^a,2^b,3^c)}\cong\mathbb{L}_{(1^c,2^b,3^a)}$, so we may assume $a\geq c$.
    \item Moreover, cases 
    \begin{itemize}
        \item $(a,b,c)=(0,b,0)$ and $(a,b,c)=(1,b,1)$ with evenly many $b$,
        \item $(a,b,c)=(2,b,0)$ with $b\geq1$ odd, and
        \item $(a,b,c)=(4,0,0)$
    \end{itemize} are already discussed in items (ii), (iii) and (iv), respectively. So, we ignore them here.
    \item  We also ignore $(a,b,c)=(2,0,2)$ as $\deg \mathbb{L}_{(1,1,3,3)}=0$.
\end{itemize}
 With these constraints in place, the set of tuples $(a,b,c)$ we left to consider and the choice of $F_{\mathrm{positive}}$ and $F_{\mathrm{negative}}$ for each such $(a,b,c)$ is listed in Table~\ref{tab:cases}. 

 \begin{table}[htbp]
\centering
\scriptsize
\renewcommand{\arraystretch}{1.8}
\setlength{\tabcolsep}{3pt}

\begin{tabular}{|c|p{6cm}|p{3cm}|p{3cm}|}
\hline
Case number & Case & $F_{\text{positive}}$ & $F_{\text{negative}}$ \\
\hline
E1 & $a\geq 3,\ b\geq 2$ & $F_{(1),(2)^2,(\mathrm{rest})}$ & $F_{(1)^3,(\mathrm{rest})}$ \\
\hline
E2 & $a\geq 3,\ b=1$ & $F_{(1),(1,1),(2),(\mathrm{rest})}$ & $F_{(1)^3,(\mathrm{rest})}$ \\
\hline
E3 & $a\geq 3,\ b=0,\ c\geq 2$ & $F_{(1),(1,1),(3,3),(\mathrm{rest})}$ & $F_{(1)^3,(\mathrm{rest})}$ \\
\hline
E4 & $a\geq 5,\ b=0,\ c\leq 1$ & $F_{(1),(1,1)^2,(\mathrm{rest})}$ & $F_{(1)^3,(\mathrm{rest})}$ \\
\hline
E5 & $a=c=2,\ b\geq 2$ & $F_{(1),(2)^2,(\mathrm{rest})}$ & $F_{(1)^2,(2,3),(2,3,\mathrm{rest})}$ \\
\hline
\end{tabular}
\caption{Choices of $F_{\mathrm{positive}}$ and $F_{\mathrm{negative}}$ curves for Case D.}
\label{tab:cases}
\end{table}
This completes the proof.
\end{proof}

\appendix
\section{Some preliminaries on VOAs}\label{sec:appendix}

\subsection{Basics of lattice VOAs}
Let $L$ be a rank $r\geq 1$ positive-definite even lattice with respect to $(\cdot|\cdot):L\times L\ra \Z$. Let $\h:=\C\otimes_{\Z} L$, extend $(\cdot|\cdot)$ to a nondegenerate $\C$-bilinear form $(\cdot|\cdot):\h\times \h\ra \C$, and let $M_{\hat{\h}}(1,0)$ be the level-one Heisenberg VOA associated to $\hat{\h}=\h\otimes \C[t,t^{-1}]\op \C\kappa$.

Let $\epsilon: L\times L\ra \<\pm 1\>$ be a $2$-cocycle of the abelian group $L$ such that $\epsilon(\al,\b)\epsilon(\b,\al)=(-1)^{(\al|\b)}$ for any $\al,\b\in L$. Write $\C^\ep[L]=\bigoplus_{\al\in L} \C e^{\al}$, where $e^{\al}$ is a formal symbol associated with $\al$ for each $\al\in L$. Define $V_L=M_{\hat{\h}}(1,0)\otimes_{\C} \C^\epsilon[L]$. Then
$$V_L=\spn\{ h_1(-n_1)\ds h_k(-n_k)e^\al: k\geq 0, \al\in L, \ h_1,\ds,h_k\in \h,\ n_{1}\geq \ds \geq n_k\geq 0 \},$$
where we omit the tensor sign $\o$ in the term $h_1(-n_1)\ds h_k(-n_k)\o e^\al$. The vertex operator $Y:V_L\ra \End(V_L)[\![z,z^{-1}]\!]$ has the following definition on the spanning elements of $V_L$:
\begin{align}
	&Y(h(-1)\vac,z):=h(z)=\sum_{n\in \Z} h(n) z^{-n-1}\quad \left(h(n)e^\al:=0,\ \forall n>0,\ h(0)e^\al:=(h|\al) e^\al\right), \label{latticeh}\\
	&Y(e^{\al},z):=E^{-}(-\al, z)E^{+}(-\al,z)e_{\al} z^{\al} \quad \left(z^\al (e^\b):=z^{(\al|\b)}e^\b,\ e_\al(e^\b):=\epsilon (\al,\b)e^{\al+\b}\right),\label{latticeealpha}\\
	&Y(h_1(-n_{1}-1)\dots h_k(-n_{k}-1)e^{\al},z):={\tiny\begin{matrix}\circ \\\circ\end{matrix}}(\partial_{z}^{(n_1)}h_1(z))\dots (\partial_{z}^{(n_k)}h_k(z))Y(e^{\al},z){\tiny\begin{matrix}\circ \\\circ\end{matrix}},\label{latticegeneral}
\end{align}
for any $k\geq 1$, $n_{1}\geq \ds\geq n_k\geq 0$, $h,h_1,\dots,h_k\in \h$, and $\al,\b\in L$, where the operators $E^{\pm}$, $\partial^{(n)}_{z}$, $e_\al$, and $z^\al$ in \eqref{latticeh}--\eqref{latticegeneral} are given as follows:
$$ E^{\pm}(-\al,z)=\exp\left(\sum_{n\in \Z_{\pm}}\frac{-\al(n)}{n}z^{-n}\right),\quad \partial_z^{(n)}=\frac{1}{n!}\frac{\partial^n}{\partial z^n},\quad e_\al(e^{\b})=\epsilon(\al,\b)e^{\al+\b},\quad z^\al(e^\b)=z^{(\al|\b)}e^{\b}.
$$

Let $L^\circ=\{\la\in \h: (\la|\al)\in \Z,\ \forall \al\in L\}$ be the dual lattice, and let $L^\circ/L=\bigsqcup_{i=1}^N (L+\la_i)$. Then, by \cite[Theorem 3.1]{dong93}, all the irreducible $V_L$-modules are given by $V_{L+\la_i}$, $1\leq i\leq N$. The module vertex operators $Y_M: V_L\ra \End(V_{L+\la_i})[\![z,z^{-1}]\!]$ are given by formulas similar to \eqref{latticeh}--\eqref{latticegeneral}.


\subsection{Basics of affine VOAs and simple-current extensions}
\label{subsubsec: GR and simple-current extension}
We first recall the definition of affine vertex operator algebras of positive integral level $k$ \cite{FZ92}.

Let $\g$ be a finite-dimensional semisimple Lie algebra with a Cartan subalgebra $\h$, and let $\Delta$ be the root system associated with $\h$, with root lattice $Q\subset \h^\ast$. Normalize the invariant bilinear form on $\g$ so that $(\theta|\theta)=2$, where $\theta$ is the highest root of $\Delta$. Let $\hat{\g}=\g\otimes \C[t,t^{-1}]\op \C\kappa$ be its affinization, with Lie bracket given by
$$
[\kappa,\hat{\g}]=0,\quad [a(m),b(n)]=[a,b](m+n)+m\delta_{m+n,0} (a|b) \kappa,\quad a,b\in \g,\ m,n\in \Z.
$$
Let $\C \vac$ be a $\hat{\g}_{\geq 0}$-module with $\kappa.\vac=k\vac$ and $\g\o \C[t].\vac=0$. The Weyl vacuum module
$$
V^{k}(\g)=\scrU(\hat{\g})\o _{\scrU(\hat{\g}_{\geq 0})} \C \vac
$$
is a VOA, with $\vac =1\o \vac$ and $\om_{\mathrm{aff}}=\frac{1}{2(h^\vee+k)}\sum_{i=1}^{\dim \g}u^i(-1)u_i(-1)\vac$, called the {\em vacuum module vertex operator algebra of level $k$}, where $h^\vee$ is the dual Coxeter number of $\Delta$, and $\{ u^i\}$ and $\{ u_i\}$ are dual orthonormal bases of $\g$.

The module $V^{k}(\g)$ has a unique maximal $\hat{\g}$-submodule $\mathcal{J}^k(\g)=\scrU(\hat{\g})x_{\theta}(-1)^{k+1}\vac$, where $x_\theta\in \g_\theta$. The quotient module $L_{k}(\g)=V^{k}(\g)/\mathcal{J}^k(\g)$ is also a VOA, called the {\em affine VOA of level $k$}. It is a strongly rational VOA \cite{DLM97}.

Given $\la\in P$, one can consider the Weyl module $V^k(\g,\lambda)=\scrU(\hat{\g})\o _{\scrU(\hat{\g}_{\geq 0})} L(\la)$, where $L(\la)$ is the irreducible highest-weight $\g$-module associated with the weight $\la$. Let $L_k(\g,\la)$ be the unique irreducible quotient $\hat{\g}$-module of $V^k(\g,\lambda)$. Then
\[
L_k(\g,\lambda),\quad \la\in P^k_+=\{\la\in P_+: (\la|\theta)\leq k\},
\]
are all the irreducible modules over the affine VOA $L_k(\g)$ up to isomorphism. Let $\rho=\frac{1}{2}\sum_{\al\in \Delta_+}\al$. The conformal weight of $L_k(\g,\la)$ is
\begin{equation}\label{eq:affineconformalwt}
	\cw(L_k(\g,\la))=\frac{(\la+2\rho|\la)}{2(h^\vee+k)},\quad \la\in P^k_+.
\end{equation}

Let $V$ be a strongly rational VOA, with irreducible modules $\mathscr{W}=\{W^0,\ds, W^N \}$ up to isomorphism. Recall that an irreducible $V$-module $S\in \mathscr{W}$ is called a {\em simple current} if the fusion tensor product of $S$ with any other irreducible $V$-module is again an irreducible $V$-module:
\[
S\boxtimes W^i= W^{n_i},\quad \mathrm{for}\quad 0\leq i\leq N.
\]

Recall the following theorem of Dong--Li--Mason:
\begin{lemma}[{\cite[Theorem 3.21]{DLM97}}]\label{thm:DLM}
	Let $V=L_k(\g)$ be the affine VOA at level $k\in \Z_{>0}$. Let $M^1,\ds , M^s$ be all the irreducible simple-current $V$-modules. Then $V_e=V\op M^1\op\ds \op M^s$ is a simple abelian intertwining operator algebra. In particular, if the conformal weights of $M^i$ are all integers, then $V_e$ is a simple VOA.
\end{lemma}

\subsubsection{The simple-current extension of $L_8(\sl_2)$}
Consider the type-$A_1$ level-$8$ affine VOA $V=L_8(\sl_2)$. It has $9$ irreducible modules:
$$L_8(i)=L_8\left(\sl_2,i\al/2\right),\quad 0\leq i\leq 8,$$
where $\al=\theta$ is the unique positive root, with $(\al|\al)=2$. The fusion tensor product is given by
\begin{equation}\label{eq:sl2fusion}
	L_8(i)\boxtimes L_8(j)=\bigoplus_{ \substack{|i-j|\leq r\leq \min(i+j,16-i-j),\\ r\equiv i+j\ (\mathrm{mod}\ 2)}} L_8(r).
\end{equation}
In particular, $J=L_8(8)$ is the unique simple-current irreducible $V$-module, with
\begin{equation}\label{L8tensor}
	L_8(8)\boxtimes L_8(i)= L_8(8-i),\quad 0\leq i\leq 8.
\end{equation}
By \eqref{eq:affineconformalwt}, the conformal weight of $J=L_8(8)$ is
$
\frac{(8\al/2+\al|8\al/2)}{2(h^\vee+k)}=\frac{(5\al|4\al)}{2(2+8)}=2\in \Z.
$
It follows from Lemma~\ref{thm:DLM} that the following direct sum space is a VOA:
\begin{equation}
	V_e=V\op J=L_8(\sl_2)\op L_8(8),
\end{equation}
which is the simple-current extension of the affine VOA $V=L_8(\sl_2)$.

Let $I\in I_V\fusion{J}{J}{V}$ be the unique nonzero intertwining operator such that
$
v^8_{(3)}u^8=\vac\in V_0,
$
where $u^8,v^8$ denote lowest-weight vectors of $J$ and $v^8_{(n)}$ denotes the $n$-th coefficient of $I(v^8,z)$.
Denote the vertex operator of $V$ by $Y_V$ and the $V$-module vertex operator for $J$ by $Y_J$. Then the vertex operator $Y_{V_e}$ is given by
\begin{equation}\label{eq:YV}
	Y_{V_e}((a,u),z)(b,v)=\left(Y_V(a,z)b+I(u,z)v, Y_J(a,z)v+e^{zL(-1)}Y_J(b,-z)u \right),
\end{equation}
where $(a,u), (b,v)\in V_e=V\op J$.

\begin{lemma}[{\cite[Corollary 5.8]{Li97}}]\label{prop:irrmodV}
	The simple-current extension VOA $V_e=L_8(\sl_2)\op L_8(8)$ has the following irreducible modules and associated conformal weights, up to isomorphism:
	\begin{align*}
		&V_e=L_8(\sl_2)\op L_8(8),&& h_{V_e}=0,\\
		&X=L_8(2)\op L_8(6),&& h_{X}=1/5,\\
		&W^+=L_8(4),\ Y_{W^+}=Y_{L_8(4)}+\Y,&& h_{W^+}=3/5,\\
		&W^-= L_8(4),\ Y_{W^-}=Y_{L_8(4)}-\Y,&& h_{W^-}=3/5,
	\end{align*}
	where $\Y\in I_V\fusion{J}{L_8(4)}{L_8(4)}$ is a normalized intertwining operator. In other words, $\mathscr{W}=\{V_e,X,W^+,W^-\}$. 
\end{lemma}




\subsubsection{Fusion rules}
The fusion rules among irreducible $V_e=L_8(\sl_2)\op L_8(8)$-modules can be calculated using the known fusion rules among irreducible $V=L_8(\sl_2)$-modules. The fusion rules among irreducible $V_e$-modules play a critical role in the global generation of $V_e$-coinvariant bundles. Since they cannot be found in the literature, we present a detailed calculation.

\begin{proposition}\label{prop:fusionsimplecurrent}
	The fusion tensor products among the irreducible $V_e$-modules are
	\begin{align*}
		&	X\boxtimes_{V_e} X \cong V_e\oplus X\oplus W^+\oplus W^-,&& 	W^+\boxtimes_{V_e} W^- \cong X,\\
		&	X\boxtimes_{V_e} W^+ \cong X\oplus W^-, && X\boxtimes_{V_e} W^- \cong X\oplus W^+,\\
		&	W^+\boxtimes_{V_e} W^+ \cong V_e\oplus W^+, && 	W^-\boxtimes_{V_e} W^- \cong V_e\oplus W^-.
	\end{align*}
\end{proposition}

\begin{proof}
	Since both $V=L_8(\sl_2)$ and $V_e$ are strongly rational VOAs, and $V\hookrightarrow V_e$ is a conformal embedding of VOAs, there is a categorical induced module functor:
	\begin{equation}
		\mathrm{Ind}^{V_e}_V(-): \mathsf{Adm}(V)\ra \mathsf{Adm}(V_e),\quad \mathrm{Ind}^{V_e}_V(M):=V_e\bt_VM,
	\end{equation}
	where the fusion tensor product is taken in the category $\mathsf{Adm}(V)$, and $V_e$ is viewed as a natural $V$-module via the embedding $V\hookrightarrow V_e$.
	Moreover, it is easy to show that the induced $V_e$-modules $\mathrm{Ind}^{V_e}_V(L_8(i))$ have the following identifications:
	\begin{align*}
		& \mathrm{Ind}^{V_e}_V(L_8(0))=\mathrm{Ind}^{V_e}_V(L_8(8))=V_e,\quad 	\mathrm{Ind}^{V_e}_V(L_8(4))=W^+\op W^-,\\
		&\mathrm{Ind}^{V_e}_V(L_8(2))=\mathrm{Ind}^{V_e}_V(L_8(6))=X.
	\end{align*}

	It follows from \eqref{eq:sl2fusion} that the $V=L_8(\mathfrak{sl}_2)$-modules have the following fusion tensor products:
	\begin{align*}
		L_8(2)\boxtimes_{V} L_8(2) &\cong L_8(0)\oplus L_8(2)\oplus L_8(4),\\
		L_8(2)\boxtimes_{V} L_8(4) &\cong L_8(2)\oplus L_8(4)\oplus L_8(6),\\
		L_8(4)\boxtimes_{V} L_8(4) &\cong L_8(0)\oplus L_8(2)\oplus L_8(4)\oplus L_8(6)\oplus L_8(8).
	\end{align*}
	Note that the induced module functor satisfies
	\begin{align*}
		\mathrm{Ind}^{V_e}_V(M\boxtimes_{V} N)&={V_e}\bt_V M\bt_V N\cong (V_e\bt_V M)\bt_{V_e} (V_e\bt_V N)\\
		&\cong \mathrm{Ind}^{V_e}_V(M)\boxtimes_{V_e} \mathrm{Ind}^{V_e}_V(N).
	\end{align*}
	Applying $\mathrm{Ind}^{V_e}_V(-)$ to the tensor products of $V$-modules and using Lemma~\ref{prop:irrmodV}, we have
	\begin{align}
		X\boxtimes_{V_e} X &\cong V_e\oplus X\oplus W^+\oplus W^-, \label{eq:XX}\\
		X\boxtimes_{V_e} (W^+\oplus W^-) &\cong 2X\oplus W^+\oplus W^-, \label{eq:XWsum}\\
		(W^+\oplus W^-)\boxtimes_{V_e} (W^+\oplus W^-) &\cong 2V_e\oplus 2X\oplus W^+\oplus W^-.
		\label{eq:WWsum}
	\end{align}
	Now we invoke the quantum dimensions of VOA-modules \cite{MS89,DJX13}. Given an irreducible module $M$ over the strongly rational VOA $V$, recall that
	\[
	\qdim_V(M)=\lim_{q\ra 1^-}\frac{\Ch_qM}{\Ch_q V}=\lim_{q\ra 1^-} \frac{\tr|_M q^{L(0)-\frac{c}{24}}}{\tr|_Vq^{L(0)-\frac{c}{24}}},
	\]
	where $c$ is the central charge of $V$. In particular, $\qdim_V(V)=1$.

	For our simple-current extension $V_e=V\op J$, both $V$ and $V_e$ are strongly rational. They also share the same central charge. By \cite[Proposition 4.9]{DJX13}, for any irreducible $V$-modules $M,N$, we have
	$$\qdim_V(M\bt_V N)=\qdim_V(M)\cdot \qdim_V(N).$$
	Now, if we view $V_e$ as an adjoint $V$-module, then
	\begin{align*}
		\qdim_{V_e}(\mathrm{Ind}^{V_e}_V(L_8(i)))&=\qdim_{V_e}(V_{e}\bt_V L_8(i))=\lim_{q\ra 1^-}\frac{\Ch_q (V_e\bt_V L_8(i))/\Ch_qV }{\Ch_q V_e/\Ch_q V}\\
		&=\frac{\qdim_V(V_e\bt_V L_8(i))}{\qdim_V(V_e)}=\frac{\qdim_V (V_e)\cdot \qdim_V(L_8(i))}{\qdim_V(V_e)}\\
		&=\qdim_V(L_8(i)).
	\end{align*}
	In particular, $\qdim_{V_e} (X)=\qdim_V(L_8(2))$ and $\qdim_{V_e}(W^+)+\qdim_{V_e}(W^-)=\qdim_V(L_8(4))$. Moreover, $\qdim_{V_e}(W^+)=\qdim_{V_e}(W^-)$, since $W^+$ and $W^-$ have the same graded dimension as $L_8(4)$.

	For $V=L_8(\mathfrak{sl}_2)$, using the formula in \cite[Example 4.18]{DJX13}, we have
	\[
	\qdim_V L_8(i)=\frac{e^{\frac{\pi \sqrt{-1}}{10}(i+1)}-e^{\frac{-\pi \sqrt{-1}}{10}(i+1)}}{e^{\frac{\pi \sqrt{-1}}{10}}-e^{\frac{-\pi \sqrt{-1}}{10}}}=\frac{\sin\!\left((i+1)\pi/10\right)}{\sin\!\left(\pi /10\right)}.
	\]
	It follows that
	$	\qdim_{V_e}(X)=\qdim_V L_8(2)=\frac{\sin(3\pi/10)}{\sin(\pi/10)}=\varphi^2,
	$	where \(\varphi=\frac{1+\sqrt5}{2}\), and
	$	\qdim_{V_e}(W^\pm)=\frac{1}{2}\qdim_V L_8(4)=\frac{1}{2}\cdot
	\frac{\sin(5\pi/10)}{\sin(\pi/10)}=\varphi.
	$
	In particular, we have
	\[
	\qdim_{V_e}(W^\pm\boxtimes_{V_e} W^\pm)=	\qdim_{V_e}(W^\pm)\cdot 	\qdim_{V_e}(W^\pm)=\varphi^2=\varphi+1.
	\]
	Note that $W^\pm\boxtimes_{V_e} W^\pm$ is a direct sum of the irreducible $V_e$-modules $V_e,X,W^\pm$. The only possible sums of these modules that have quantum dimension
	\(\varphi+1\) are
	\[
	V_e\oplus W^+,\quad V_e\oplus W^-,\quad X.
	\]
	Thus, $W^\pm\boxtimes_{V_e} W^\pm\in \{V_e\oplus W^+,V_e\oplus W^-,X\}$.
	For simplicity of notation, write
	$a=W^+ $ and $ b=W^-.$
	Then \eqref{eq:WWsum} becomes
	\begin{equation}\label{eq:4.14}
		a^2+b^2+2ab = 2V_e+2X+a+b,
	\end{equation}
	with $ab\in \{V_e\op a, V_e\op b,X \}$. The first two cases are impossible. Indeed, if \(ab= V_e\oplus a\), then
	\[
	a^2+b^2 +2V_e+2a= 2V_e+2X+a+b.
	\]
	However, there is only one copy of the irreducible $V_e$-module $a$ on the right-hand side, which is a contradiction. Similarly, $ab\neq V_e\op b$. Hence, $ab=X$, i.e.,
	\begin{equation}\label{eq:fusionW+withW-}
		W^+\bt_{V_e} W^-=W^-\bt_{V_e} W^+=X.
	\end{equation}
	Cancelling $2ab=2X$ from both sides of \eqref{eq:4.14}, we have
	\[
	a^2+b^2 \cong 2V_e+a+b.
	\]
	Since each of $a^2$ and $b^2$ has quantum dimension $\varphi+1$, it follows that either
	\[
	\begin{cases}
		a^2=V_e\op a\quad  \mathrm{and}\quad b^2=V_e\op b,\\
		\mathrm{or}\\
		a^2=V_e\op b\quad \mathrm{and}\quad b^2=V_e\op a.
	\end{cases}
	\]
	We claim that the second case cannot occur. Indeed, suppose otherwise. Then $W^+\bt_{V_e} W^+=V_e\op W^-$ and $W^-\bt_{V_e} W^-=V_e\op W^+$. It follows that
	\begin{align*}
		W^+\bt_{V_e} X&\cong(W^+\bt_{V_e} W^+)\bt_{V_e} W^- =(V_e\op W^-)\bt_{V_e} W^-=W^-\op  V_e\op W^+,\\
		W^-\bt_{V_e} X&\cong (W^-\bt_{V_e} W^-)\bt_{V_e} W^+ =(V_e\op W^+)\bt_{V_e} W^+=W^+\op  V_e\op W^-.
	\end{align*}
	Then $X\bt_{V_e}(W^+\op W^-)=2V_e\op 2W^+\op 2W^-$. This contradicts \eqref{eq:XWsum}. Therefore, we must have
	\begin{equation}\label{eq:fusionWpm}
		W^+\bt_{V_e} W^+=V_e\op W^+\quad \mathrm{and}\quad W^-\bt_{V_e} W^-=V_e\op W^-.
	\end{equation}
	Finally, using \eqref{eq:fusionWpm} and \eqref{eq:fusionW+withW-}, we have
	\[
	W^+\bt_{V_e}X\cong (W^+\bt_{V_e} W^+)\bt_{V_e} W^-=(V_e\op W^+)\bt_{V_e} W^-=W^-\op X.
	\]
	Similarly, $W^-\bt_{V_e} X\cong W^+\op X$. This finishes the proof of the claimed fusion tensor products, in view of \eqref{eq:XX}, \eqref{eq:fusionWpm}, and \eqref{eq:fusionW+withW-}.
\end{proof}
\subsection{Tensor product of VOAs}\label{sec:tensor}
Recall the notion of tensor product VOAs and their irreducible modules introduced by Frenkel-Huang-Lepowsky.
\begin{lemma}[{\cite[Section 2.5, Proposition 4.7.2]{FHL93}}]\label{lm:tensor}
	Let $(V^1,Y_{V^1},\om_1,\vac_1)$ and $(V^2,Y_{V^2},\om_2,\vac_2)$ be two VOAs. The tensor product space $V^1\o V^2$ has a VOA structure with
	\begin{equation}\label{eq:tensor}
		Y_{V^1\o V^2}(a^1\o a^2,z)(b^1\o b^2)=Y_{V^1}(a^1,z)b^1\o Y_{V^2}(a^2,z)b^2,\quad a^1,b^1\in V^1, a^2,b^2\in V^2,
	\end{equation}
	$\vac_{V^1\o V^2}=\vac_1\o \vac_2$ and $\om_{V^1\o V^2}=\om_1\o \vac_2+\vac_1\o \om_2$.

	Furthermore, the irreducible ordinary modules over $V^1\o V^2$ are of the form $M^1\o M^2$, where $M^i$ is an irreducible ordinary module over $V^i$ for $i=1,2$. In particular, if $V^1$ and $V^2$ are both simple VOAs, then $V^1\o V^2$ is a simple VOA.
\end{lemma}

Moreover, if $V^1$ and $V^2$ are both strongly rational, then $V^1\o V^2$ is also strongly rational \cite{DMZ94,DLM2}. In particular, let $\mathscr{W}^1=\{M^{11},\ds, M^{1m}\}$ and $\mathscr{W}^2=\{M^{21},\ds, M^{2n}\}$ be the sets of all irreducible modules over $V^1$ and $V^2$, respectively. Then
$$\mathscr{W}^1\bt \mathscr{W}^2:=\{M^{1i}\o M^{2j}: 1\leq i\leq m,\ 1\leq j\leq n\}$$
is the set of all irreducible modules over $V^1\o V^2$.

\bibliographystyle{alpha} 
\bibliography{ref3}

\end{document}